\documentclass[12pt,a4paper]{article}
\usepackage{amscd,amsmath,amssymb,amsfonts,times}
\usepackage{mathrsfs}
\usepackage[all]{xy}
\numberwithin{equation}{subsection}
    \newtheorem{thm}{Theorem}[subsection]
    \newtheorem{lem}[thm]{Lemma}
    
    \newtheorem{prop}[thm]{Proposition}
    \newtheorem{cor}[thm]{Corollary}
    \newtheorem{cond}[thm]{Condition}
    \newtheorem{claim}[thm]{Claim}
    
    \newtheorem{defn}[thm]{Definition}
    \newtheorem{exmp}[thm]{Example}
    \newtheorem{rem}[thm]{Remark}

\newcommand{\qed}
{\mbox{}\nolinebreak$\square$\medbreak\par}
\newenvironment{pf}{\par\smallskip\noindent\emph{Proof.}}{\hfill\qed\par\smallskip}
\newenvironment{pf*}[1]{\par\smallskip\noindent\emph{#1.}}{\hfill\qed\par\smallskip}
\begin{document}
\title{Syntomic cohomology and Beilinson's Tate conjecture for $K_2$}
\author{Masanori Asakura and Kanetomo Sato}
\date\empty
\maketitle
{\small
\tableofcontents
}\bigskip

\def\can{\text{\rm can}}
\def\cano{\text{\rm canonical}}
\def\cd{\text{\rm cd}}
\def\ch{\text{\rm ch}}
\def\Coker{\mathrm{Coker}}
\def\Cor{\text{\rm Cor}}
\def\crys{\mathrm{crys}}
\def\dec{\text{\rm dec}}
\def\dlog{d{\text{\rm log}}}
\def\dR{\text{\rm d\hspace{-0.2pt}R}}
\def\Eis{{\text{\rm E}}}
\def\et{\text{\rm \'et}}  % etale
\def\exp{\text{\rm exp}}
\def\Frac{\text{\rm Frac}}
\def\phami{\phantom{-}}
\def\id{\text{\rm id}}              % identity
\def\Image{\text{\rm Im}}        % image
\def\ind{\text{\rm ind}}
\def\Hom{\text{\rm Hom}}  
\def\ker{\text{\rm Ker}}          % kernel
\def\mf{\text{\rm mapping fiber of}}
\def\Pic{\text{\rm Pic}}
\def\CH{\text{\rm CH}}
\def\NS{\text{\rm NS}}
\def\End{\text{\rm End}}
\def\pr{\text{\rm pr}}
\def\Proj{\text{\rm Proj}}
\def\ord{\text{\rm ord}}
\def\qis{\text{\rm qis}}
\def\rank{\text{\rm rank}}
\def\reg{\text{\rm reg}}
\def\res{\text{\rm res}}
\def\Res{\text{\rm Res}}
\def\Spec{\text{\rm Spec}}
\def\syn{\text{\rm syn}}
\def\tors{\text{\rm tors}}
\def\tr{\text{\rm tr}}
\def\cont{\text{\rm cont}}
\def\zar{\text{\rm zar}}
\def\bA{{\mathbb A}}
\def\bC{{\mathbb C}}
\def\C{{\mathbb C}}
\def\G{{\mathbb G}}
\def\bE{{\mathbb E}}
\def\bF{{\mathbb F}}
\def\F{{\mathbb F}}
\def\bH{{\mathbb H}}
\def\bJ{{\mathbb J}}
\def\bN{{\mathbb N}}
\def\bP{{\mathbb P}}
\def\P{{\mathbb P}}
\def\bQ{{\mathbb Q}}
\def\Q{{\mathbb Q}}
\def\bR{{\mathbb R}}
\def\R{{\mathbb R}}
\def\bZ{{\mathbb Z}}
\def\Z{{\mathbb Z}}

\def\cC{{\mathscr C}}
\def\cD{{\mathscr D}}
\def\cE{{\mathscr E}}
\def\cK{{\mathscr K}}
\def\cO{{\mathscr O}}
\def\O{{\mathscr O}}
\def\cR{{\mathscr R}}
\def\cS{{\mathscr S}}
\def\cU{{\mathscr U}}
\def\cV{{\mathscr V}}
\def\cX{\hspace{-0.8pt}{\mathscr X}\hspace{-1.2pt}}

\def\fo{\mathfrak o}

\def\XR{\cX}
\def\XK{X}
\def\XZp{\cX}
\def\XQp{X}
\def\XbK{\ol X}

\def\UR{\cU}
\def\UK{U}
\def\UZp{\cU}
\def\UQp{U}
\def\UbK{\ol U}
\def\UF{\UR_s}

\def\SR{\cV}
\def\SK{S}
\def\SbK{\ol S}
\def\SF{\SR_s}

\def\CR{\cC}
\def\CK{C}
\def\CZp{\cC}
\def\CQp{C}
\def\CbK{\ol C}

\def\DR{\cD}
\def\DK{D}
\def\DZp{\cD}
\def\DQp{D}
\def\DbK{\ol D}
\def\tDZp{\widetilde{\cD}}

%                                 Greece
%
\def\ep{\epsilon}
\def\vG{\varGamma}
\def\vg{\varGamma}
%
%
%
%
%
%                                 simple
%
%
%
\def\lra{\longrightarrow}
\def\lla{\longleftarrow}
\def\Lra{\Longrightarrow}
\def\hra{\hookrightarrow}
\def\lmt{\longmapsto}
\def\ot{\otimes}
\def\op{\oplus}
%                              decolation
\def\wt#1{\widetilde{#1}}
\def\wh#1{\widehat{#1}}
\def\spt{\sptilde}
\def\ol#1{\overline{#1}}
\def\ul#1{\underline{#1}}
\def\us#1#2{\underset{#1}{#2}}
\def\os#1#2{\overset{#1}{#2}}
\def\lim#1{\us{#1}{\varinjlim}}
\def\plim#1{\varprojlim_{#1}}
%
%                         extra-ordinary symbol
%
%
%
\def\Gm{{\mathbb G}_{\hspace{-1pt}\mathrm{m}}}
\def\Ga{{\mathbb G}_{\hspace{-1pt}\mathrm{a}}}
\def\zp{{\bZ_p}}
\def\qp{{\bQ_p}}
\def\qzp{\qp/\zp}
\def\ql{{\bQ_\ell}}
\def\zl{{\bZ_\ell}}
\def\qzl{\ql/\zl}
\def\qz'{{\bQ}/{\bZ}'}
\def\isom{\hspace{9pt}{}^\sim\hspace{-16.5pt}\lra}
\def\lisom{\hspace{10pt}{}^\sim\hspace{-17.5pt}\lla}
%{\hspace{9pt}\sptilde\hspace{-16.5pt}\lra}

\def\hA{\wh{A}}
\def\uA{\mathscr A_0}
\def\uuA{\mathscr A}
\def\uB{\mathscr B_0}
\def\uuB{{\mathscr B}}
\def\uK{{{\mathscr R}^\flat_K}}
\def\uoK{{{\mathscr R}^\flat_{\ol K}}}
\def\uuK{{\uuR_K}}
\def\uR{{{\mathscr R}^\flat}}
\def\uuR{{\mathscr R}}
\def\tD{\widetilde{D}}
\def\tT{\widetilde{T}}
\def\tZ{\widetilde{Z}}

\def\psr{R\hspace{.5pt}[\hspace{-1.57pt}[q_0]\hspace{-1.58pt}]}
\def\psro{R\hspace{.5pt}[\hspace{-1.57pt}[q_1]\hspace{-1.58pt}]}
\def\psf{R\hspace{.5pt}(\hspace{-1.57pt}(q_0)\hspace{-1.58pt})}
\def\psfo{R\hspace{.5pt}(\hspace{-1.57pt}(q_1)\hspace{-1.58pt})}
\def\psru{\hspace{.5pt}[\hspace{-1.57pt}[u]\hspace{-1.58pt}]}
\def\psrt{\hspace{.5pt}[\hspace{-1.57pt}[t]\hspace{-1.58pt}]}
\def\psrti{\hspace{.5pt}[\hspace{-1.57pt}[t_i]\hspace{-1.58pt}]}
\def\psrq{\hspace{.5pt}[\hspace{-1.57pt}[q_0]\hspace{-1.58pt}]}
\def\psfu{\hspace{.5pt}(\hspace{-1.57pt}(u)\hspace{-1.58pt})}
\def\psfq{\hspace{.5pt}(\hspace{-1.57pt}(q_0)\hspace{-1.58pt})}
\def\psrqq{\hspace{.5pt}[\hspace{-1.57pt}[q]\hspace{-1.58pt}]}
\def\psrqqq{\hspace{.5pt}[\hspace{-1.57pt}[q_0,q_0^{-1}]\hspace{-1.58pt}]}
\def\psfqq{\hspace{.5pt}(\hspace{-1.57pt}(q)\hspace{-1.58pt})}
\def\psrqu{\hspace{.5pt}[\hspace{-1.57pt}[q_0,u]\hspace{-1.58pt}]}
\def\psfti{\hspace{.5pt}(\hspace{-1.57pt}(t_i)\hspace{-1.58pt})}
\def\psftt{\hspace{.5pt}(\hspace{-1.57pt}(t-1)\hspace{-1.58pt})}
\def\psfqi{\hspace{.5pt}(\hspace{-1.57pt}(q_i)\hspace{-1.58pt})}

\def\tet{{\tau^\et_\infty}}
\def\tdR{{\tau^\dR_\infty}}
\def\tsyn{{\tau^\syn_\infty}}
\def\htsyn{\wh{\tau}_{\infty}^{\,\syn}}
\def\htet{\wh{\tau}_{\infty}^{\,\et}}
\def\ttsyn{\wt{\tau}_{\infty}^{\syn}}
\def\ttet{\wt{\tau}_{\infty}^{\et}}
%
%
%
%                             Hodge-Witt sheaf
%
\def\mwitt#1#2#3{W_{\hspace{-2pt}#2}{\hspace{1pt}}\omega_{#1}^{#3}}
\def\witt#1#2#3{W_{\hspace{-2pt}#2}{\hspace{1pt}}\Omega_{#1}^{#3}}
\def\mlogwitt#1#2#3{W_{\hspace{-2pt}#2}{\hspace{1pt}}\omega_{{#1},{\log}}^{#3}}
\def\tlogwitt#1#2#3{W_{\hspace{-2pt}#2}{\hspace{1pt}}\tom_{{#1},{\log}}^{#3}}
\def\logwitt#1#2#3{W_{\hspace{-2pt}#2}{\hspace{1pt}}\Omega_{{#1},{\log}}^{#3}}
\def\loglogwitt#1#2{W_{\hspace{-2pt}n}{\hspace{1pt}}\Omega^{#2}_{Y^{(#1)}\hspace{-1.5pt},\hspace{1pt}E^{[#1]}\hspace{-1pt},\hspace{1pt}\log}}

\section{Introduction}
The Tate conjecture for Chow groups is one of the most important question in arithmetic geometry. It appears in various areas explicitly or implicitly and often plays a key role. Although there have been significant progresses on this conjecture  (e.g., for abelian varieties), very few things are known about its analogue for higher $K$-theory raised by Beilinson (cf.\ \cite{jan} 5.19). In this paper we are mainly concerned with the Tate conjecture for $K_2$, which asserts that for a nonsingular variety $U$ over a number field $F$ (or more generally
a finitely generated field over a prime field)
the \'etale chern class map \[ c_\et:K_2(U)\ot\qp\lra H_\et^{2}(\ol U,\qp(2))^{G_F} \] is surjective. Here $G_F:={\mathrm{Gal}}(\ol{F}/F)$ is the absolute Galois group, $\ol U$ denotes $U \otimes_F \ol F$ and the superscript ${G_F}$ denotes the fixed part by $G_F$.

In this paper, we focus on an elliptic surface $\pi:X\to C$ over a $p$-adic local field $K$ which is absolutely unramified.
Let $D=\sum \, D_i$ be the sum of the split multiplicative fibers of $\pi$, and put $U=X-D$.
Assume that $X$ and $C$ have projective smooth models $\cX$ and $\cC$ over the integer ring $R$ of $K$, respectively,
  and that the closure $\cD \subset \cX$ of $D$ has normal crossings. Put $\cU:=\cX-\cD$.
Then we introduce a space of {\it formal Eisenstein series}
\[ \Eis(\cX,\cD)_\zp \subset \vg(\cX,\Omega_{\cX/R}^2(\log\cD)) \]
 (see \S\ref{sect5-2} for details), where the right hand side is the space of global 2-forms with log poles along $\cD$.
One of our main results asserts that
\begin{equation}\label{eq0-0}
\Image\big(H^2_\syn(\cX(\cD),\cS_\zp(2)) \to \vg(\cX,\Omega_{\cX/R}^2(\log\cD))\big)\subset \Eis(\cX,\cD)_\zp,
\end{equation}
where $H^*_\syn(\cX(\cD),\cS_\zp(2))$ denotes the syntomic cohomology of $\cX$ with log poles along $\cD$
 due to Kato and Tsuji.
From this result, we will further deduce inequalities
\begin{equation}\label{eq0-1} \dim_\qp \, c_\et (K_2(U)\otimes \bQ_p) \le
\dim_\qp \, H^2_\et(\ol U,\bQ_p(2))^{G_K} \le \rank_\zp \, \Eis(\cX,\cD)_\zp 
\end{equation}
using the fact that the Fontaine-Messing map $H^2_\syn(\cX(\cD),\cS_\zp(2)) \to H^2_\et(\ol U,\zp(2))^{G_K}$ 
is surjective (\S\ref{sect5-3}).
The inclusion \eqref{eq0-0} is an extension of Beilinson's theorem on Eisenstein symbols in the following sense.
When $\pi:X\to X(\Gamma)$ is the universal family of elliptic curves over a modular curve $X(\Gamma)$,
the space $\Eis(\cX,\cD)_\zp$ consists of the (usual) Eisenstein series of weight 3 (\cite{A} \S 8.3).
Beilinson proved that it is spanned by the dlog image of the Eisenstein symbols in $K_2(U)$ (\cite{bei-modular}), so that
\[ \Image\big(K_2(U)\ot\qp \to \vg(X,\Omega_{X/K}^2(D))\big)= \Eis(\cX,\cD)_{\qp}=\qp^{\,\# \text{ of cusps}} \]
under our notation.
We can thus regard \eqref{eq0-0} as a partial extension of Beilinson's theorem to arbitrary elliptic surfaces.
A distinguished feature is that it gives a new upper bound of the rank of $H^2_\et(\ol U,\bQ_p(2))^{G_F}$
 in view of the fact that the rank of $\Eis(\cX,\cD)_\zp$ is strictly less than the number of split multiplicative fibers in some cases (\S \ref{application-sect},
see also \cite{tata} \S5).
The proof of \eqref{eq0-0} is quite different from that of the theorem of Beilinson.
A key ingredient is the $p$-adic Hodge theory, in particular a detailed computation on the syntomic cohomology of Tate curves over 2-dimensional complete local rings.

The inequality \eqref{eq0-1} 
is a key tool in our application to Beilinson's Tate conjecture for $K_2$.
In fact, if one can construct enough elements in $K_2(U)$ (e.g.\ by symbols) % (even in ad hoc manner!)
 so that the dimension of $c_\et(K_2(U))$ is
 equal to the upper bound, then the equalities hold in \eqref{eq0-1} and 
one can conclude that $c_\et$ is surjective.
In \S \ref{beilinson-sect} we will give a number of such examples.
It is remarkable that as a consequence
one has non-trivial elements in the {\it Selmer group}
$H^1_f(G_K, H^2_\et(\ol X,\Q_p(2))/\NS(\ol{X}))$ of Bloch-Kato \cite{BK2}
with large rank (\S \ref{selmer-sect}).

In \S \ref{mordellsect}, 
we present another application of \eqref{eq0-0} to the finiteness of torsion $0$-cycles.
It is a folklore conjecture that $\CH^m(X)$ is finitely generated $\Z$-module
% and hence the torsion part $\CH^m(X)_\tors$ is finite,
for a projective smooth variety $X$ over a number field,
which is a widely open problem unless $m=1$.
The finiteness of torsion part $\CH^m(X)_{\mathrm{tor}}$ supports the question.
When the base field is a $p$-adic local field, 
the Chow group is no longer finitely generated. 
Although most people had believed that the finiteness of torsion cycles remains true even for $p$-adic local fields, counter-examples were found recently (\cite{RS}, \cite{ASai}). On the other hand all such examples are {\it not} defined over number fields,
i.e. the minimal fields of definition are not number fields. Therefore we are naturally lead to the following modified question:
\begin{quote}
{\it Let $X_0$ be a projective nonsingular variety over a number field $F$.
Let $K$ be a completion of $F$ at a finite place $v$ and put $X=X_0\otimes_F K$.
Then is $\CH^m(X)_\tors$ finite?}
\end{quote}
It is in fact a crucial question whether the $p$-primary torsion part $\CH^m(X)\{p\}$ is finite.
When $m=2$, the finiteness of $\CH^2(X)\{p\}$ is reduced to the study of the $p$-adic regulator on $K_1$ by a recent work of Saito and the second author \cite{sato-saito}, that is, $\CH^2(X)\{p\}$ is finite if the $p$-adic regulator map \[ \varrho:K_1(X)^{(2)}\ot\qp\lra H^1_g(G_K, H^2_\et(\ol X,\qp(2)))\] is surjective onto the $g$-part of Bloch-Kato \cite{BK2} (3.7).
When $H^2(X,\O_X)=0$, this map is well-known to be surjective even when the minimal field of definition is not a number field (cf.\ \cite{CTR1}, \cite{CTR2}, \cite{S} 3.6). However when $H^2(X,\O_X) \ne 0$, the question becomes more difficult, and nobody has found an affirmative or negative example so far. In fact, the Selmer group 
$H^1_f(G_K, H^2_\et(\ol X,\Q_p(2))/\NS(\ol{X}))$ of Bloch-Kato is no longer zero (cf.\ Lemma \ref{lem6-2-2}), so one needs to construct {\it integral indecomposable} elements of $K_1$ generating the Selmer group over $\qp$, which is a crucial difficulty there.
We achieve it in the same way as in \S \ref{selmer-sect}
and give the first example of a surface $X$ over a $p$-adic local field with $H^2(X,\O_X) \ne 0$
such that the $p$-adic regulator map $\varrho$ is surjective and hence
the torsion of $\CH^2(X)$ is finite
(Theorem \ref{finite1}, Corollary \ref{finite2}, \S\ref{sect6-3}).

\medskip

This paper is organized as follows.
In \S \ref{sect1}, we review and fix the notation for (log) syntomic cohomology and Tate curves.
In \S \ref{sect2}, we state the main result on Tate curves (Theorem \ref{thm2-1}) and prove it admitting  a key commutative key diagram.
\S \ref{sect3} is devoted to the proof of the key diagram.
In \S \ref{knownsect} we prove the main results on elliptic surfaces over $p$-adic fields.
In \S \ref{application-sect}, we apply them to the study
of Beilinson's Tate conjecture for $K_2$.
In \S \ref{mordellsect} we give an example of elliptic $K$3 surface over $\qp$ with finitely many torsion $0$-cycles.

\newpage 
\section{Preliminaries}\label{sect1}

For a scheme $X$ over a ring $A$ and an $A$-algebra $B$
we write $X_B:=X\times_A B$.
For an integer $n$ which is invertible on $X$,
 $\bZ/n(1)$ denotes the \'etale sheaf $\mu_n$ of $n$-th roots of unity.
We often write $\bZ/n(m)$\,($m \in \bN$) for the \'etale sheaf $\mu_{n}^{\otimes m}$.
For a function $f \in \vG(X,\cO_X)$ which is not a zero divisor, we put
\[ X[f^{-1}] := {\boldsymbol \Spec}(\cO_X[T]/(fT-1)),  \]
which is the maximal open subset of $X$ where $f$ is invertible.

\subsection{Syntomic cohomology}\label{sect1-1}
Let $p$ be a prime number. For a scheme $T$, we put
\[ T_n := T \otimes \bZ/p^n. \]
% In this subsection, we always assume $0 \le r \le p-1$.
\begin{defn}\label{def1-1-0}
Let $T$ be a scheme.
\begin{enumerate}
\item[{\rm (1)}]
A morphism $\varphi : T \to T$ over $\bF_p$ is called the absolute Frobenius endomorphism,
 if the underlying morphism of topological spaces is the identitiy map and
 the homomorphism $\varphi^* : \cO_T \to \varphi_*\cO_T=\cO_T$ sends $x \mapsto x^p$.
Here the equality $\varphi_*\cO_T=\cO_T$ means the natural identification.
\item[{\rm (2)}]
A morphism $\varphi : T \to T$ over $\bZ_p$ is called a Frobenius endomorphism,
 if $\varphi \otimes \bZ/p : T_1 \to T_1$ is the absolute Frobenius endomorphism.
\end{enumerate}
\end{defn}
\par
Let $X$ be a scheme which is flat over $\bZ_p$.
Assume the following condition:
\begin{cond}\label{cond1-1-1}
There exists a closed immersion $X \hra Z$ satisfying the following conditions
 {\rm(}cf.\ {\rm\cite{Ka} 2.4):}
\begin{enumerate}
\item[{\rm (1)}]
$Z_n$ has $p$-bases over $\bZ/p^n$ locally for any $n \ge 1$ {\rm(}loc.\ cit.\ Definition {\rm1.3)}.
\item[{\rm (2)}]
$Z$ has a Frobenius endomorphism, or more weakly,
$Z_n$ has a Frobenius endomorphism $\varphi_n$ for each $n \ge 1$
 and the morphism $\varphi_{n+1} \otimes_{\bZ/p^{n+1}} \bZ/p^n$ agrees with $\varphi_n$ for any $n \ge 1$.
\item[{\rm (3)}]
Let $D_n$ be the PD-envelope of $X_n$ in $Z_n$ compatible with the canonical PD-structure on
 the ideal $(p) \subset \bZ/p^n$.
For $i \ge 1$, let $J_{D_n}^{[i]} \subset \cO_{D_n}$ be the $i$-th divided power of the ideal
$J_{D_n}:=\ker(\cO_{D_n} \to \cO_{X_n})$.
For $i \le 0$, we put $J_{D_n}^{[i]}:=\cO_{D_n}$.
% The above conditions (1)--(4) verify the exactness of
Then the following sequence is exact for any $m,n \ge 1$ and any $i \ge 0${\rm :}
\addtocounter{equation}{2}
\begin{equation}\label{eq1-1-1}
 J_{D_{m+n}}^{[i]} \os{\times p^m}\lra J_{D_{m+n}}^{[i]} \os{\times p^n}\lra J_{D_{m+n}}^{[i]} \lra J_{D_n}^{[i]} \lra 0.
\end{equation}
\end{enumerate}
\end{cond}
We give some examples of $Z$ which satisfy the condition (1):
\stepcounter{thm}
\begin{exmp}\label{ex1-1-1}
Let $W=W(k)$ be the Witt ring of a perfect field $k$, and let $W \psrt$ be the formal power series ring over $W$ with $t$ an indeterminate. Then $W_n \psrt := W \psrt/(p^n)$ has a $p$-basis over $\bZ/p^n$ for any $n \ge 1$. Indeed, $W_n \psrt$ is flat over $\bZ/p^n$ and $t \in k \psrt$ is a $p$-basis of $k \psrt$ over $\bF_p$. Therefore $t \in W_n \psrt$ is a $p$-basis over $\bZ/p^n$ by {\rm \cite{KaL}} Proposition {\rm1.4}. More generally, $Z_n$ has $p$-bases locally over $\bZ/p^n$ in the following cases{\rm:}
\begin{enumerate}
\item[{\rm (i)}]
$Z$ is smooth of finite type over $W \psrt$.
\item[{\rm (ii)}]
$Z$ is flat of finite type over $W \psrt$ and $Z_1$ is a regular semistable family over $k \psrt$.
\end{enumerate}
\end{exmp}
The following fact provides a sufficient condition for $X \hra Z$ to satisfy the condition (3):
\begin{prop}[{\bf Fontaine-Messing\,/\,Kato \cite{Ka} Lemma 2.1}]\label{prop1-1-kato}
Let $X \hra Z$ be a closed immersion satisfying the conditions {\rm (1)} and {\rm (2)} in Condition {\rm \ref{cond1-1-1}}.
Assume that $Z$ is locally noetherian, and
 $\ker(\cO_{Z,x} \to \cO_{X,x})$ is generated by an $\cO_{Z,x}$-regular sequence for any $x \in X_1$.
Then $X \hra Z$ satisfies the condition {\rm (3)} in Condition {\rm \ref{cond1-1-1}}.
\end{prop}
\par
For $0 \le r \le p-1$, we define the syntomic complex $\cS_n(r)_{X,Z}$ (with respect to the embedding $X \hra Z$)
 as follows.
Let $\bJ^{[r]}_{n,X,Z}$ be the complex of sheaves on $(X_1)_\et$
\[ J_{D_n}^{[r]} \os{d}\lra J_{D_n}^{[r-1]}\otimes_{\cO_{Z_n}} \Omega^1_{Z_n} \os{d}\lra \dotsb \os{d}\lra
   J_{D_n}^{[r-q]}\otimes_{\cO_{Z_n}} \Omega^q_{Z_n} \os{d}\lra \dotsb, \]
where $J_{D_n}^{[r]}$ is placed in degree $0$. Put $\bE_{n,X,Z}:=\bJ^{[0]}_{n,X,Z}$.
By the assumption that $0 \le r \le p-1$, the Frobenius endomorphism on $Z_{n+r}$ induces a homomorphism of complexes
\[ f_r:=\ol{p^{-r} \cdot \varphi_{n+r}^*} :  \bJ^{[r]}_{n,X,Z} \lra \bE_{n,X,Z} \]
(see \cite{Ka} p.\ 411 for details).
\begin{defn}\label{def1-1-1}
For $0 \le r \le p-1$, we define the complex $\cS_n(r)_{X,Z}$ on $(X_1)_\et$ as the mapping fiber of
\[ 1-f_r:  \bJ^{[r]}_{n,X,Z} \lra \bE_{n,X,Z}. \]
More precisely, the degree $q$-part of $\cS_n(r)_{X,Z}$ is
\[ (J_{D_n}^{[r-q]}\otimes_{\cO_{Z_n}} \Omega^q_{Z_n}) \oplus  (\cO_{D_n}\otimes_{\cO_{Z_n}} \Omega^{q-1}_{Z_n}) \]
and the differential operator is given by
\[ (x,y) \longmapsto (dx, x-f_r(x)-dy). \]
We define the syntomic cohomology of $X$ with coefficients in $\cS_n(r)$ as
 the hypercohomology groups of this complex$:$
\[ H^*_{\syn}(X,\cS_n(r)):=\bH^*_{\et}(X,\cS_n(r)_{X,Z}). \]
This notation is well-defined,
 because the image of the complex $\cS_n(r)_{X,Z}$ in the derived category
  is independent of $X \hra Z$ as in Condition {\rm \ref{cond1-1-1}}
 {\rm(\cite{Ka}} p.\ {\rm 412)}. See also Remark {\rm \ref{rem1-log-2}} below.
\end{defn}
\begin{prop}\label{prop1-1-0}
For $m,n \ge 1$, there is an exact sequence of complexes on $(X_1)_\et$
\[ \cS_{m+n}(r)_{X,Z} \os{\times p^m}\lra \cS_{m+n}(r)_{X,Z} \os{\times p^n}\lra \cS_{m+n}(r)_{X,Z} \lra \cS_n(r)_{X,Z} \lra 0. \]
Consequently, there is a short exact sequence
\[ 0 \lra \cS_m(r)_{X,Z} \os{\ul {p}^n}\lra \cS_{m+n}(r)_{X,Z} \lra \cS_n(r)_{X,Z} \lra 0. \]
\end{prop}
\begin{pf}
The assertion follows from the exactness of \eqref{eq1-1-1}
 and the fact that $\Omega^q_{Z_n}$ is free over $\bZ/p^n$ (\cite{Ka} Lemma 1.8).
\end{pf}
%\stepcounter{thm}
\begin{rem}\label{rem1-1-1}
Assume that the identity map $X \to X$ satisfies the conditions {\rm (1)}--{\rm (3)} in Condition {\rm \ref{cond1-1-1}}.
Then we have $D_n=X_n$, i.e.,
\[ \bJ^{[r]}_{n,X,X}=\Omega_{X_n}^{\bullet \ge r}\quad \hbox{ and } \quad
 \bE_{n,X,X}=\Omega_{X_n}^{\bullet}, \]
and there is a short exact sequence of complexes
\addtocounter{equation}{5}
\begin{equation}\label{eq1-1-2}
\xymatrix{
 0 \ar[r] & \cS_n(r)_{X,X} \ar[r] & \Omega_{X_n}^{\bullet \ge r} \ar[r]^-{1-f_r} & \Omega_{X_n}^{\bullet} \ar[r] & 0.}
\end{equation}
\end{rem}
% \medskip
\stepcounter{thm}
\begin{rem}
Syntomic complexes can be defined in a more general situation by a gluing argument in the derived category
 {\rm(}see {\rm \cite{KaV}} Remark {\rm 1.8} and {\rm \cite{Ka}} Lemma {\rm2.2)}.
% However we will not need this fact in this paper .
% We will use syntomic cohomology of such a smooth schemes over $p$-adic integer rings only in \S\ref{exmpsect2}.
\end{rem}
\subsection{Syntomic cohomology with log poles}\label{sect1-log}
The notation remain as in \S\ref{sect1-1}.
The aim of this subsection is to define of the syntomic cohomology of a regular scheme $X$
 with log poles along a simple normal crossing divisor $\cD$ following Tsuji (\cite{Ts2} \S2),
 which is necessary to formulate Theorem \ref{thm2-1} below.
The case that $\cD$ is empty corresponds to the syntomic cohomology defined in \S\ref{sect1-1}.
To give precise arguments, we will use the terminology in log geometry \cite{KaL}.
% The reader who does not want to get into the details on log geometry can skip this subsection.
\begin{defn}
Let $(T,M_T)$ be a log scheme.
\begin{enumerate}
\item[{\rm (1)}]
{\rm (\cite{KaL}} Definition {\rm 4.7)}
A morphism $\varphi : (T,M_T) \to (T,M_T)$ over $\bF_p$ is called the absolute Frobenius endomorphism,
 if the underlying morphism $T \to T$ is the absolute Frobenius endomorphism in the sense of Definition {\rm\ref{def1-1-0}\,(1)}
  and the homomorphism $\varphi^* : M_T \to \varphi_*M_T=M_T$ is the multiplication by $p$.
Here the equality $\varphi_*M_T=M_T$ means the natural identification obtained from the fact
 that the underlying morphism of $\varphi$ of topological spaces is the identity map.
\item[{\rm (2)}]
A morphism $\varphi : (T,M_T) \to (T,M_T)$ over $\bZ_p$ is called a Frobenius endomorphism, if
 $\varphi\otimes \bZ/p : (T_1,M_{T_1}) \to (T_1,M_{T_1})$ is the absolute Frobenius endomorphism.
Here $M_{T_1}$ denotes the inverse image log structure of $M_T$ onto $T_1$ {\rm(}loc.\ cit.\ {\rm(1.4))}.
\end{enumerate}
\end{defn}
\par
Let $(X,M)$ be a fine log scheme such that $X$ is flat over $\bZ_p$ (as a usual scheme).
The main example we are concerned with is the following case:
\begin{exmp}\label{ex1-2-1}
Let $X$ be a regular scheme which is flat over $\bZ_p$.
Let $\cD$ be a simple normal crossing divisor on $X$, which may be empty.
Put $U:=X-\cD$ and let $j : U \hra X$ be the natural open immersion.
We define the sheaf $M$ of monoids on $X_\et$ as
\[ M := \cO_X \cap j_*\cO_{U}^\times. \]
The canonical map $M \to \cO_X$ is a fine log structure on $X$, which we call the log structure associated with $\cD$.
\end{exmp}

For $n \ge 1$, we write $M_n$ for the inverse image log structure of $M$ onto $X_n$.
We assume the following condition, which is a logarithmic variant of Condition \ref{cond1-1-1}:
\begin{cond}\label{cond1-1-2}
There exist exact closed immersions $i_n : (X_n,M_n) \hra (Z_n,M_{Z_n})$ of log schemes for $n \ge 1$
 which satisfy the following conditions for all $n \ge 1$\,{\rm(}cf.\ {\rm\cite{Ts2} \S2):}
\begin{enumerate}
\item[{\rm (0)}]
$(Z_n,M_{Z_n})$ is fine. We have
\[ (Z_n,M_{Z_n}) \isom (Z_{n+1},M_{Z_{n+1}}) \otimes_{\bZ/p^{n+1}} \bZ/p^n \]
as log schemes, and the morphism $i_{n+1} \otimes_{\bZ/p^{n+1}} \bZ/p^n$ agrees with $i_n$.
\item[{\rm (1)}]
$(Z_n,M_{Z_n})$ has $p$-bases locally over $\bZ/p^n$ with the trivial log structure {\rm(}loc.\ cit.\ Definition {\rm1.4)}.
\item[{\rm (2)}]
$(Z_n,M_{Z_n})$ has a Frobenius endomorphism $\varphi_n$
 and the morphism $\varphi_{n+1} \otimes_{\bZ/p^{n+1}} \bZ/p^n$ agrees with $\varphi_n$.
\item[{\rm (3)}]
Let $(D_n,M_{D_n})$ be the PD-envelope of $(X_n,M_n)$ in $(Z_n,M_{Z_n})$
 which is compatible with the canonical PD-structure on the ideal $(p) \subset \bZ/p^n$
 {\rm (\cite{KaL}} Definition {\rm(5.4))}.
For $i \ge 1$, let $J_{D_n}^{[i]} \subset \cO_{D_n}$ be the $i$-th divided power of the ideal
 $J_{D_n}=\ker(\cO_{D_n} \to \cO_{X_n})$. For $i \le 0$, we put $J_{D_n}^{[i]}:=\cO_{D_n}$.
Then the following sequence is exact for any $m,n \ge 1$ and any $i \ge 0${\rm :}
\addtocounter{equation}{3}
\begin{equation}\label{eq1-log-1}
 J_{D_{m+n}}^{[i]} \os{\times p^m}\lra J_{D_{m+n}}^{[i]} \os{\times p^n}\lra J_{D_{m+n}}^{[i]} \lra J_{D_n}^{[i]} \lra 0.
\end{equation}
\end{enumerate}
\end{cond}
We give some examples of $(Z_n,M_{Z_n})$ which satisfy the condition (1):
\stepcounter{thm}
\begin{exmp}\label{ex1-2-2}
Let $W \psrt$ be as in Example {\rm \ref{ex1-1-1}}.
We endow $\Spec(W \psrt)$ with a pre-log structure $\bN \to W \psrt$ by sending $1$ to $t$,
  and write $N$ for the associated log structure on $\Spec(W \psrt)$.
Then $(W_n \psrt,N_n)$ has a $p$-basis over $\bZ/p^n$ for any $n \ge 1$,
  which one can check in a similar way as in Example {\rm \ref{ex1-1-1}}.
More generally, a fine log scheme $(Z_n,M_{Z_n})$ which is log smooth over $(W_n \psrt,N_n)$
 has $p$-bases locally over $\bZ/p^n$.
Indeed, $(Z_n,M_{Z_n})$ has $p$-bases locally over $(W_n \psrt,N_n)$ by {\rm\cite{Ts2}} Lemma {\rm 1.5}.
Hence  $(Z_n,M_{Z_n})$ has $p$-bases locally over $\bZ/p^n$ by loc.\ cit.\ Proposition {\rm 1.6\,(2)}.
\end{exmp}
The following fact is a logarithmic variant of Propostion \ref{prop1-1-kato}:
\begin{prop}[{\bf Tsuji \cite{Ts2} Corollary 1.9}]\label{prop1-log-tsuji}
Let $\{ i_n : (X_n,M_n) \hra (Z_n,M_{Z_n}) \}_{n \ge 1}$ be a system of exact closed immersions
 satisfying the conditions {\rm (0)}--{\rm (2)} in Condition {\rm \ref{cond1-1-2}}.
Assume that $Z_n$ is locally noetherian, and that
 $\ker(\cO_{Z_n,x} \to \cO_{X_n,x})$ is generated by an $\cO_{Z_n,x}$-regular sequence
 for any $x \in X_1$.
Then $\{ i_n : (X_n,M_n) \hra (Z_n,M_{Z_n}) \}_{n \ge 1}$ satisfies the condition {\rm (3)} in Condition {\rm \ref{cond1-1-2}}.
\end{prop}
\par
For $0 \le r \le p-1$, we define a log syntomic complex $\cS_n(r)_{(X,M),(Z_*,M_{Z_*})}$ as follows.
Let $\bJ^{[r]}_{n,(X,M),(Z_*,M_{Z_*})}$ be the complex of sheaves on $(X_1)_\et$
\[ J_{D_n}^{[r]} \os{d}\lra J_{D_n}^{[r-1]}\otimes_{\cO_{Z_n}} \omega^1_{(Z_n,M_{Z_n})} \os{d}\lra \dotsb \os{d}\lra
   J_{D_n}^{[r-q]}\otimes_{\cO_{Z_n}} \omega^q_{(Z_n,M_{Z_n})} \os{d}\lra \dotsb, \]
where $J_{D_n}^{[r]}$ is placed in degree $0$ and $\omega^*_{(Z_n,M_{Z_n})}$ denotes the differential module
 of $(Z_n,M_{Z_n})$ \cite{KaL} (1.7). The arrows $d$ denote the derivations defined in \cite{Ts2} Corollary 1.10.
Put \[ \bE_{n,(X,M),(Z_*,M_{Z_*})}:=\bJ^{[0]}_{n,(X,M),(Z_*,M_{Z_*})}. \]
By the assumption that $0 \le r \le p-1$, the Frobenius endomorphism on $(Z_{n+r},M_{Z_{n+r}})$
 induces a homomorphism of complexes
\[ f_r:=\ol{p^{-r} \cdot \varphi_{n+r}^*} :  \bJ^{[r]}_{n,(X,M),(Z_*,M_{Z_*})} \lra \bE_{n,(X,M),(Z_*,M_{Z_*})} \]
(see \cite{Ts2} p.\ 540 for details).
\begin{defn}\label{def1-2-1}
For $0 \le r \le p-1$, we define the complex $\cS_n(r)_{(X,M),(Z_*,M_{Z_*})}$ on $(X_1)_\et$ as the mapping fiber of
\[ 1-f_r:  \bJ^{[r]}_{n,(X,M),(Z_*,M_{Z_*})} \lra \bE_{n,(X,M),(Z_*,M_{Z_*})} \]
{\rm(}cf.\ Definition {\rm\ref{def1-1-1})}.
We define the syntomic cohomology of $(X,M)$ with coefficients in $\cS_n(r)$ as
 the hypercohomology groups of this complex$:$
\[ H^*_{\syn}((X,M),\cS_n(r)):=\bH^*_{\et}(X,\cS_n(r)_{(X,M),(Z_*,M_{Z_*})}). \]
This notation is well-defined, because the image of the complex $\cS_n(r)_{(X,M),(Z_*,M_{Z_*})}$ in the derived category
  is independent of embedding systems as in Condition {\rm \ref{cond1-1-2}} by Remark {\rm \ref{rem1-log-2}} below.
In the case of Example {\rm \ref{ex1-2-1}}, we put
\[ H^*_{\syn}(X(\cD),\cS_n(r)):=H^*_{\syn}((X,M),\cS_n(r)) \]
for simplicity.
\end{defn}
\begin{rem}\label{rem1-log-2}
Take another embedding system $\{i'_n:(X_n,M_n) \hra (Z'_n,M_{Z'_n})\}_{n \ge 1}$ as in Condition {\rm \ref{cond1-1-2}},
and consider the embedding system
\[ \big\{ i_n \times i'_n: (X_n,M_n) \hra (Z_n,M_{Z_n})\times_{\bZ/p^n}(Z'_n,M_{Z'_n})=:(Z''_n,M_{Z''_n}) \big\}{}_{n \ge 1}. \]
We define a Frobenius endomorphism on $(Z''_n,M_{Z''_n})$\,{\rm(}$n \ge 1${\rm)} as the fiber product of 
 those of $(Z_n,M_{Z_n})$ and $(Z'_n,M_{Z'_n})$.
Then this embedding system satisfies the conditions {\rm (0)}--{\rm (3)} in Condition {\rm \ref{cond1-1-2}} as well
 {\rm(}see {\rm \cite{Ka}} Proof of Lemma {\rm 2.2}, and {\rm \cite{Ts2}} Proposition {\rm 1.8)},
 and there are natural quasi-isomorphisms of complexes on $(X_1)_\et$
\[ \cS_n(r)_{(X,M),(Z_*,M_{Z_*})} \os{\qis}\lra \cS_n(r)_{(X,M),(Z''_*,M_{Z''_*})}
 \os{\qis}\longleftarrow \cS_n(r)_{(X,M),(Z'_*,M_{Z'_*})} \] by {\rm \cite{Ts2}} Corollary {\rm 1.11}.
Hence the image of the complex $\cS_n(r)_{(X,M),(Z_*,M_{Z_*})}$ in the derived category is independent of embedding systems.
Moreover, this fact verifies that log syntomic cohomology groups are contravariantly functorial for morphisms $(X,M) \to (X',M')$ of log schemes which satisfy Condition {\rm \ref{cond1-1-2}} {\rm(}see also {\rm\cite{KaV}} p.\ {\rm 212)}. 
\end{rem}
\begin{defn}% [{\bf \cite{Ts1} \S2.2\,/\,\cite{KaV}}]
\label{def1-2-2}
For $r,r' \ge 0$ with $r+r' \le p-1$, we define a product structure
\[ \cS_n(r)_{(X,M),(Z_*,M_{Z_*})} \otimes \cS_n(r')_{(X,M),(Z_*,M_{Z_*})} \lra  \cS_n(r+r')_{(X,M),(Z_*,M_{Z_*})} \]
by
\[ (x,y) \otimes (x',y') \longmapsto (xx',(-1)^qxy'+f_{r'}(x')y), \]
where
\begin{align*}
(x,y) \in \cS_n(r)^q_{(X,M),(Z_*,M_{Z_*})} &= (J_{D_n}^{[r-q]}\otimes_{\cO_{Z_n}} \omega^q_{(Z_n,M_{Z_n})}) \oplus  (\cO_{D_n}\otimes_{\cO_{Z_n}} \omega^{q-1}_{(Z_n,M_{Z_n})}), \\
(x',y') \in \cS_n(r')^{q'}_{(X,M),(Z_*,M_{Z_*})} &= (J_{D_n}^{[r'-q']}\otimes_{\cO_{Z_n}} \omega^{q'}_{(Z_n,M_{Z_n})}) \oplus  (\cO_{D_n}\otimes_{\cO_{Z_n}} \omega^{q'-1}_{(Z_n,M_{Z_n})}).
\end{align*}
\end{defn}
\par
The following proposition follows from the same arguments as for Proposition \ref{prop1-1-0}
 (use \cite{Ts2} Corollary 1.9 instead of \cite{Ka} Lemma 1.8):
\begin{prop}
For $\{i_n : (X_n,M_n) \hra (Z_n,M_{Z_n})\}_{n \ge 1}$ as before and integers $m,n \ge 1$,
 there is a short exact sequence of complexes on $(X_1)_\et$
\[ 0 \lra \cS_m(r)_{(X,M),(Z_*,M_{Z_*})} \os{\ul {p}^n}\lra
 \cS_{m+n}(r)_{(X,M),(Z_*,M_{Z_*})} \lra \cS_n(r)_{(X,M),(Z_*,M_{Z_*})} \lra 0. \]
\end{prop}
\par The following standard fact relates the syntomic cohomology with de Rham cohomology.
\begin{prop}\label{prop1-log-1}
For $0 \le r \le p-1$ and $i \ge 0$, there is a canonical map
\[ \ep^i : H^i_\syn((X,M),\cS_n(r)) \lra  \bH^i_\et\big(X_n,\omega^{\bullet \ge r}_{(X_n,M_n)}\big), \]
which is compatible with product structures and contravariantly functorial in $(X,M)$.
\end{prop}
\begin{pf}
Fix a system of embeddings $(X_n,M_n) \hra (Z_n,M_{Z_n})$\,{\rm(}$n \ge 1${\rm)} as before.
There are natural maps
\[ \cS_n(r)_{(X,M),(Z_*,M_{Z_*})} \lra \bJ^{[r]}_{n,(X,M),(Z_*,M_{Z_*})} \lra \omega^{\bullet \ge r}_{(X_n,M_n)}, \]
where the first arrow arises from the definition of the syntomic complex,
  and the second arrow is induced by the natural maps
\[  \cO_{X_n} \otimes_{\cO_{D_n}} \big(J_{D_n}^{[r-q]}\otimes_{\cO_{Z_n}} \omega^q_{(Z_n,M_{Z_n})}\big)
  \lra \begin{cases} 0 \quad & (q \le r-1) \\ \omega^q_{(X_n,M_n)} \quad &(q \ge r). \end{cases} \]
Let $\ep_{n,(X,M),(Z_*,M_{Z_*})}$ be the composite of the above natural maps of complexes, and define the desired map $\ep^i$
 as that induced by $\ep_{n,(X,M),(Z_*,M_{Z_*})}$.
One can easily check that the image of $\ep_{n,(X,M),(Z_*,M_{Z_*})}$ in the derived category is independent of the choice
 of an embedding system, by repeating the arguments in Remark \ref{rem1-log-2}. Thus we obtain the proposition.
\end{pf}
\begin{rem}\label{rem1-log-1}
If the identity maps $(X_n,M_n) \to (X_n,M_n)$\,{\rm(}$n \ge 1${\rm)} satisfy the conditions {\rm (0)}--{\rm (3)} in Condition {\rm \ref{cond1-1-2}},
 then we have $(D_n,M_{D_n})=(X_n,M_n)$, i.e.,
\[ \bJ^{[r]}_{n,(X,M),(X_*,M_*)}=\omega^{\bullet \ge r}_{(X_n,M_n)}\quad \hbox{ and } \quad
 \bE_{n,(X,M),(X_*,M_*)}=\omega^{\bullet}_{(X_n,M_n)}, \]
and there is a short exact sequence of complexes
\addtocounter{equation}{8}
\begin{equation}\label{eq1-log-2}
% 0 \lra \cS_n(r)_{(X,M),(X_*,M_*)} \lra \omega^{\bullet \ge r}_{(X_n,M_n)} \os{1-f_r}\lra \omega^{\bullet}_{(X_n,M_n)} \lra 0
\xymatrix{
 0 \ar[r] & \cS_n(r)_{(X,M),(X_*,M_*)} \ar[r]^-{\ep_n} & \omega_{(X_n,M_n)}^{\bullet \ge r} \ar[r]^-{1-f_r}
 & \omega_{(X_n,M_n)}^{\bullet} \ar[r] & 0.}
\end{equation}
for $0 \le r \le p-1$.
\end{rem}
\stepcounter{thm}
We provide a spectral sequence computing syntomic cohomology with log poles,
 which will be used in \S\ref{sect3-4} below.
\begin{prop}\label{prop1-1-1}
Let $X,\cD$ and $M$ be as in Example {\rm \ref{ex1-2-1}}, assume that $\cD$ is flat over $\bZ_p$. 
Let $\{ \cD_i \}_{i \in I}$ be the irreducible components of $\cD$.
Put $X^{(0)}:=X$ and
\[ X^{(m)}:=\coprod_{\{i_1,i_2,\dotsb,i_m \} \subset I}\ \cD_{i_1} \times_X \dotsb \times_X \cD_{i_m} \]
for $m \ge 1$, where for each subset $\{i_1,i_2,\dotsb, i_m\} \subset I$, the indices are pairwise distinct.  
Assume the following conditions$:$
\begin{enumerate}
\item[{\rm (i)}] The identity maps $(X_n,M_n) \to (X_n,M_n)$ satisfy the conditions {\rm (0)}--{\rm (3)}
 in Condition {\rm \ref{cond1-1-2}} for all $n \ge 0$.
\item[{\rm (ii)}] The identity maps $X^{(m)} \to X^{(m)}$ satisfy the conditions {\rm (1)}--{\rm (3)}
 in Condition {\rm \ref{cond1-1-1}} for all $m \ge 0$.
\item[{\rm (iii)}] The given Frobenius endomorphisms on $X_{n+r}$ and $(X^{(0)})_{n+r}$ are compatible under
    the canonical finite morphism $(X^{(0)})_{n+r} \to X_{n+r}$.
\item[{\rm (iv)}] For any $m \ge 0$,
 the given Frobenius endomorphisms on $(X^{(m)})_{n+r}$ and $(X^{(m+1)})_{n+r}$ are compatible under
    the canonical finite morphism $(X^{(m+1)})_{n+r} \to (X^{(m)})_{n+r}$.
\end{enumerate}
Fix an ordering on the set $I$.
Then for $0 \le r \le p-1$, there is a spectral sequence of syntomic cohomology groups
\stepcounter{equation}
\[ E_1^{a,b}=H^{2a+b}_{\syn}(X^{(-a)},\cS_n(a+r)) \Lra H^{a+b}_{\syn}(X(\cD),\cS_n(r)). \]
\end{prop}
\begin{pf}
By (i) and (ii), the syntomic complexes $\cS_n(r)_{(X,M),(X_*,M_*)}$ and $\cS_n(r)_{X^{(m)},X^{(m)}}$
are defined for $0 \le r \le p-1$ and $m \ge 0$.
These complexes are computed as in Remarks \ref{rem1-1-1} and \ref{rem1-log-1}.
By (iii) and (iv), there is a natural `filtration' on $\cS_n(r)_{(X,M),(X_*,M_*)}$ as follows:
\begin{eqnarray*}
 & 0 \lra \cS_n(r)_{X,X} \os{\alpha_0}\lra \cS_n(r)_{(X,M),(X_*,M_*)} \lra C^\bullet_{1,n} \lra 0 \\
 & 0 \lra \cS_n(r)_{X^{(m)},X^{(m)}}[-m] \os{\alpha_m}\lra C^{\bullet}_{m,n} \lra C^\bullet_{m+1,n} \lra 0
  & \quad (m \ge 1) \\
 & C^\bullet_{m,n} = 0 & \quad (\hbox{if }\;X^{(m)}=\emptyset),
\end{eqnarray*}
where we have omitted the indication of the direct image of sheaves under the canonical finite morphisms
  $(X^{(m)})_1 \to X_1$.
The arrow $\alpha_0$ denotes the natural inclusion of complexes, and the arrow $\alpha_m$ for $m \ge 1$
  is given by the alternate sum of the inverse of Poincar\'e residue mappings whose signs are determined by the fixed
 ordering on $I$.
The spectral sequence in question is obtained from this filtration.
\end{pf}
\par
The following fact relates the syntomic cohomology with \'etale cohomology and plays an important role in this paper.
\begin{thm}[{\bf Tsuji \cite{Ts1} \S3.1}]\label{thm1-1-1}
Let $X,\cD$ and $M$ be as in Example {\rm \ref{ex1-2-1}} {\rm(}$\cD$ may be empty{\rm)}.
Assume that $(X,M)$ satisfies Condition {\rm \ref{cond1-1-2}} and that there exists a henselian local ring $\cR$ which is faithfully flat over $\bZ_p$ and such that $X$ is proper over $\cR$. Put $U:=X-\cD$.
Then for $0 \le r \le p-2$ and $i \ge 0$, there is a canonical homomorphism
\[ c^i : H^i_\syn(X(\cD),\cS_n(r)) \lra H^i_\et(U[p^{-1}],\bZ/p^n(r)), \]
which is compatible with product structures and contravariantly functorial in $(X,M)$.
% See Remark {\rm \ref{rem1-log-2}} for the functoriality of log syntomic cohomology.
\end{thm}
\begin{pf}
Since $X$ is regular and flat over $\bZ_p$, we have a short exact sequence
\[\xymatrix{ 0 \ar[r] & \bZ_p\,t^{\{ r\}} \ar[r] & Fil_p^rA_{\crys}(\ol {A^h}) \ar[rr]^-{1-\frac{\varphi}{p^r}} && A_{\crys}(\ol {A^h}) \ar[r] & 0 }\]
for an affine open subset $\Spec(A) \subset X$ (see \cite{Ts1} p.\ 245).
We obtain a morphism
\[ c: \cS_n(r)_{(X,M)} \lra \iota^*Rj_*\bZ/p^n(r) \qquad (X_1 \os{\iota}\lra X \os{j}\lla U[p^{-1}]) \]
 in the derived category of \'etale sheaves on $X_1$ by repeating the arguments in loc.\ cit.\ pp.\ 316--322, and obtain the desired map $c^i$ by the proper base-change theorem: 
\[ c^i : H^i_\syn(X(\cD),\cS_n(r)) \os{c}\lra H^i_\et(X_1,\iota^*Rj_*\bZ/p^n(r)) \lisom H^i_\et(U[p^{-1}],\bZ/p^n(r)). \]
See \cite{Ts2} pp.\ 544--545 for the functoriality.
\end{pf}

\subsection{Symbol maps}\label{sect1-1'}
For a scheme $Z$, we write $\cO(Z)$ for $\vG(Z,\cO_Z)$, for simplicity.
We first review \'etale symbol maps.
Let $X$ be a scheme and let $n$ be a positive integer which is invertible on $X$.
Then we have a short exact sequence on $X_\et$, called {\it the Kummer sequence}
\[\xymatrix{ 0 \ar[r] & \bZ/n(1) \ar[r] & \cO_X^\times \ar[r]^{\times n} & \cO_X^\times \ar[r] & 0. }\]
See the beginning of \S\ref{sect1} for the definition of the \'etale sheaf $\bZ/n(1)$.
Taking \'etale cohomology groups, we get a connecting map
\begin{equation}\label{eq1-2-0}
\xymatrix{ \cO(X)^\times/n \; \ar@{^{(}->}[r]& H^1_{\et}(X,\bZ/n(1)). }
\end{equation}
We write $\{ x \}^\et$ for the image of $x \in \cO(X)^\times$ under this map.
Taking cup products, we obtain a map
\begin{equation}\label{eq1-2-1}
 (\cO(X)^\times)^{\otimes r}/n \lra H^r_{\et}(X,\bZ/n(r)),
\end{equation}
which sends $x_1 \otimes x_2 \otimes \dotsb \otimes x_r$\,(each $x_i \in \cO(X)^\times$) to
 $\{x_1\}^\et \cup \{x_2\}^\et \cup \dotsb \cup\{ x_r \}^\et$.
By an argument of Tate \cite{T} Proposition 2.1,
 this map annihilates Steinberg relations in $(\cO(X)^\times)^{\otimes r}$, i.e.,
  the elements of the form
\[ x_1 \otimes x_2 \otimes \dotsb \otimes x_r \; \hbox{ with } \; x_i+x_j=0 \hbox{ or } 1 \;
  \hbox{ for some } i \ne j \]
 map to $0$ under the map \eqref{eq1-2-1}. Consequently, we get a map
\begin{equation}\label{eq1-2-2}
 K^M_r(\cO(X))/n \lra H^r_{\et}(X,\bZ/n(r)),
\end{equation}
which we call {\it the \'etale symbol map}.
When $X$ is the spectrum of a field, we often call
  this map {\it the Galois symbol map}.
\addtocounter{thm}{3}
\begin{rem}\label{rem1-2-1}
\begin{enumerate}
\item[{\rm (1)}]
Since we have $H^1_\et(X,\cO_X^\times) \simeq \Pic(X)$ by Hilbert's theorem $90$,
  the map \eqref{eq1-2-0} is bijective if $X$ is the spectrum of a UFD {\rm (}e.g., a field{\rm)}. 
\item[{\rm (2)}]
If $r=2$ and $X=\Spec(F)$ with $F$ a field, then the map \eqref{eq1-2-2} is bijective by
  the Merkur'ev-Suslin theorem {\rm \cite{MS}}.
\end{enumerate}
\end{rem}
\par
\def\gp{{\mathrm{gp}}}
\def\cM{{\mathscr M}}
We next review syntomic symbol maps (\cite{FM} p.\ 205, \cite{KaV} Chapter I \S3, \cite{Ts1} \S2.2, \cite{Ts2} p.\ 542).
Let $(X,M)$ be a log scheme which is flat over $\bZ_p$ and satisfies Condition \ref{cond1-1-2}.
Fix an embedding system $\{ i : (X_n,M_n) \hra (Z_n,M_{Z_n}) \}_{n \ge 1}$ as in Condition \ref{cond1-1-2}.
We define the complex $C_n$ as
\[ C_n :=\big(1+J_{D_n} \lra M_{D_n}^\gp \big) \qquad \hbox{($1+J_{D_n}$ is placed in degree $0$)}, \]
where for a sheaf of $\cM$ of commutative monoids, $\cM^\gp$ denotes the associated sheaf of abelian groups.
We define the map of complexes
\[ s : C_{n+1}  \lra  \cS_n(1)_{(X,M),(Z_*,M_{Z_*})} \]
as the map
\[  s^0 : 1+J_{D_{n+1}} \lra J_{D_n}, \quad a \longmapsto \log(a) \]
in degree $0$, and the map
\[ \begin{CD}
 s^1 : M_{D_{n+1}}^\gp @. \ \lra \ @. \big(\cO_{D_n} \otimes_{\cO_{Z_n}}\omega^1_{(Z_n,M_{Z_n})}\big) \oplus\cO_{D_n}  \\
 \qquad b @.  \longmapsto @. \big(\dlog(b), p^{-1}\log(b^p\varphi_{n+1}(b)^{-1})\big)
\end{CD} \]
in degree $1$.
Here $\varphi_{n+1}(b)b^{-p}$ belongs to $1+p\cO_{D_{n+1}}$ and
 the logarithm $\log(b^p\varphi_{n+1}(b)^{-1})$ belongs to $p\cO_{D_{n+1}}$.
The notation `$p^{-1}$' means the inverse image under the isomorphism
\[ p : \cO_{D_n} \isom p\cO_{D_{n+1}} \qquad \hbox{(cf.\ \eqref{eq1-log-1})}. \]
One can easily check that the maps $s^0$ and $s^1$ yield a map of complexes.
Since there is a natural quasi-isomorphism $C_{n+1} \isom M_{n+1}^\gp[1]$,
 the map $s$ induces a morphism
\[   M_{n+1}^\gp[1] \lra \cS_n(1)_{(X,M),(Z_*,M_{Z_*})} \]
 in the derived category, which is independent of the choice of embedding systems.

Now suppose further that $X,\cD$ and $M$ be as in Example \ref{ex1-2-1},
 and let $j : U:=X-\cD \hra X$ be the natural open immersion.
Then we have $M^\gp=j_*\cO_U^\times$, and obtain a map
\stepcounter{equation}
\begin{equation}\label{eq1-2-3} \cO(U_{n+1})^\times \lra H^1_{\syn}(X(\cD),\cS_n(1)). \end{equation}
We often write $\{ x \}^\syn \in H^1_{\syn}(X(\cD),\cS_n(1))$ for the image of $x \in \cO(U_{n+1})^\times$.
This map and the product structure of syntomic cohomology (Definition \ref{def1-2-2}) give rise to a map
\begin{equation}\label{eq1-2-4}  K^M_r(\cO(U_{n+1})) \lra H^r_{\syn}(X(\cD),\cS_n(r)) \end{equation}
for $0 \le r \le p-1$ (cf.\ \cite{KaV} Proposition 3.2), which we call {\it the syntomic symbol map}.
\addtocounter{thm}{2}
\begin{rem}\label{rem1-2-2}
When $M$ is the trivial log structure $\cO_X^\times$, we obtain a symbol map
\[ K^M_r(\cO(X_{n+1})) \lra H^r_{\syn}(X,\cS_n(r)) \]
for $0 \le r \le p-1$.
If the identity map $X \to X$ satisfies the conditions {\rm (1)}--{\rm (3)} in Condition {\rm \ref{cond1-1-1}},
then the connecting homomorphism induced by \eqref{eq1-1-2}
\[ H^0_\dR(X_n) \lra H^1_{\syn}(X,\cS_n(1)) \]
sends $a \in H^0_\dR(X_n)$ to $\{1+pa\}^\syn$, where $1+pa$ is well-defined in $\cO(X_{n+1})^\times$.
One can easily check this fact directly from the definition of the symbol map.
\end{rem}
\begin{thm}[{Tsuji \bf \cite{Ts1} Proposition 3.2.4}]\label{prop1-2-1}
Let $X,\cD$ and $M$ be as in Example {\rm \ref{ex1-2-1}}, and assume that $(X,M)$ satisfies the assumptions in Theorem {\rm \ref{thm1-1-1}}.
Put $U:=X-\cD$.
Then there is a commutative diagram
\[ \xymatrix{
K^M_r(\cO(U)) \ar@{->>}[r] \ar[rd] & K^M_r(\cO(U_{n+1})) \ar[rr]^-{\eqref{eq1-2-4}}
 & & H^r_{\syn}(X(\cD),\cS_n(r)) \ar[d]^{c^r} \\
& K^M_r(\cO(U)[p^{-1}]) \ar[rr]^-{\eqref{eq1-2-2}} & & H^r_{\et}(U[p^{-1}],\bZ/p^n(r)) }
\]
for $0 \le r \le p-2$, where $c^r$ denotes the canonical map in Theorem {\rm \ref{thm1-1-1}}.
\end{thm}
\begin{pf}
The case $r=1$ follows from the same arguments as in \cite{Ts1} Proposition 3.2.4.
The general case follows from the previous case by the compatibility of these arrows with product structures
 (Theorem \ref{thm1-1-1}).
\end{pf}
\par
We state a syntomic analogue of the facts in Remark \ref{rem1-2-1}, which will be used in \S\ref{sect3-4} below.
\begin{thm}\label{lem1-2-1}
Let $R$ be a henselian discrete valuation ring whose fraction field $L$ has characteristic zero and whose residue field $F$ has characteristic $p$.
% Assume that $p$ is a prime element of $R$ and that $R$ has a Frobenius endomorphism.
Then the syntomic symbol map
\[ K^M_r(R)/p^n \lra H^r_{\syn}(R,\cS_n(r))\]
is surjective for $r=2$ {\rm(}and $p \ge 5${\rm)}, and bijective for $r=1$ {\rm(}and $p \ge 3${\rm)}.
\end{thm}
\begin{pf}
We prove only the case $r=2$. The case $r=1$ follows from a similar (and simpler) arguments as below and the details are left to the reader. We first show that the following sequence of Milnor $K$-groups is exact:
\addtocounter{equation}{3}
\begin{equation}\label{eq1-2-5} K^M_2(R) \lra K^M_2(L) \os{\partial}\lra K^M_1(F) \lra 0, \end{equation}
where the first arrow is the natural pull-back of symbols. The arrow $\partial$ is the boundary map of Milnor $K$-groups, which is obviously surjective. We show the exactness at $K^M_2(L)$.
Indeed, we have a localization sequence of algebraic $K$-groups \[ K_2(R) \lra K_2(L) \os{d}\lra K_1(F) \] and natural isomorphisms $K_2(L) \simeq K^M_2(L)$ and $K_1(F) \simeq F^\times$. The arrow $d$ agrees with $\partial$ (up to a sign) under these isomorphisms. Moreover the natural map $K^M_2(R) \to K_2(R)$ is surjective (\cite{Sr} p.\ 17 Remark). Hence the sequence \eqref{eq1-2-5} is exact at $K^M_2(L)$.

Put $\eta:=\Spec(F)$, and consider the following commutative diagram with exact rows:
\[ \xymatrix{ &K^M_2(R)/p^n \ar[r] \ar[d] &
K_2^M(L)/p^n \ar[r]^-{\partial} \ar[d]_{\wr\hspace{-1pt}}^{\eqref{eq1-2-2}} & F^\times/p^n \ar[r] \ar[d]^{\dlog}_{\wr\hspace{-1pt}} & 0 \\ 0 \ar[r]  &H^2_{\syn}(R,\cS_n(2)) \ar[r]^-{c^2} & H^2_{\et}(L,\bZ/p^n(2)) \ar[r]^-{\partial'} & H^0_{\et}(\eta,\logwitt{\eta} n 1), }\]
where the upper row is the exact sequence \eqref{eq1-2-5} modulo $p^r$.
The arrow $\partial'$ is the boundary map of Galois cohomology (\cite{KaH} \S1) and the lower row is exact by Kurihara \cite{Ku}. The central vertical arrow is bijective by Remark \ref{rem1-2-1}\,(2). The right vertical arrow is bijective by the short exact sequence \[\xymatrix{ 0 \ar[r] & \cO_\eta^\times \ar[r]^-{\times p^n} & \cO_\eta^\times \ar[r]^-{\dlog} &
 \logwitt{\eta} n 1 \ar[r] & 0 }\] on $\eta_{\et}$ and Remark \ref{rem1-2-1}\,(1).
Thus we obtain the lemma. \end{pf}
We end this section with the following standard fact on Milnor $K$-groups, which will be used in the proof of Proposition \ref{prop2-1} below:
\stepcounter{thm}
\begin{prop}\label{lem2-2}
Let  $R$ be a henselian discrete valuation ring whose residue field $F$ has characteristic $p$. Let $L$ be the fraction field of $R$. Assume that $F$ is infinite and that $\Omega_{F,\log}^2=0$. Then the natural map $K_2^M(R)/p^n \to K_2^M(L)/p^n$ is injective for any $n \ge 1$.
\end{prop}
\begin{pf}
Since $F$ is infinite, we have $K^M_2(R)=K_2(R)$ by a theorem of van der Kallen (\cite{Sr} p.\ 17 Remark).
%Stein, Matsumoto, and van der Kallen (cf.\ \cite{Sr} Theorem 1.14, ): 
There is a localization exact sequence of algebraic $K$-groups
\[ K_2(F) \lra K_2(R) \lra K_2(L) \lra K_1(F) (=F^\times). \]
The last arrow is surjective, as $R$ is a principal ideal domain. We decompose this sequence into two exact sequences
\begin{align*}
   K_2(F) \lra & \, K_2(R) \lra M \lra 0, \\
0 \lra M \lra & \, K_2(L) \lra F^\times \lra 0.
\end{align*}
Since $F^\times$ is $p$-torsion-free, we have $M/p^n \hra K_2(L)/p^n$.
On the other hand, we have \[ K_2(F)/p^n = K_2^M(F)/p^n = 0 \]
by the assumption that $\Omega_{F,\log}^2=0$ (\cite{BK} Theorem 2.1) and $K_2(R)/p^n \simeq M/p^n$. Hence
\[ K_2^M(R)/p^n = K_2(R)/p^n \simeq M/p^n \hra  K_2(L)/p^n= K_2^M(L)/p^n \]
as required.
\end{pf}
By Theorem \ref{lem1-2-1} and Proposition \ref{lem2-2}, we obtain the following consequence, which will not be used in the rest of this paper:
\begin{cor}
Let $R$ be as in Theorem {\rm\ref{lem1-2-1}}. Assume $p \ge 5$ and that the residue field $F$ satisfies the assumptions in Proposition {\rm\ref{lem2-2}}. Then the symbol map $K^M_2(R)/p^n \to H^2_{\syn}(R,\cS_n(2))$ is bijective.
\end{cor}

\subsection{Tate curve}\label{sect1-2}
Let $B$ be a noetherian complete local ring with $6^{-1} \in B$.
Let $q \in B$ be an element which is contained in the maximal ideal of $B$ and not nilpotent. Put
\[ A := B[q^{-1}]. \]
{\it The Tate curve $E=E_q$ over $A$ with period $q$} is the projective completion in $\bP^2_A$ of the affine curve
\[ y^2+xy=x^3+a_4(q)x+a_6(q) \quad \hbox{ on } \;\; \Spec(A[x,y]), \]
where $a_4(q)$ and $a_6(q) \in B$ are defined as follows:
\[ a_4(q)=-5\sum_{n=1}^\infty \frac{n^3q^n}{1-q^n},\quad
a_6(q)=-\sum_{n=1}^\infty \frac{(5n^3+7n^5)q^n}{12(1-q^n)}. \]
We review the Tate parameterization of $E$.
Let $O \in E(A)$ be the infinite point.
The series
\begin{align*}
x(\alpha)&=\sum_{n\in\bZ}\; \frac{q^n\alpha}{(1-q^n\alpha)^2}-2\,\sum_{n\geq1} \; \frac{q^n}{(1-q^n)^2} \\
y(\alpha)&=\sum_{n\in\bZ}\; \frac{(q^n\alpha)^2}{(1-q^n\alpha)^3}+\sum_{n\geq1}\; \frac{q^n}{(1-q^n)^2}
\end{align*}
converge for all $\alpha \in A^\times-q^\bZ$. They induce an injective homomorphism
\[ A^\times/q^\bZ \lra E(A), \quad \alpha \longmapsto
\begin{cases}
(x(\alpha):y(\alpha):1) & (\alpha \not\in q^\bZ) \\
O & (\alpha \in q^\bZ),
\end{cases}
\]
which is bijective if $B$ is a complete discrete valuation ring (\cite{Si} V Theorem 3.1).
This assignment is not algebraic, but we have the following algebraic byproduct.
Let $u$ be an indeterminate, and put
\[ \uB := \varprojlim_n \ B[u,u^{-1}]/(q^n) \quad \hbox{ and } \quad \uA:=\uB[q^{-1}]. \]
\begin{prop}\label{prop1-1}
There is an injective homomorphism of $B$-algebras
\[ B[x,y]/(y^2+xy-x^3-a_4(q)x-a_6(q))  \lra \uB [(1-u)^{-1}] \]
given by the `$q$-adic presentation' of $(x(u),y(u))${\rm:}
\begin{align*}
 x & \longmapsto \frac{u}{(1-u)^2}+\sum_{d \ge 0} \; \sum_{m|d} \; m\big(u^m+u^{-m}-2\big) q^d \\
 y & \longmapsto \frac{u^2}{(1-u)^3}+\sum_{d \ge 0}\; \sum_{m|d} \;
 m\bigg(\frac{m-1}{2}\,u^m-\frac{m+1}{2}\,u^{-m}+1\bigg)q^d.
\end{align*}
Moreover, this ring homomorphism induces a dominant morphism of schemes
\[ \beta_0 : \Spec(\uA) \to E.  \]
\end{prop}
\begin{pf}
The assertion follows from the same arguments as in \cite{Si} p.\ 425 Proof of Theorem 3.1\,(c).
\end{pf}

\begin{defn}\label{def1-4-1}
\begin{enumerate}
\item[{\rm (1)}]
We define the canonical invariant $1$-form $\omega_E \in \vG(E,\Omega^1_E)$ as
\[  \omega_E:=\frac{dx}{2y+x}, \]
which maps to $u^{-1}du$ under the pull-back map
\[\xymatrix{ \beta_0^* : \vG(E,\Omega^1_E) \; \ar@{^{(}->}[r] & \; \Omega^1_{\uA}. }\]
% \begin{equation}\label{CI}  \omega_E \longmapsto \frac{du}{u}. \end{equation}
\item[{\rm (2)}]
We define the theta function $\theta(u,q) \in \uB$ as
\addtocounter{equation}{2}
\begin{equation}\label{def:theta}
\theta(u,q):=(1-u)\prod_{n=1}^\infty(1-q^nu)(1-q^nu^{-1}),
\end{equation}
which satisfies
\[ \theta(qu,q)=\theta(u^{-1},q)=-u^{-1}\theta(u,q). \]
\end{enumerate}
\end{defn}
\addtocounter{thm}{1}
The following proposition follows from standard facts on theta functions and theta divisors, whose details are left to the reader.
\begin{prop}\label{prop1-2}
Put $K:=\Frac(B)$, and let $K(E)$ be the function field of $E$.
\begin{enumerate}
\item[{\rm (1)}]
A function $f(u) \in \Frac(\uB)$ given by a finite product
\stepcounter{equation}
\begin{equation}\label{11}
f(u)=c \, \prod_i \; \frac{\theta(\alpha_iu,q)}{\theta(\beta_iu,q)} \qquad (c,\alpha_i,\beta_i \in A^\times)
\end{equation}
is $q$-periodic, if $\prod_i\ {\alpha_i}/{\beta_i}=1$.
\item[{\rm (2)}]
The ring homomorphism in Proposition {\rm \ref{prop1-1}} induces a natural inclusion
\[\xymatrix{ \displaystyle \left\{ c \, \prod_i \; \frac{\theta(\alpha_iu,q)}{\theta(\beta_iu,q)} \ \bigg| \ c,\alpha_i,\beta_i \in A^\times \text{ with } \prod_i \; {\alpha_i}/{\beta_i}=1 \right\} \; \ar@{^{(}->}[r] & K(E)^\times. }\]
\end{enumerate}
\end{prop}
\stepcounter{thm}
\begin{prop}\label{prop1-3}
Let $f$ be the structural morphism $E[p^{-1}] \to \Spec(A[p^{-1}])$.
% By Tate's parameterization of Tate curves, 
Then there is an exact sequence of \'etale sheaves on $\Spec(A[p^{-1}])$
\stepcounter{equation}
\begin{equation}\label{eq2-0}
0 \lra \bZ/p^n(1) \lra R^1f_*\bZ/p^n(1) \lra \bZ/p^n \lra 0.
\end{equation}
\end{prop}
\begin{pf}
See \cite{DR} VII.1.13.
\end{pf}
\subsection{Frobenius endomorphism on Tate curves}\label{sect1-3}
Let $p$ be a prime number at least $5$, and let $B,A$ and $q$ be as in the beginning of \S\ref{sect1-2}.
We assume here that $B$ is a $\bZ_p$-algebra, and that $B$ has a Frobenius endomorphism $\phi$.
Here a Frobenius endomorphism means a
   ring endomorphism compatible with the absolute Frobenius endomorphism on $B/(p)$.
Let $E_q=E_{q,A}$ be the Tate curve over $A$ with period $q$.
There is a canonical morphism
\[ \can :  E_{\phi(q)} \lra E_q \]
induced by the ring homomorphism
\begin{align*} A[x,y]/(y^2+xy-x^3-\,&a_4(q)x-a_6(q)) \\
 & \lra A[x,y]/(y^2+xy-x^3-a_4(\phi(q))x-a_6(\phi(q)))
 \end{align*}
sending
\[ x \mapsto x,\qquad y \mapsto y,\qquad a \mapsto \phi(a) \quad(a \in A). \]
We define $W:E_{\phi(q)} \to E_q$ as the composite
\[\xymatrix{ W : E_{\phi(q)} \ar[r]^-{\can} & E_q  \ar[r]^-{[p]} & E_q, }\]
where $[p]$ denotes the multiplication by $p$ (with respect to the group structure on $E_q$).
\begin{rem}\label{rem1-1}
Under the map in Proposition {\rm \ref{prop1-1}}, the map $W$ corresponds to the endomorphism of $\Frac(\uB)$ sending
\[ u \mapsto u^p,\qquad b \mapsto \phi(b) \quad (b \in B). \]
\end{rem}
\par
In what follows, we assume that \[ \phi(q)=q^p \quad \hbox{ and } \quad \phi(a)=a \quad (a\in \bZ_p). \] There is an isogeny over $\Spec(A)$ \[ \iota :  E_{q^p} \lra E_q \] corresponding to the identity map of $\Frac(\uB)$, which is an analogue of  the natural projection $\bC^\times/q^{p\bZ}\to \bC^\times/q^\bZ$. The morphism $W$ factors through $\iota$. Indeed, $\iota$ is surjective and we have \[ \ker(\iota) \subset \ker([p]:E_{q^p} \to E_{q^p}) \] as finite \'etale group schemes over $\Spec(A)$. Thus $W$ gives rise to an endomorphism
\stepcounter{equation}
\begin{equation}\label{eq1-2} \varphi : E_q \lra E_q. \end{equation}
We show that $\varphi$ is a Frobenius endomorphism:
\addtocounter{thm}{1}
\begin{lem}\label{lem1-1}
Put $(E_q)_1:=E_q \otimes_A A/(p)=E_q \otimes \bF_p$.
Then the morphism $\varphi_1 : (E_q)_1 \to (E_q)_1$ induced by $\varphi$
 is the absolute Frobenius endomorphism of $(E_q)_1$.
\end{lem}
\begin{pf}
Define the Frobenius endomorphism $\phi$ on $\bZ_p \psfqq := \bZ_p \psrqq [q^{-1}]$
 as $\phi(q)=q^p$ and $\phi(a)=a$ for $a \in \bZ_p$.
The natural ring homomorphism $\bZ_p \psfqq \to A$ is compatible with the Frobenius endomorphisms.
Note that $E_q=E_{q,A}$ is obtained from $E'=E_{q,\bZ_p\psfqq}$ by scalar extension, and that we have
\[ \varphi =\varphi' \otimes \phi \]
  under the identification $E_q=E' \otimes_{\bZ_p\psfqq} A$.
Here $\varphi'$ mean the map \eqref{eq1-2} defined for $E'$.
Therefore it is enough to consider the case that $A=\bZ_p \psfqq$.
Then $\varphi_1$ agrees with $\varphi''$, the map \eqref{eq1-2} defined for $E''=E_{q,\bF_p\psfqq}$.
Under the morphism $\beta_0$ in Proposition \ref{prop1-1} with $B=\bF_p\psrqq$, $\varphi''$ corresponds to the absolute Frobenius endomorphism of the field
\[ \uA=\uB[q^{-1}]=\left(\varprojlim_n \ \bF_p\psrqq[u,u^{-1}]/(q^n)\right)[q^{-1}] \] (cf.\ Remark \ref{rem1-1}), which shows that $\varphi''$ is the absolute Frobenius morphism.
\end{pf}

\newpage
\section{De Rham regulator of Tate curves}\label{sect2}
In this section, we assume $p \ge 5$.
The main result of this section will be stated in Theorem \ref{thm2-1} below.
\subsection{Setting}\label{sect2-1}
Let $R$ be a $p$-adic integer ring which is unramified over $\bZ_p$,
 and let $k$ be the residue field of $R$. Let $q_0$ be an indeterminate.
We define the rings $A$ and $B$ as
\[  B :=\psr, \qquad A :=B[q_0^{-1}]. \]
$B$ is a 2-dimensional regular complete local domain, and $A$ is a Dedekind domain.
Let $\hA$ be the $p$-adic completion of $A$: 
\[ \hA := \varprojlim_{n \ge 1} \ A/p^n, \]
which is a complete discrete valuation ring whose maximal ideal is generated by $p$.
Put \[ L:=\hA\,[p^{-1}]. \]
Let $r$ be a positive integer prime to $p$, and put $q : = q_0^{r}$.
Let $E=E_q$ be the Tate curve over $\Spec(A)$ with period $q$.
Let $\cE'$ be the projective curve over $\Spec(B)$
 defined by the following homogeneous equation in $\bP^2_B=\Proj(B[x,y,z])$:
\[ \cE' : y^2z+xyz=x^3+a_4(q)\,xz^2+a_6(q)\,z^3 \]
(see \S\ref{sect1-1} for $a_4(q)$ and $a_6(q)$),
which is a projective flat model of $E$ over $B$.
By blowing-up $\cE'$ along the locus $\{x=y=q_0=0\}$ up to $(r-1)$-times, we get a regular scheme $\cE$,
 which is a generalized elliptic curve in the sense of \cite{DR} II.1.12.
The divisor $\cD:=\{q_0 = 0\} \subset \cE$ is the standard N\'eron $r$-gon over $\Spec(R)$,
 and the structural morphism $\cE \to \Spec(B)$ is smooth outside of the intersection loci of
 two distinct irreducible components of $\cD$.
There is a cartesian diagram
\[ \xymatrix{ E \; \ar@{^{(}->}[r] \ar[d]_{\pi_E} \ar@{}[rd]|{\square} & \cE  \ar[d]^{\pi_\cE} \\
\Spec(A) \; \ar@{^{(}->}[r] & \Spec(B), } \]
where $\pi_E$ is projective smooth and $\pi_\cE$ is projective flat.
The horizontal arrows are open embeddings.
\begin{lem}\label{lem2-1-1}
Let $M$ be the log structure on $\cE$ associated with $\cD$, as in Example {\rm \ref{ex1-2-1}}.
Then the identity maps $(\cE_n,M_n) \to (\cE_n,M_n)$\,$(n \ge 1)$ satisfy the conditions {\rm (0)}--{\rm (3)} in
  Condition {\rm \ref{cond1-1-2}}.
Consequently, the log syntomic cohomology $H^*_{\syn}(\cE(\cD),\cS_n(r))$ is defined
 for $0 \le r \le p-1$ {\rm (}Definition {\rm \ref{def1-2-1})}.
\end{lem}
\begin{pf}
The condition (0) is obvious, and (3) follows from Proposition \ref{prop1-log-tsuji}.
The condition (1) follows from Example \ref{ex1-2-2}, and (2) follows from Lemma \ref{lem3-1-1} below.
\end{pf}
\begin{rem}\label{rem2-1-1}
Similarly, the identity map $\cE \to \cE$ satisfies the conditions
 {\rm (1)}--{\rm (3)} in Condition {\rm \ref{cond1-1-1}}, by Example {\rm \ref{ex1-1-1}} and Lemma {\rm\ref{lem3-1-1}} below.
Consequently, the syntomic cohomology $H^*_{\syn}(\cE,\cS_n(r))$\,$(0 \le r \le p-1)$ is computed by
  the distinguished triangle
\[\xymatrix{ 0 \ar[r] & \cS_n(r)_{\cE,\cE} \ar[r] & \Omega_{\cE_n}^{\bullet \ge r} \ar[r]^-{1-f_r} &
 \Omega_{\cE_n}^{\bullet} \ar[r] & 0 }\]
by Remark {\rm\ref{rem1-1-1}}.
\end{rem}
We put
\[ H^*_{\syn}(-,\cS_{\bZ_p}(r)):=\varprojlim_{n \ge 1} H^*_{\syn}(-,\cS_n(r)), \qquad H^*_{\et}(-,\bZ_p(r)):=\varprojlim_{n \ge 1} H^*_{\et}(-,\bZ/p^n(r)). \]
We define the map of K\"ahler differential forms
\[ \tdR : \vG(E,\Omega^2_E) \lra \Omega^1_A \] as the composite of natural isomorphisms
\[ \vG(E,\Omega^2_E) \simeq \vG(E,\Omega^1_A \otimes_A \Omega^1_{E/A}) \simeq \Omega^1_A \otimes_A \vG(E,\Omega^1_{E/A}) \] and the map
\[\Omega^1_A \otimes_A \vG(E,\Omega^1_{E/A}) \lra \Omega^1_A, \quad \eta \otimes \omega_E \mapsto \eta \quad(\eta \in \Omega^1_A).  \] Here $\omega_E$ denotes the canonical invariant $1$-form on $E$ defined in Definition \ref{def1-4-1}\,(1), and we have used the fact that $\vG(E,\Omega^1_{E/A})$ is a free $A$-module generated by $\omega_E$.
% Let $\cD$ be the regular divisor on $\cE$ defined by $q_0$.
%
%
\subsection{Main result on Tate curves}\label{sect2-2}
Let $\zeta_1,\dotsc,\zeta_d\in R$ be roots of unity which form a basis of $R$ over $\bZ_p$,
 whose existence is verified by the assumption that $R$ is absolutely unramified. An arbitrary formal Laurant power series 
$f(q_0) \in R\psrqqq$ is expanded into a power series of the form
\begin{equation}\label{eq2-3}
 f(q_0)  = \sum_{j \leq 0}\, b_j \, q_0^j + \sum_{j \ge 1} \bigg( \frac{a_{1j}\,\zeta_1\, q_0^j}{1-\zeta_1\,q_0^j}+\cdots+ \frac{a_{dj}\,\zeta_d\,q_0^j}{1-\zeta_d\,q_0^j}\bigg) \quad (a_{ij} \in \bZ_p,\; b_j \in R)
\end{equation}
and the coefficients $a_{ij}$ and $b_j$ are uniquely determined by $f(q_0)$.
We say that $f(q_0)$ is a {\it formal power series of Eisenstein type} (in $R\psrqqq$) if 
\begin{enumerate}
\renewcommand{\labelenumi}{\;\;(E\theenumi)}
\item
$c_k=0$ for $k<0$ and $c_0\in \zp$,
\item
$a^{(j)}_k\in k^2\zp$ for all $j$ and $k\geq1$.
\end{enumerate}
The condition (E2) does not depend on the choice of 
$\zeta_i$ (\cite{A} Lemma 3.4).
Moreover, $f(q_0) \in R\psrqqq$ is of Eisenstein type if and only if so is it in $R'\psrqqq$
for an $p$-adic integer ring $R'$ which is unramified over $R$.
\begin{rem}\label{rem2-1}
The reason we call ``Eisenstein" is the following fact.
Suppose that $f(q_0)$ is the $q_0$-expansion of a modular form of weight $3$ at a cusp.
Then $f(q_0)$ is of Eisenstein type if and only if it is a linear combination of the usual Eisenstein series of weight 3
{\rm(\cite{A}} \S{\rm 8.3)}.
\end{rem}

The main result of this section deals with the image of a `de Rham regulator map' from syntomic cohomology
\[ \begin{CD} \reg_{\dR} : H^2_{\syn}(\cE(\cD),
\cS_{\bZ_p}(2)) @>{\ep^2}>> \displaystyle \varprojlim_{n \ge 1} \ \vG(E_n,\Omega^2_{E_n}) @>{\tdR}>> \Omega^1_{\hA}  @>{f \cdot \frac{dq_0}{q_0} \mapsto f}>> \hA, \end{CD}\]
where $\tdR$ denote the map defined in \S\ref{sect2-1}. See Proposition \ref{prop1-log-1} for $\ep^2$.
\stepcounter{thm}
\begin{thm}\label{thm2-1}
Assume that $f(q_0) \in \hA$ is contained in the image of $\reg_{\dR}$.
Then it is a formal power series of Eisenstein type.
\end{thm}
%\begin{rem}\label{rem2-1}
%\begin{enumerate}
%\item[{\rm (1)}]
%The composite of the syntomic Chern class map and $\reg_{\dR}${\rm :} \[\xymatr%ix{ K_2(E) \ar[r] & H^2_{\syn}(E,\cS_{\bZ_p}(2)) \ar[r]^-{\reg_{\dR}} & \hA }\]% is close to the regulator map to Deligne-Beilinson cohomology.
%This is why we call $\reg_{\dR}$ the de Rham regulator.
%\item[{\rm (2)}]
%The assertion of Theorem {\rm\ref{thm2-1}} is obvious if $(j,p)=1$.
%\end{enumerate}
%\end{rem}
\subsection{A commutative diagram}
In the rest of this section, we prove Theorem \ref{thm2-1} assuming a commutative diagram \eqref{eq2-1} below.
For $n \ge 1$, let $\tet$ be the composite of canonical maps of \'etale cohomology groups
\[ \tet : H^2_{\et}(E_L,\bZ/p^n(2)) \lra H^1_{\et}(L,H^1_\et(E_{\ol L},\bZ/p^n(2))) \lra H^1_{\et}(L,\bZ/p^n(1)). \]
Here the second arrow is induced by the map $H^1_\et(E_{\ol L},\bZ/p^n(1)) \to \bZ/p^n$ in \eqref{eq2-0}. The first arrow is obtained from the Hochschild-Serre spectral sequence
\[ E_2^{a,b}=H^a_\et(L,H^b_\et(E_{\ol L},\bZ/p^n(2))) \Lra H^{a+b}_{\et}(E_L,\bZ/p^n(2)) \]
and the fact that \[ E_2^{0,2} \simeq H^0_{\et}(L, \bZ/p^n(1))=0, \] where we have used the assumption that $R$ is unramified over $\bZ_p$.
\par
 A key ingredient of the proof of Theorem {\rm\ref{thm2-1}} is the following commutative diagram:
\begin{equation}\label{eq2-1}
\xymatrix{
H^2_{\et}(E_L,\bZ_p(2)) \ar[rr]^-{\tet} & & H^1_{\et}(L,\bZ_p(1)) \\
H^2_{\syn}(\cE(\cD),\cS_{\bZ_p}(2)) \ar[u]^{c^2} \ar[rr]^-{\tsyn} \ar[d]_{\ep^2} & &
H^1_{\syn}(A,\cS_{\bZ_p}(1)) \ar[u]_{c^1} \ar[d]^{\dlog}\\
 \varprojlim{}_{n \ge 1} \ \vG(E_n,\Omega^2_{E_n}) \ar[rr]^-{\tdR} & & \Omega^1_{\hA}. }
\end{equation}
Here $c^1$ and $c^2$ are canonical maps in Theorem \ref{thm1-1-1}, and $\dlog$ is given by logarithmic differentials (see Theorem \ref{lem1-2-1} for the isomorphism $s$):
\begin{equation*}
\begin{CD} H^1_{\syn}(A,\cS_{\bZ_p}(1)) \to H^1_{\syn}(\hA,\cS_{\bZ_p}(1)) \os{s}\simeq \displaystyle \varprojlim_{n \ge 1} \ \hA^\times/p^n @>{f \mapsto \frac{df}{f}}>> \displaystyle \varprojlim_{n \ge 1} \ \Omega^1_{\hA/p^n} =
 \Omega^1_{\hA}\;. \end{CD}\end{equation*}
We will construct $\tsyn$ in \S\ref{sect3-2} and prove the commutativity of the squares in \S\ref{sect3-4} below.
\par
As a preliminary of the proof of Theorem {\rm\ref{thm2-1}}, we prove Lemma \ref{lem2-1} below.
For $n \ge 1$, let $\nu$ be the isomorphism
\[  \nu :  H^1_{\et}(L,\bZ/p^n(1)) \simeq L^\times/p^n. \]
in Remark \ref{rem1-2-1}\,(1).
\stepcounter{thm}
\begin{lem}\label{lem2-1}
For any $n \ge 0$, the composite map
\[\begin{CD} H^2_{\et}(E_L,\bZ/p^n(2)) @>{\nu \circ \tet}>> L^\times/p^n @>{a \mapsto \{a,q_0\}}>> K_2^M(L)/p^n \end{CD}\]
is zero.
\end{lem}
\begin{pf}
Since $q=q_0^{r}$ and $(r,p)=1$ by definition, we may replace the second arrow with the assignment $a \mapsto \{a,q\}$.
We consider the following commutative diagram, whose top row is an exact sequence arising from \eqref{eq2-0}:
\[\xymatrix{
H^1_{\et}(L,H^1_{\et}(E_{\ol L},\bZ/p^n(2))) \ar[r] & H^1_{\et}(L,\bZ/p^n(1)) \ar[r]^{(1)} & H^2_{\et}(L,\bZ/p^n(2)) \\
H^2_{\et}(E_L,\bZ/p^n(2)) \ar[r]^-{\nu \circ \tet} \ar[u] & L^\times/p^n \ar@{<-}[u]^{\wr\hspace{-1pt}}_{\nu} & K_2^M(L)/p^n \ar[u]^{\wr\hspace{-1pt}}_{(2)}. }\]
Here the arrow (2) is a Galois symbol map \eqref{eq1-2-2} (see also Remark \ref{rem1-2-1}).
%which are bijective by Hilbert's theorem 90 and the Merkur'ev-Suslin theorem \cite{MS}, respectively.
We claim that the arrow (1) maps $a \in H^1_{\et}(L,\bZ/p^n(1))$ to $a \cup \{ q \}^\et$ up to a sign.
Indeed, the connecting homomorphism \[ \bZ/p^n=H^0_{\et}(L,\bZ/p^n) \lra H^1_{\et}(L,\bZ/p^n(1)) \]
associated with \eqref{eq2-0} sends $1$ to $\{ q \}^\et$, and the claim follows from a straight-forward computation on cup products. The lemma follows from these facts.
\end{pf}
\par
\subsection{Proof of Theorem \ref{thm2-1}}
We prove of Theorem \ref{thm2-1}, assuming the diagram \eqref{eq2-1}.
Assume that $f(q_0) \in \hA$ lies in the image of $\reg_{\dR}$.
By the lower square of the diagram \eqref{eq2-1}, there exists $h(q_0) \in \hA^\times$ such that
\begin{equation}\label{eq2-4}
 \frac{dh(q_0)}{h(q_0)}=f(q_0)\,\frac{dq_0}{q_0} \in \Omega^1_{\hA}\;.
\end{equation}
We fix such an $h(q_0)$ in what follows.
By Lemma \ref{lem2-1} and the upper square of the diagram \eqref{eq2-1}, 
  $h(q_0)$ must satisfy
\begin{equation}\label{eq2-5}
 \{h(q_0),q_0\}=0 \in K^M_2(L)/p^n \; \hbox{ for any $n \ge 0$}.
\end{equation}
Now expand $f(q_0)$ into a series of the form \eqref{eq2-3}:
\begin{align*}
f(q_0) & = \sum_{j \leq 0}\, b_j\,q_0^j + 
\sum_{j \ge 1}\bigg(
\frac{a_{1j}\,\zeta_1\, q_0^j}{1-\zeta_1\,q_0^j}+\cdots+
\frac{a_{dj}\,\zeta_d\, q_0^j}{1-\zeta_d\,q_0^j} \bigg) \qquad (a_{ij} \in \bZ_p,\; b_j \in R) \\ 
\intertext{and expand $h(q_0)$ into an infinite product of the following form:}
h(q_0) & = c \, q_0^m
\prod_{j \ge 1} \,(1-\zeta_1\,q_0^j)^{a'_{1j}}\cdots (1-\zeta_d\,q_0^j)^{a'_{dj}}
\qquad (a'_{ij} \in \bZ_p,\; c \in R^\times,\;m\in\bZ).
\end{align*}
By \eqref{eq2-4}, we see that
\begin{equation}\label{eq2-6} a_{ij}=-ja'_{ij}, \end{equation}
hence $a_{ij}$ is divisible by $j$ for any $j \ge 1$.
We next  write down what the equation \eqref{eq2-5} yields,
  using the explicit reciprocity law of higher-dimensional regular local rings due to Kato \cite{Ka}:
\addtocounter{thm}{3}
\begin{prop}\label{prop2-1}
$a'_{ij}$ is divisible by $j$ in $\bZ_p$ for any $j \ge 0$.
\end{prop}
Theorem \ref{thm2-1} immediately follows from this proposition and \eqref{eq2-6}.
\par\medskip
\begin{pf*}{\it Proof of Proposition \ref{prop2-1}}
Note that the symbol $\{h(q_0),q_0\}$ is defined in $K_2^M\big(\hA \big)$.
Since $\Omega^2_{k\psfq}=0$, we have
\stepcounter{equation}
\begin{equation}\label{eq2-8}
 \{h(q_0),q_0\} = 0 \in K_2^M\big(\hA \big)/p^n \; \hbox{ for any $n \ge 1$},
\end{equation}
by \eqref{eq2-5} and Proposition \ref{lem2-2}. Let $\phi : \hA \to \hA$ be the Frobenius endomorphism of $\hA$ defined by
 the canonical Frobenius automorphism on $R$ and the assignment $q_0 \mapsto q_0^p$.
Define the function $\ell_\phi : \hA \to \hA$ as
\[\ell_\phi(a) := \frac{1}{\,p \,} \log\bigg(\frac{\phi(a)}{a^p}\bigg). \]
It follows from \eqref{eq2-8} and the explicit reciprocity law (\cite{Ka} Corollary 2.9), that
\[ \ell_\phi(h(q_0))\frac{dq_0}{q_0}-\ell_\phi(q_0)\frac{dh(q_0)}{h(q_0)} \in d\hA. \]
Since $\ell_\phi(q_0)=0$ by definition, there exists $\alpha \in \hA$ satisfying
\begin{equation}\label{eq2-9}
 \ell_\phi(h(q_0))=q_0 \frac{d\alpha}{dq_0}.
\end{equation}
Because $\ell_\phi(fg)=\ell_\phi(f)+\ell_\phi(g)$ and
\[ \ell_\phi(1-rq_0)=\sum_{(n,p)=1} \, \frac{(rq_0)^n}{n} \qquad (r \in R) \]
by definition, we have
\[ \ell_\phi(h(q_0))= \ell_\phi(c) + \sum_{j \ge 1}\sum_{(n,p)=1}
\bigg( a'_{1j}\frac{(\zeta_1\,q_0^j)^n}{n} + \dotsb +  a'_{dj}\frac{(\zeta_d\,q_0^j)^n}{n} \bigg). \]
Comparing these coefficients with those in the right hand side of \eqref{eq2-9},
 one can easily check that $a'_{ij}$ is divisible by $j$ in $\bZ_p$.
\end{pf*}
Thus we obtained Theorem \ref{thm2-1}, assuming the diagram \eqref{eq2-1}.
\newpage

\section{Construction of the key diagram}\label{sect3}
The notation remains as in \S\ref{sect2-1}.
In this section, we construct the key homomorphism
\[\tsyn : H_{\syn}^2(\cE(\cD),\cS_{\bZ_p}(2))\to H_{\syn}^1(A,\bZ_p(1)) \]
to establish the commutative diagram \eqref{eq2-1}.
\subsection{Preliminary}\label{sect3-1}
We define the Frobenius map $\phi : B \to B$ as
\begin{equation}\label{eq3-1}
  \phi(q_0)=q_0^p \quad \hbox{and} \quad \phi(a)=\sigma(a) \quad (a \in R),
\end{equation}
where $\sigma$ denotes the canonical Frobenius automorphism of $R$.
Let $\uuB$ (resp.\ $\uuR$) be the $(q_0,p)$-adic completion of $B[u,u^{-1}]$
 (resp.\ $p$-adic completion of $R[u,u^{-1}]$):
\[ \uuB := \varprojlim_n \ B[u,u^{-1}]/(q_0,p)^n,
\qquad \uuR :=\varprojlim_n \ R[u,u^{-1}]/(p^n). \]
We extend the Frobenius map $\phi$ on $B$ to $\uuB$ by defining $\phi(u):=u^p$.
We define the Frobenius map $\phi$ on $\uuR$ as $\phi(u):=u^p$ and $\phi(a)=\sigma(a)$\,$(a \in R)$.
The ring homomorphism in Proposition \ref{prop1-1} induces a morphism of schemes
\begin{equation}\label{eq3-3} \beta : \Spec(\uuB) \lra \cE,  \end{equation}
where $\cE$ is as we defined in \S\ref{sect2-1}.
For a scheme (or a ring) $Z$ and $n \ge 1$, we put
\[ Z_n:=Z \otimes \bZ/(p^n).\]
We have a Frobenius endomorphism $\varphi: E \to E$ by \eqref{eq3-1} and the construction in \S\ref{sect1-3}.
\addtocounter{thm}{2}
\begin{lem}\label{lem3-1-1}
For any integer $n \ge 1$, the Frobenius endomorphism $\varphi_n : E_n \to E_n$ extends to a
  Frobenius endomorphism $\varphi_n : (\cE_n,M_n) \to (\cE_n,M_n)$.
Furthermore, $\varphi_{n+1} \otimes \bZ/p^n$ agrees with $\varphi_n$ for any $n \ge 1$.
\end{lem}
This lemma has been used in Lemma \ref{lem2-1-1}.
\par
\medskip
\begin{pf}
Let $Z$ be the singular locus of $\cD$ (see \S\ref{sect2-1} for the definition of $\cD$).
Note that $X := \cE-Z$ is a commutative group scheme over $\Spec(B)$.
In particular, the multiplication by $p$ is defined on $X$ with respect to its group structure,
 and there is an endomorphism
\[ \varphi : X \lra X \]
defined in a similar way as for $\varphi$ on $E$.
By Lemma \ref{lem1-1}, $\varphi$ is a Frobenius endomorphism of $X$, i.e.,
$\varphi_1 : X_1 \to X_1$ induced by $\varphi$ is the absolute Frobenius of $X_1$,
because $X_1$ and $E_1=(E_q)_1$ have the same function field.
Since $\cD$ is defined by $q_0$ and $\varphi$ sends $q_0$ to $q_0^p$, the map $\varphi : X \to X$ induces
 a Frobenius endomorphism
\[ \varphi : (X,M_X) \lra (X,M_X) \]
of log schemes, where $M_X$ denotes the log structure on $X$ associated with the regular divisor $\cD-Z$
  (see Example \ref{ex1-2-1}).

We start the proof of the lemma.
Because a Frobenius endomorphism on $\cE_n$ is the identity map on the underlying topological space,
 we define an endomorphism $\psi : \cO_{\cE_n} \to \cO_{\cE_n}$
 which lifts the absolute endomorphism of $\cO_{\cE_1}$ and which is compatible with
  $\varphi_n^*$ on $\cO_{X_n}$.
Let $\alpha : X \hra \cE$ and $\alpha_n : X_n \hra \cE_n$ be the natural open embeddings.
We first show that the natural adjunction map
\stepcounter{equation}
\begin{equation}\label{eq3-4-4} \cO_{\cE_n} \lra \alpha_{n*}\cO_{X_n} \end{equation}
is bijective.
Consider the following commutative diagram with exact rows:
\[ \xymatrix{ 0 \ar[r] & \cO_\cE \ar[r]^-{\times p} \ar[d]_{\wr\hspace{-1pt}} & \cO_\cE \ar[r] \ar[d]_{\wr\hspace{-1pt}} &
 \cO_{\cE_1} \ar[r] \ar[d]_{\wr\hspace{-1pt}} & 0 \\
   0 \ar[r] & \alpha_*\cO_X \ar[r]^-{\times p} & \alpha_*\cO_X \ar[r] & \alpha_{1*}\cO_{X_1} \ar[r] &
   R^1\alpha_*\cO_X \ar[r]^-{\times p} & R^1\alpha_*\cO_X. }\]
Here the vertical arrows are natural adjunction maps, which are bijective by the facts that
 $\cE$ and $\cE_1$ are regular and that the complement $Z$ (resp.\ $Z_1$)
 has pure codimension $2$ in $\cE$ (resp.\ in $\cE_1$).
Hence $R^1\alpha_*\cO_X$ is $p$-torsion-free and we obtain the bijectivity of \eqref{eq3-4-4} by
  repeating a similar argument using $\cO_{\cE_n}$ and $\alpha_{n*}\cO_{X_n}$
 instead of $\cO_{\cE_1}$ and $\alpha_{1*}\cO_{X_1}$.
We define $\psi : \cO_{\cE_n} \to \cO_{\cE_n}$ as the composite map
\[\xymatrix{  \cO_{\cE_n} \simeq \alpha_{n*}\cO_{X_n} \ar[rr]^-{\alpha_{n*}(\varphi_n^*)} & &
 \alpha_{n*}(\cO_{X_n}) \simeq \cO_{\cE_n}.  } \]
It is clear that $\psi$ lifts the absolute Frobenius endomorphism of $\cO_{\cE_1}$.
Thus we obtain a Frobenius endomorphism $\varphi_n$ on $\cE_n$ which extends $\varphi_n$ on $X_n$, and
 it is easy to see that this induces a Frobenius endomorphism $\varphi_n$ on $(\cE_n,M_n)$.
By the construction, this morphism is the only morphism that extends $\varphi_n$ on $(X_n,M_{X_n})$,
 and the compatibility assertion $\varphi_{n+1} \otimes \bZ/p^n=\varphi_n$
 follows from the fact that $\varphi_n$ on $(X_n,M_{X_n})$ is induced by $\varphi$ on $(X,M_X)$. 
\end{pf}

\subsection{Syntomic residue mapping}\label{sect3-2}
To construct $\tsyn$, we first review the definition of residue mappings.
\begin{defn}
Let $m$ be a positive interger, and put $C:=B/(q_0,p)^m$ and $C\psfu:=C \psru [u^{-1}]$. We define residue mappings
\[ \Res_{u=0} : \Omega^r_{C\psfu} \lra \Omega^{r-1}_{C} \] as follows.
For $r=1$, we define $\Res_{u=0}$ as the composite map
\[ \Omega^1_{C\psfu} \os{\can}\lra \Omega^1_{C\psfu/C} \simeq C\psfu \cdot du \lra C \]
where the last arrows sends $c_iu^idu$\,{\rm ($c_i \in C$)} to $0$ if $i \ne -1$, and to $c_{-1}$ if $i=-1$.
For $r \ge 2$, we define $\Res_{u=0}$ as
\[\begin{CD} \Omega^r_{C\psfu} \os{\can}\lra \Omega^{r-1}_C\otimes \Omega^1_{C\psfu/C} @>{\id \otimes \Res_{u=0}}>> \Omega^{r-1}_C.
 \end{CD}\]
\end{defn}
Since $\uuA$ has a Frobenius endomorphism $\phi$ (cf.\ \S\ref{sect3-1}), the syntomic cohomology groups $H^*_\syn(\uuA,\cS_n(2))$ are computed by the complex
\begin{align*}
\cS_n(2)_{\uuA,\uuA} & := \mf \ \left(1-f_2: \Omega^{\bullet\geq 2}_{\uuA_n} \lra \Omega^{\bullet}_{\uuA_n}\right) \\
\intertext{ \hspace{-6pt} (see Remark \ref{rem1-1-1}). Similarly, $H^*_\syn(A,\cS_n(1))$ are computed by the complex}
\cS_n(1)_{A,A} & := \mf \ \left(1-f_1: \Omega^{\bullet\geq 1}_{A_n} \lra \Omega^{\bullet}_{A_n} \right).
\end{align*}
We define a residue mapping
\[ \varrho_{\uuB}^r : \Omega^{r}_{\uuB} \lra \displaystyle \varprojlim_m \ \Omega^{r-1}_{B/(q_0,p)^m}=\Omega^{r-1}_B \]
as the projective limit, with respect to $m$, of the composite map:
\[\xymatrix{ \Omega^{r}_{\uuB} \ar[rr]^-{\cano} && \Omega^{r}_{B\psfu/(q_0,p)^m} \ar[rr]^-{\Res_{u=0}} && \Omega^{r-1}_{B/(q_0,p)^m}. }\]
The map $\varrho_{\uuB}^r$ induces a map
\stepcounter{equation}
\begin{equation}\label{eq3-2-0} \varrho_{\uuB,n}^r : \Omega^r_{\uuB_n} \lra \Omega^{r-1}_{B_n}.\end{equation}
Inverting $q_0$, we get a map
\begin{equation}\label{eq3-2-1}  \varrho^r_n : \Omega^r_{\uuA_n} \lra \Omega^{r-1}_{A_n}. \end{equation}
The following lemma is straight-forward and left to the reader:
\addtocounter{thm}{2}
\begin{lem}
The following square is commutative$:$
\[ \xymatrix{ \Omega^r_{\uuA_n} \ar[rr]^-{f_2} \ar[d]_{\varrho^r_n} & & \Omega^r_{\uuA_n} \ar[d]^{\varrho^r_n} \\
 \Omega^{r-1}_{A_n} \ar[rr]^-{f_1} & & \Omega^{r-1}_{A_n}. } \]
\end{lem}
By this lemma, the maps $\varrho_n^\bullet$ induce a homomorphism of complexes
\[ \cS_n(2)_{\uuA,\uuA} \lra \cS_n(1)_{A,A}[-1] \quad \hbox{ for \, $n \ge 1$} \]
and a residue map of syntomic cohomology
\stepcounter{equation}
\begin{equation}\label{eq3-2-htsyn}
\htsyn : H^2_\syn(\uuA,\cS_{\bZ_p}(2)) \lra H^1_\syn(A,\cS_{\bZ_p}(1)). 
\end{equation}
We define the required arrow $\tsyn$ in the diagram \eqref{eq2-1} as the composite map
\[\xymatrix{ \tsyn : H^2_\syn(\cE(\cD), \cS_{\bZ_p}(2)) \ar[r]^-{\beta^*} & H^2_\syn(\uuA,\cS_{\bZ_p}(2)) \ar[r]^-{\htsyn} & H^1_\syn(A,\cS_{\bZ_p}(1)), }\] where $\beta^*$ denotes the pull-back by $\beta : \Spec(\uuB) \to \cE$ in \eqref{eq3-3}.
\begin{rem}
By \eqref{eq3-2-0}, we obtain a residue map
\[ \htsyn : H^2_\syn(\uuB,\cS_{\bZ_p}(2)) \lra H^1_\syn(B,\cS_{\bZ_p}(1)). \]
By a similar construction, we obtain a residue map
\[ \htsyn : H^2_\syn(\uuR,\cS_{\bZ_p}(2)) \lra H^1_\syn(R,\cS_{\bZ_p}(1)). \]
These maps will be used in {\rm\S\ref{sect3-5}} below.
\end{rem}
\subsection{\'Etale residue mapping}\label{sect3-3}
Let $R[u]_{(p)}$ be the localization of $R[u]$ at the prime ideal $(p)$, 
and denote its $p$-adic completion by $\uR$:
\[ \uR := \varprojlim_n \ R[u]_{(p)}/(p^n), \] 
which is a complete discrete valuation ring with residue field $k(u)$.
There is a natural embedding of $R[u,u^{-1}]$-algebras
\[\xymatrix{ \uuR \, \ar@{^{(}->}[r]  & \uR. }\]
Let $K$ be the fraction field of $R$, and put
\[ \uK := \uR \otimes_R K= \uR [p^{-1}]. \]
The following lemma will be used in Proposition \ref{prop3-3-1} below.
\begin{lem}\label{lem3-3-0}
For $n \ge 1$ and $f\in (\uR)^\times$ {\rm(}resp.\ $f\in (\uK)^\times${\rm)},
 there are $f_0\in R[u]_{(p)}^\times$ and $g\in \uR$
{\rm(}resp.\ $f_0\in K(u)^\times$ and $g\in \uK${\rm)} with $f=f_0(1+p^{n}g)$. Consequently the natural maps
\[ R[u]_{(p)}^\times/p^n\lra (\uR)^\times/p^n,\quad K(u)^\times/p^n\lra (\uK)^\times/p^n \]
are surjective for any $n \ge 1$.
\end{lem}
\begin{pf}
Exercise.
\end{pf}
We construct here an auxiliary residue mapping \[ \ttet : H^2_{\et}(\uK,\bZ/p^n(2)) \lra H^1_{\et}(K,\bZ/p^n(1)), \]
which will be useful in \S\ref{sect3-5} below. Let $\ol K$ be the algebraic closure of $K$, and put $\uoK:=\uK \otimes_K \ol K$. We have a $G_K$-equivariant homomorphism
\begin{equation}\label{eq3-3-1} (\uoK)^\times/p^n \lra \bZ/p^n,\qquad f \longmapsto \Res_{u=0}\, \frac{df}{f}, \end{equation}
where $\Res_{u=0}$ denotes a residue map defined in a similar way as for $\varrho_{\uuB,n}^1$ in \eqref{eq3-2-0}.
Since $H^1_\et(\uoK,\bZ/p^n(1)) \simeq (\uoK)^\times/p^n$, this map induces a $G_K$-equivariant homomorphism
\begin{equation}\label{eq3-3-2} H^1_\et(\uoK,\bZ/p^n(1)) \lra \bZ/p^n. \end{equation}
We define the desired map $\ttet$ by the composition
\[\begin{CD} \ttet : H^2_{\et}(\uK,\bZ/p^n(2)) \lra H^1_{\et}(K,H^1_\et(\uoK,\bZ/p^n(2))) @>{\eqref{eq3-3-2}}>> H^1_{\et}(K,\bZ/p^n(1)), \end{CD}\] where the first arrow is an edge homomorphism of the Hochschild-Serre spectral sequence
\[E_2^{a,b}=H^a_{\et}(K,H^b_\et(\uoK,\bZ/p^n(2))) \Lra H^{a+b}_{\et}(\uK,\bZ/p^n(2)) \]
and we have used the following fact:
\addtocounter{thm}{2}
\begin{prop}\label{prop3-3-1}
The natural restriction map \[ H^2_{\et}(\uK,\bZ/p^n(2)) \lra E_2^{0,2}=H^2_\et(\uoK,\bZ/p^n(2))^{G_K} \] is zero for any $n \ge 1$.
\end{prop}
\begin{pf}
We show that the natural map
\[ K^M_2(\uK)/p^n \lra K^M_2(\uoK)/p^n \]
is zero, which implies the assertion by Remark \ref{rem1-2-1}\,(2).
Let $K(u)$ be the rational fuction field in $u$ over $K$, which is a subfield of $\uK$.
By Lemma \ref{lem3-3-0} the natural map
\[K_2^M(K(u))/p^n \lra K_2^M(\uK)/p^n\]
is surjective.
Hence it is enough to show that the natural map
\[ K^M_2(K(u))/p^n \lra K^M_2({\ol K}(u))/p^n \]
is zero. One can see $K^M_2({\ol K}(u))/p^n=0$ in the following way.
Since \[ \{{\ol K}^\times,{\ol K}(u)^\times\}=0 \quad \hbox{ in \;\; $K^M_2({\ol K}(u))/p^n$}, \] it is enough to show that
$\{u-a,b-u\}=0$. If $a=b$, it is clear. If $a\neq b$,
one may replace $(u-a)$ with $(u-a)/(b-a)$
and $(b-u)$ with $(b-u)/(b-a)$. Then one has
\[ \left\{\frac{u-a}{b-a},\frac{b-u}{b-a}\right\}=0 \quad \hbox{ as \;\; $\dfrac{u-a}{b-a}+\dfrac{b-u}{b-a}=1$}, \]
which shows the assertion.
\end{pf}
\par
The following proposition will be used in \S\ref{sect3-5} below.
\begin{prop}\label{prop3-3-2}
Let $m$ be a positive integer, and let $s_m : B \to R$ be the homomorphism of $R$-algebras sending $q_0$ to $p^m$.
Consider the following cartesian diagram of schemes{\rm:}
\[ \xymatrix{
\Spec(\uR) \ar[r]^-{\gamma} & \Spec(\uuB_{(m)}) \ar[d] \ar[r]^-{\alpha_m} \ar@{}[rd]|{\square} & \cE_{(m)} \ar[r] \ar[d] \ar@{}[rd]|{\square} & \Spec(R) \ar[d]^{s_m} \\
&  \Spec(\uuB) \ar[r]^-{\beta} & \cE \ar[r] & \Spec(B),
} \]
where $\uuB_{(m)}$ and $\cE_{(m)}$ are defined by this diagram, and $\gamma$ is induced by the natural inclusion $\uuR \hra \uR \simeq \uuB_{(m)}$. Then the following diagram is commutative{\rm:}
\[ \xymatrix{H^2_{\et}(E_{(m)},\bZ/p^n(2)) \ar[rr]^-{\alpha_m^*} \ar[rd]_{\tet} && H^2_{\et}(\uK,\bZ/p^n(2)) \ar[ld]^-{\ttet} \\ & H^1_{\et}(K,\bZ/p^n(1)), } \]
where $\tet$ is defined in a similar way as for $\tet$ in the diagram \eqref{eq2-1}.
\end{prop}
\begin{pf}
Put $E_{\ol K}:=E_{(m)} \otimes_K \ol K$, which is indepenet of $m$.
By the construction of these maps, the assertion is reduced to the commutativity of a diagram of $G_K$-modules
\[ \xymatrix{H^1_{\et}(E_{\ol K},\bZ/p^n(1)) \ar[rr]^-{\gamma^*\circ\alpha_m^*} \ar[rd]_{\eqref{eq2-0}} && H^1_{\et}(\uoK,\bZ/p^n(1)) \ar[ld]^-{\eqref{eq3-3-2}} \\ & \bZ/p^n, } \] which has been shown in \cite{A0} Lemma 4.2.
\end{pf}

\subsection{Commutativity of the key diagram}\label{sect3-4}
We prove that the squares in \eqref{eq2-1} are commutative. The lower square of \eqref{eq2-1} commutes by the construction of $\tsyn$, which is rather straight-forward and left to the reader.
We prove the commutativity of the upper square of \eqref{eq2-1}, which will be finished in \S\ref{sect3-5} below.
\par
We first reduce the problem to the case that the positive integer $r$ fixed in \S\ref{sect2-1} is at least $3$.
Put $q_1:=(q_0)^{1/3}$, and let $E'$ be the Tate curve $E_{q,A[q_1]}$ over $A[q_1]=\psfo$.
Let $\cE':=\cE_{q,B[q_1]}$ be the regular proper model of $E'$ over $B':=B[q_1]=\psro$ defined in a similar way as for $\cE$ (see \S\ref{sect2-1}).  Let $\cD' \subset \cE'$ be the divisor defined by $q_1$. Consider the following diagram:
{\scriptsize 
\[\xymatrix{H^2_{\syn}(\cE(\cD),\cS_{\bZ_p}(2)) \ar[rd]^{\alpha} \ar[rrr]^-{\tsyn} \ar[ddd]_{c_2} \ar@{}[rrrd]|{(1)} 
 \ar@{}[rddd]|{(2)} &&& H^1_{\syn}(B,\cS_{\bZ_p}(1)) \ar[ld]_{\alpha} \ar[ddd]^{c^1} \ar@{}[lddd]|{(3)} \\
 & H^2_{\syn}(\cE'(\cD'),\cS_{\bZ_p}(2)) \ar[r]^-{\tsyn} \ar[d]_{c^2} & H^1_{\syn}(B',\cS_{\bZ_p}(1)) \ar[d]^{c^1} & \\
 & H^2_{\et}(\cE'[p^{-1}],\bZ_p(2)) \ar[r]^-{\tet} & H^1_{\et}(B'[p^{-1}],\bZ_p(1)) & \\
H^2_{\et}(\cE[p^{-1}],\bZ_p(2))  \ar[rrr]^-{\tet} \ar[ru]^{\alpha} \ar@{}[rrru]|{(4)} &&& H^1_{\et}(B[p^{-1}],\bZ_p(1)) \ar@{_{(}->}[ul]_{\beta} }\]
}Here the arrows $\alpha$ and $\beta$ are pull-back maps, and $\beta$ is injective by a standard norm argument. The squares (1) and (4) commute by the construction of $\tsyn$ and $\tet$, respectively. The squares (2) and (3) commute by the functoriality of $c^2$ and $c^1$, respectively (Theorem \ref{thm1-1-1}). Therefore the commutativity of the outer rectangle is reduced to that of the central square, i.e., we are reduced to the case that $r \ge 3$ (because we have $q=q_0^r=q_1^{3r}$ and $(p,3r)=1$ by the assumptions on $p$ and $r$).

We may thus assume $r \geq 3$.
Then we introduce an element $\xi_{\cE} \in H^2_{\syn}(\cE(\cD),\cS_{\bZ_p}(2))$ as follows.
Let $Z_i\subset \cE$ be the section defined by $u=q_0^i$.
Put \[ \cD_i:=Z_i\cap \cD \quad \hbox{ and } \quad Z:=\sum_{i=0}^{r-1} \ Z_i. \]
Noting $H^0_\syn(Z_i(\cD_i),\cS_{\bZ_p}(1))=0$,
we obtain, from Proposition \ref{prop1-1-1}, an exact sequence
\begin{equation}\label{eq2-0-0}
0 \to H^2_\syn(\cE(\cD),\cS_{\bZ_p}(2)) \to H^2_\syn(\cE(\cD+Z),\cS_{\bZ_p}(2))
\os{\partial}{\to} \bigoplus_{i=0}^{r-1} \ H^1_\syn(Z_i(\cD_i),\cS_{\bZ_p}(1)).
\end{equation}
Let $0<a<b<r$ be integers. Put
\[ f(u):=\frac{\theta(q_0^au)^{r}}{\theta(u)^{r-a}\theta(qu)^a}
=(-u)^a \frac{\theta(q_0^au)^r}{\theta(u)^r} \;\; \hbox{ and } \;\;
g(u):=\frac{\theta(q_0^bu)^{r}}{\theta(u)^{r-b}\theta(qu)^b}
=(-u)^b\frac{\theta(q_0^b u)^r}{\theta(u)^r},
\]
which are rational functions on $E=E_{q,A}$ by Proposition \ref{prop1-2}.
Put
\[ \xi'_\cE:=\left\{
\frac{f(u)}{f(q_0^{-b})},
\frac{g(u)}{g(q_0^{-a})}
\right\}^\syn=\left\{
\frac{f(u)}{f(q_0^{-b})}\right\}^\syn\cup
\left\{\frac{g(u)}{g(q_0^{-a})}
\right\}^\syn\in H^2_\syn(\cE(\cD+Z),\cS_{\bZ_p}(2)), \]
and the braces $\{ - \}^\syn$ denote the syntomic symbol in \eqref{eq1-2-3} and \eqref{eq1-2-4} with $U:=\cE-\cD-Z$.
% \[ \{ - \}^\syn : \vG(E,\cO_E^\times) \lra H^1_{\syn}(\cE(\cD),\cS_{\bZ_p}(1)). \]
One can easily check that the boundary $\partial(\xi'_\cE)$ agrees with
the syntomic symbol of the tame symbol of
\[ \left\{ \frac{f(u)}{f(q_0^{-b})}, \frac{g(u)}{g(q_0^{-a})} \right\} \in K^M_2(\cO(U)). \]
The tame symbol vanishes at each $Z_i-D_i$ by a straight-forward computation.
Therefore one has $\partial(\xi'_\cE)=0$, and $\xi'_\cE$ defines 
an element $\xi_{\cE}\in H^2_\syn(\cE(\cD),\cS_{\bZ_p}(2))$ 
by the exact sequence \eqref{eq2-0-0}.
\stepcounter{thm}
\begin{prop}\label{prop3-1}
$H^2_{\syn}(\cE(\cD),\cS_{\bZ_p}(2))$ is generated by the subgroups
\[ H^2_{\syn}(\cE,\cS_{\bZ_p}(2)), \;\;\;
 R_{\dR} \cup \{ q_0 \}^\syn \;\;\; \hbox{ and } \;\;\; \bZ_p \, \xi_{\cE}, \]
where $R_{\dR}$ denotes the image of the natural map
\[ \displaystyle R=\varprojlim_{n \ge 1} \ H^0_{\dR}(\cE_n) \lra H^1_{\syn}(\cE,\cS_{\bZ_p}(1)) \]
 induced by \eqref{eq1-1-1}, and the braces $\{ - \}^\syn$ denote the syntomic symbol in {\rm \eqref{eq1-2-3}:}
\[ \{ - \}^\syn : \vG(E,\cO_E^\times) \lra H^1_{\syn}(\cE(\cD),\cS_{\bZ_p}(1)). \]
\end{prop}
\par
This proposition will be proved in \S\ref{sect3-4'} below.
We explain how to prove the commutativity of \eqref{eq2-1}.
It is enough, by Proposition \ref{prop3-1}, to check the commutativity for the elements of
\[ R_{\dR} \cup \{ q_0\}^\syn, \;\;\; \bZ_p \, \xi_{\cE} \;\; \hbox{ and } \;\;
 H^2_{\syn}(\cE,\cS_{\bZ_p}(2)). \]
The commutativity for $R_{\dR} \cup \{ q_0\}^\syn$ is clear.
Indeed, we have
\[ \tsyn(R_\dR \cup \{ q_0 \}^\syn)=0  \;\; \hbox{ in } \;  H^1_{\syn}(A,\cS_{\bZ_p}(1)) \]
by the construction of $\tsyn$. On the other hand, we have
\[ \tet \circ c^2(R_\dR\cup \{ q_0 \}^\syn)=\tet (\{1+pR\}^\et \cup \{ q_0 \}^\et) = 0
 \;\; \hbox{ in } \; H^1_{\et}(L,\bZ_p(1)) \]
by Remark \ref{rem1-2-2} and Theorem \ref{prop1-2-1}.
We will show
\stepcounter{equation}
\begin{equation}\label{eq3-1-0} 
c^1 \circ \tsyn(\xi_{\cE})=
\tet \circ c^2(\xi_{\cE}) \quad
\mbox{ in }H^1_{\et}(L,\bZ_p(1)) \simeq L^\times/p^n
\end{equation}
in \S\ref{sect3-5-0} below.
As for the commutativity for the elements of $H^2_{\syn}(\cE,\cS_{\bZ_p}(2))$,
 it is equivalent to the commutativity of the following square:
\begin{equation}\label{eq3-4}
\xymatrix{
 H^2_{\syn}(\cE,\cS_{\bZ_p}(2)) \ar[rr]^-{\tsyn} \ar[d]_{c^2} & & H^1_{\syn}(B,\cS_{\bZ_p}(1)) \ar[d]^{c^1} \\
 H^2_{\et}(\cE[p^{-1}],\bZ_p(2)) \ar[rr]^-{\tet} & & H^1_{\et}(B[p^{-1}],\bZ_p(1)), }
\end{equation}
where the arrow $\tet$ is defined by the composite map of sheaves on $B[p^{-1}]_\et$
\[\xymatrix{ R^1f_{\cE*}\bZ/p^n(2) \ar[r] & j_*R^1f_*\bZ/p^n(2) \ar[rr]^-{\eqref{eq2-0}} &&
 j_*\bZ/p^n(1) = \bZ/p^n(1) \qquad(n \ge 1) }\]
and similar arguments as for $\tet$ in the diagram \eqref{eq2-1}.
Here $f_{\cE}$, $f$ and $j$ are as follows:
\[ \xymatrix{ E[p^{-1}] \; \ar@{^{(}->}[r] \ar[d]_{f} & \cE[p^{-1}]  \ar[d]^{f_\cE} \\
  \Spec(A[p^{-1}]) \; \ar@{^{(}->}[r]^-{j} & \Spec(B[p^{-1}]). }\]
We will prove the commutativity of \eqref{eq3-4} in \S\ref{sect3-5} below.
\subsection{Proof of Proposition \ref{prop3-1}}\label{sect3-4'}
Put
\[ \cE^{(0)}:=\cE,\qquad \cE^{(1)}:=\cD^{(1)}, \qquad \cE^{(2)}:=\cD^{(2)}, \]
where $\cD^{(1)}$ denotes the disjoint union of the irreducible components of $\cD$,
 and $\cD^{(2)}$ denotes the disjoint union of the intersections of two distinct irreducible components of $\cD$.
Since $\cD$ is the standard N\'eron $r$-gon over $\Spec(R)$, we have
\[ \cE^{(1)}=\coprod_{j=1}^{r} \ \bP^1_R \quad \hbox{ and } \quad \cE^{(2)}=\coprod_{j=1}^{r} \ \Spec(R). \]
Fix an integer $n \ge 1$.
One can easily check that $\cE$ and $\cE^{(m)}$\,($m=0,1,2$) satisfy
 the conditions (i)--(iv) in Proposition \ref{prop1-1-1}, using Lemma \ref{lem2-1-1} and Remark \ref{rem2-1-1}.
Hence there is a spectral sequence
\begin{equation}\label{eq3-4-0}
E_1^{a,b}=H^{2a+b}_{\syn}(\cE^{(-a)},\cS_n(a+2)) \Lra H^{a+b}_{\syn}(\cE(\cD),\cS_n(2)).
\end{equation}
Since we have
\begin{equation}\label{eq3-4-1}
 E_1^{a,b}=0  \;\; \hbox{ unless $-2 \le a \le 0$ and $2a+b \ge 0$}, 
\end{equation}
the quotient group $C_n := H^2_{\syn}(\cE(\cD),\cS_n(2))/E_\infty^{0,2}$ fits into an exact sequence
\[ 0 \lra E_2^{-1,3} \lra C_n \lra E_2^{-2,4} \lra E_2^{0,3}. \]
We will prove the following lemma:
\addtocounter{thm}{2}
\begin{lem}\label{lem3-4-1}
\begin{enumerate}
\item[{\rm (1)}]
Put $R_n:=R/(p^n)$. Then $E_2^{-1,3}$ agrees with the diagonal subgroup of
\[ E_1^{-1,3} = \bigoplus_{j=1}^{r} \ H^1_{\syn}(\bP^1_R,\cS_n(1)) \simeq \bigoplus_{j=1}^{r} \ R_n. \]
\item[{\rm (2)}]
$E_2^{-2,4}$ agrees with the diagonal subgroup of
\[ E_1^{-2,4} = \bigoplus_{j=1}^{r} \ H^0_{\syn}(R,\cS_n(0)) \simeq \bigoplus_{j=1}^{r} \ \bZ/p^n. \]
Moreover, the edge homomorphism
\[ H^2_{\syn}(\cE(\cD),\cS_n(2)) \lra E_1^{-2,4} \]
sends $\xi_{\cE}$ to $(r,r,\dotsc,r)$ up to a sign.
\end{enumerate}
\end{lem}
\par
We first prove Proposition \ref{prop3-1}, admitting this lemma. Consider the spectral sequence \eqref{eq3-4-0}. Since $H^0_\syn(\bP^1_R,\cS_n(1))=0$ by \eqref{eq1-1-2}, we have
\[\xymatrix{ E_1^{0,2} = E_2^{0,2} \ar@{=}[r]^-{\eqref{eq3-4-1}} & E_\infty^{0,2} }\]
and hence a short exact sequence
\[ 0 \lra H^2_{\syn}(\cE,\cS_n(2))) \lra H^2_{\syn}(\cE(\cD),\cS_n(2))) \lra C_n \lra 0. \]
Taking projective limit with respect to $n \ge 1$, we get an exact sequence
\[ 0 \lra H^2_{\syn}(\cE,\cS_{\bZ_p}(2))) \lra H^2_{\syn}(\cE(\cD),\cS_{\bZ_p}(2))) \lra \varprojlim_{n \ge 1} \, C_n. \]
By Lemma \ref{lem3-4-1}, there is a short exact sequence
\addtocounter{equation}{1}
\begin{equation}\label{eq3-4-2}
0 \lra R_n \lra C_n \lra \bZ/p^n \lra 0,
\end{equation}
which shows that $\varprojlim{}_{n \ge 1} \, C_n$ is generated by elements which lift to either
\[ R_\dR \cup \{ q_0 \}^\syn \;\; \hbox{ or } \;\;
 \bZ_p \,\xi_{\cE} \;\; \big(\subset  H^2_{\syn}(\cE(\cD),\cS_{\bZ_p}(2))\big). \]
Thus we obtain the proposition.
We prove Lemma \ref{lem3-4-1} in what follows.
\par \medskip \smallskip \noindent
{\it Proof of Lemma \ref{lem3-4-1}.}
The assertion on the the edge homomorphism $H^2_{\syn}(\cE(\cD),\cS_n(2)) \to E_1^{-2,4}$
 follows from computations on symbols (cf.\ \cite{A0} Remark 5.5).
% is straight-forward and left to the reader.
Thus $E_2^{-2,4}$ contains the diagonal subgroup of $E_1^{-2,4}$, by the assumption that $p$ is prime to $r$.
Similarly, it is easy to see that $E_2^{-1,3}$ contains the diagonal subgroup of $E_1^{-1,3}$.
It remains to show that $E_2^{-1,3}$ and $E_2^{-2,4}$ are contained in the diagonal subgroups of
 $E_1^{-1,3}$ and $E_1^{-2,4}$, respectively.

Recall that $\cD$ is the standard N\'eron $r$-gon over $\Spec(R)$.
Let $\cD_1,\dotsc,\cD_{r}$ be the irreducible components of $\cD$, which are all isomorphic to $\bP^1_R$.
Changing the ordering of these components if necessary, we suppose that $\cD_1$ meets $\cD_{r}$ and $\cD_2$,
 and that $\cD_j$ meets $\cD_{j-1}$ and $\cD_{j+1}$ for $j=2,\dotsc,r-1$
 (hence $\cD_{r}$ meets $\cD_{r-1}$ and $\cD_1$).
Put
\[ T_j := \cD_j \times_\cE \cD_{j+1} \quad(j=1,\dotsc,r-1),\qquad T_{r} := \cD_{r} \times_\cE \cD_1, \]
which are all isomorphic to $\Spec(R)$. Let \[\xymatrix{ i : \displaystyle \cE^{(1)}=\coprod_{j=1}^{r} \ \cD_j \ar[r] & \cE }\] be the natural finite morphism. Since $\cD_j \simeq \bP^1_R$ for each $1 \le j \le r$, we have isomorphisms
\begin{align*}
 H_{\syn}^1(\cD_j,\cS_n(1)) & \simeq H_{\dR}^0(\cD_{j,n}) \simeq R_n,  \\
 H_{\syn}^2(\cD_j,\cS_n(1)) & \simeq \ker(1-\sigma: R_n \to R_n)=\bZ/p^n, \\
 H_{\syn}^3(\cD_j,\cS_n(2)) & \simeq H_{\dR}^2(\cD_{j,n}) \simeq R_n  \qquad \hbox{(trace isomorphism)},
\end{align*}
by \eqref{eq1-1-2}. Lemma \ref{lem3-4-1} follows from the following claims (i) and (ii):
\begin{enumerate}
\item[(i)]
{\it The composite map
\[ E_1^{-1,3}=\bigoplus_{j=1}^{r} \ R_n = \bigoplus_{j=1}^{r} \ H_{\dR}^0(\cD_{j,n})
 \os{i_*}\lra H_{\dR}^2(\cE_n) \os{i^*}\lra \bigoplus_{j=1}^{r} \ H_{\dR}^2(\cD_{j,n}) \simeq \bigoplus_{j=1}^{r} \ R_n \]
is given by the $r \times r$ matrix
\[  \begin{pmatrix} -2 & \phami 1 & \phami 0 & \dotsb & \phami 0 & \phami 1 \\
 \phami 1 & -2  & \phami 1 & \ddots & \phami \vdots & \phami 0 \\
 \phami 0 & \phami 1  & -2 & \ddots & \phami 0 & \phami \vdots \\
 \phami \vdots & \ddots  & \ddots & \ddots & \phami 1 & \phami 0 \\
 \phantom{\vdots}\;\, 0 & \dotsb & \phami 0 & \phami 1 & -2 & \phami 1 \\
 \phantom{\vdots}\;\, 1 & \phami 0  & \dotsb & \phami 0 & \phami 1 & -2 \\
 \end{pmatrix} \;\; .  \]
Consequently, $E_2^{-1,3}$ is contained in the diagonal subgroup of $E_1^{-1,3}$
 by the assumption that $r$ is prime to $p$
{\rm (}note also that the above composite map agrees with the composite map
$E_1^{-1,3} \to E_1^{0,3} \to \bigoplus_{j=1}^{r} \ H^3_{\syn}(\cD_j,\cS_n(2))\simeq \bigoplus_{j=1}^{r} \ R_n${\rm )}.
}
\item[(ii)]
{\it The edge homomorphism
\[E_1^{-2,4}=\bigoplus_{j=1}^r \ \bZ/p^n[T_j] \lra E_1^{-1,4}=\bigoplus_{j=1}^r \ H_{\syn}^2(\cD_j,\cS_n(1))
 \simeq \bZ/p^n[D_j] \]
sends $[T_j]$ to
\[ \begin{cases} [D_j]-[D_{j+1}] & \quad (1 \le j \le r-1) \\ [D_1]-[D_{r}] & \quad (j = r), \end{cases} \]
where the signs arise from the construction of the spectral sequence \eqref{eq3-4-0}.
Consequently, $E_2^{-2,4}$ is contained in the diagonal subgroup of $E_1^{-2,4}$. }
\end{enumerate}
These claims follows from standard facts on intersection theory of divisors.
The details are straight-forward and left to the reader. % (cf.\ \cite{A} (4.17), (4.21), (4.22)).
This completes the proof of Lemma \ref{lem3-4-1} and Proposition \ref{prop3-1}.
\subsection{Proof of \eqref{eq3-1-0}}\label{sect3-5-0}
This subsection is devoted to showing \eqref{eq3-1-0}. More precisely, we will prove
\begin{align*} 
c^1 \circ \tsyn(\xi_{\cE})&=
(-1)^{a(r-b)}
q_0^{a(b-a)(b-r)}
\left(
\frac{\theta(q_0^{b})^{b}}{\theta(q_0^{b-a})^{b-a}\theta(q_0^{a})^{a}}
\right)^r
\left(
\frac{S(q_0^{b})}{S(q_0^{b-a})S(q_0^{a})}
\right)^{r^2}\\
& =
\tet \circ c^2(\xi_{\cE}) 
\end{align*}
in $L^\times/p^n$, where we put
\[S(\alpha):=\prod_{k=1}^\infty
\left(\frac{1-\alpha q^k}{1-\alpha^{-1} q^k}\right)^k.\]
Since the equality for $\tet$ has been shown in \cite{A0} Lemma 7.4,
we prove the equality for $\tsyn$.
It follows from the construction of the syntomic residue
map (\S \ref{sect3-2}) that we may truncate $\theta(u)$
with respect to $q_0$.
Therefore we can calculate $\tsyn$ by factorizing the infinite products.
It is enough to check the following:
{\allowdisplaybreaks
\begin{align}
\label{tau0}
&\htsyn\{a,b\}=1,&a,~b\in A^\times\\
\label{tau1}
&\htsyn\{a,u\}=a, & a\in A^\times\\
\label{tau2}
&\htsyn\{a,g(u)\}=1, &a\in A^\times,\;g(u)\in R\psrqu^\times\\
\label{tau4}
&\htsyn\{u,g(u)\}=g(0)^{-1}, & g(u)\in R\psrqu^\times\\
\label{tau5}
&\htsyn\{u,1-bu^{-1}\}=1, & b\in(p,q_0)\\
\label{tau6}
&\htsyn\{g(u),h(u)\}=1, & g(u),h(u)\in R\psrqu^\times\\
\label{tau7}
&\htsyn\{g(u),1-bu^{-1}\}=g(0)^{-1}g(b), & 
g(u)\in R\psrqu^\times,\;b\in(p,q_0)\\
\label{tau8}
&\htsyn\{1-bu^{-1},1-cu^{-1}\}=1, & 
~b,c\in(p,q_0)
\end{align}
}where $\htsyn$ is the residue map \eqref{eq3-2-htsyn}.
We prove only \eqref{tau4} and \eqref{tau7}. The other equalities are simpler and left to the reader.
By the definition of syntomic symbols (cf.\,\S\ref{sect1-1'}), we have
\[ \{u\}^\syn=\left(\frac{du}{u},0\right),\quad \{g(u)\}^\syn=\left(\frac{dg(u)}{g(u)},
\frac{\,1\,}{p}\log\frac{\varphi(g(u))}{g(u)^p}\right). \]
By the definition of the product of syntomic cohomology (cf.\,Definition \ref{def1-2-2}), we have
\begin{align*}
\{u, g(u)\}^\syn&=\left(\frac{du}{u}\frac{dg(u)}{g(u)},
-\frac{\,1\,}{p}\log\frac{\varphi(g(u))}{g(u)^p}\frac{du}{u}\right)\\
&=\left(-\frac{dg(u)}{g(u)}\frac{du}{u},
-\frac{\,1\,}{p}\log\frac{\varphi(g(u))}{g(u)^p}\frac{du}{u}\right).
\end{align*}
Therefore
\[ \hat{\tau}^\syn_\infty[u, g(u)]=\left(-\frac{dg(0)}{g(0)},-\frac{\,1\,}{p}\log\frac{\varphi(g(0))}{g(0)^p}\right) =-[g(0)], \]
which completes the proof of \eqref{tau4}.
One can reduce \eqref{tau7} to \eqref{tau4} in the following way.
It is enough to show $\htsyn\{g(u),u-b\}=g(b)$.
Consider an endomorphism $w:\uuA\to\uuA$ of topological ring
such that $w(u)=u+b$. Then we have $\htsyn \circ w=\htsyn$ as 
$\varrho^r_n \circ w=\varrho^r_n$ (cf.\,\eqref{eq3-2-1}).
Thus one has \[ \htsyn\{g(u),u-b\}=\htsyn w\{g(u),u-b\}=\htsyn\{g(u+b),u\}=g(b). \]
as required.
\subsection{Commutativity of \eqref{eq3-4}}\label{sect3-5}
We prove the commutativity of \eqref{eq3-4}.
Let $\uR$ and $\uK$ be as we defined in \S\ref{sect3-3}.
Since $B$ is a regular local ring, $B[p^{-1}]$ is a UFD. By Remark \ref{rem1-2-1}\,(1), we have
\[ H^1_{\et}(B[p^{-1}],\bZ_p(1)) \simeq \varprojlim_{n} \ B[p^{-1}]^\times/p^n. \]
\begin{lem}\label{lem3-2}
For $m \ge 1$, let $s_m: B \to R$ be the $R$-homomorphism sending $q_0$ to $p^m$.
Then the map
\[ \prod_{m \ge 1} s_m \; : \;
 \varprojlim_{n} \; B[p^{-1}]^\times/p^n \lra  \prod_{m \ge 1} \; \varprojlim_{n}\; K^\times/p^n \]
 is injective.
\end{lem}
\begin{pf}
Exercise (left to the reader).
\end{pf}
We consider the following cartesian squares of schemes for each $m \ge 1$:
\[ \xymatrix{
\Spec(\uR) \ar[r]^-{\gamma} & \Spec(\uuB_{(m)}) \ar[d]_{\sigma_m} \ar[r]^-{\alpha_m} \ar@{}[rd]|{\square} & \cE_{(m)} \ar[r] \ar[d]_{t_m} \ar@{}[rd]|{\square} & \Spec(R) \ar[d]_{s_m} \\
&  \Spec(\uuB) \ar[r]^{\beta} & \cE \ar[r] & \Spec(B),
} \]
where $\uuB_{(m)}$ and $\cE_{(m)}$ are defined by this diagram, and $\gamma$ is induced by the natural inclusion $\uuR \hra \uR \simeq \uuB_{(m)}$. See \S\ref{sect3-1} for $\uuB$ and $\beta$. We defined $\cE_{(m)}$ (resp.\ $\Spec(\uuB_{(m)})$) by the lower (resp.\ upper) cartesian square. By Lemma \ref{lem3-2}, the commutativity of \eqref{eq3-4} is reduced to showing
 that of the following diagram for all $m \ge 1$:
\stepcounter{equation}
\begin{equation}\label{eq3-5}
\xymatrix{
H^2_{\syn}(\cE,\cS_{\bZ_p}(2)) \ar[rr]^-{t_m^* \circ c^2} \ar[d]_{\tsyn}
 & & H^2_{\et}(\cE_{(m)}[p^{-1}],\bZ_p(2)) \ar[d]^{\tet} \\
  H^1_{\syn}(B,\cS_{\bZ_p}(1)) \ar[rr]^-{s_m^* \circ c^1} & & H^1_{\et}(K,\bZ_p(1)). }
\end{equation}
\stepcounter{thm}
\begin{lem}\label{lem3-3}
The diagram \eqref{eq3-5} factors as follows$:$
\[\xymatrix{ H^2_{\syn}(\cE,\cS_{\bZ_p}(2)) \ar[r]^-{c^2} \ar[d]_{\beta^*} \ar@{}[rd]|{\bigcirc\;\;}
 & H^2_{\et}(\cE[p^{-1}],\bZ_p(2)) \ar[r]^-{t_m^*}  \ar[d]_{\beta^*} \ar@{}[rd]|{\bigcirc\;\;}
 & H^2_{\et}(\cE_{(m)}[p^{-1}],\bZ_p(2)) \ar[d]^{\alpha_m^*} \\
H^2_{\syn}(\uuB,\cS_{\bZ_p}(2)) \ar[r]^-{c^2} \ar[d]_{\htsyn} & H^2_{\et}(\uuB[p^{-1}],\bZ_p(2)) \ar[r]^-{\sigma_m^*} &
H^2_{\et}(\uuB_{(m)}[p^{-1}],\bZ_p(2)) \ar[d]^{\htet} \\
H^1_{\syn}(B,\cS_{\bZ_p}(1)) \ar[rr]^-{s_m^* \circ c^1} & & H^1_{\et}(K,\bZ_p(1)), } \]
and the upper squares are commutative. Here $\htsyn$ denotes the residue map constructed in \S\ref{sect3-2} and $\htet$ is defined as the projective limit, with respect to $n \ge 1$, of the composite map
\[\xymatrix{ \htet :  H^2_{\et}(\uuB_{(m)}[p^{-1}],\bZ/p^n(2)) \ar[r]^-{\gamma^*}&  H^2_{\et}(\uK,\bZ/p^n(2))
 \ar[r]^-{\ttet} & H^1_{\et}(K,\bZ/p^n(1)). }\]
See {\rm\S\ref{sect3-3}} for $\ttet$.
\end{lem}
\begin{pf}
The factorization $\tsyn=\htsyn \circ \beta^*$ follows from the construction of these maps,
 and the factorization $\tet=\htet \circ \alpha_m^*$ follows from Proposition \ref{prop3-3-2}.
The second assertion follows from the functoriality of $c^2$ and that of \'etale cohomology.
\end{pf}
\begin{lem}\label{lem3-4}
In the diagram of Lemma {\rm\ref{lem3-3}}, the lower square factors as follows{\rm:}
\[ \xymatrix{
H^2_{\syn}(\uuB,\cS_{\bZ_p}(2)) \ar[r]^-{\sigma^*_m} \ar[d]_{\htsyn} \ar@{}[rd]|{\bigcirc \;\;} &
H^2_{\syn}(\uuB_{(m)},\cS_{\bZ_p}(2)) \ar[r]^-{c^2} \ar[d]_{\htsyn} &
H^2_{\et}(\uuB_{(m)}[p^{-1}],\bZ_p(2)) \ar[d]^{\htet} \\
H^1_{\syn}(B,\cS_{\bZ_p}(1)) \ar[r]^-{s^*_m} & H^1_{\syn}(R,\cS_{\bZ_p}(1)) \ar[r]^-{c^1} & H^1_{\et}(K,\bZ_p(1)),
} \]
and the left square commutes.
\end{lem}
\begin{pf}
The first assertion follows from the functoriality of $c^2$ and $c^1$.
The commutativity of the left square follows from the construction of $\htsyn$'s (see \S\ref{sect3-2}).
\end{pf}
\begin{lem}\label{lem3-5}
In the diagram of Lemma {\rm \ref{lem3-4}}, the right square factors as follows{\rm:}
\[ \xymatrix{
H^2_{\syn}(\uuB_{(m)},\cS_{\bZ_p}(2)) \ar[rr]^-{c^2} \ar[d]_{\gamma^*} \ar@{}[rrd]|{\bigcirc \;\;} & &
 H^2_{\et}(\uuB_{(m)}[p^{-1}],\bZ_p(2)) \ar[d]^{\gamma^*} \\
H^2_{\syn}(\uR,\cS_{\bZ_p}(2)) \ar[rr]^-{c^2} \ar[d]_{\ttsyn} & &
 H^2_{\et}(\uK,\bZ_p(2)) \ar[d]^{\ttet} \\
H^1_{\syn}(R,\cS_{\bZ_p}(1)) \ar[rr]^-{c^1}  & & H^1_{\et}(K,\bZ_p(1)),
} \]
and the upper square commutes.
Here $\ttsyn$ denotes the syntomic residue map defined in a similar way as for $\htsyn$.
Note also that $\uK=\uR[p^{-1}]$.
\end{lem}
\begin{pf}
The first assertion follows from the construction of $\htsyn$ and $\htet$.
The commutativity of the upper square follows from the functoriality of $c^2$ (Theorem \ref{thm1-1-1}).
\end{pf}
By Lemmas \ref{lem3-3}--\ref{lem3-5}, the commutativity of \eqref{eq3-5}
 is reduced to that of the lower square of the diagram in  Lemma \ref{lem3-5},
  which is further reduced to that of the following square for all $n \ge 1$:
\addtocounter{equation}{3}
\begin{equation}\label{eq3-6}
\xymatrix{
H^2_{\syn}(\uR,\cS_n(2)) \ar[rr]^-{c^2} \ar[d]_{\ttsyn} & &
 H^2_{\et}(\uK,\bZ/p^n(2)) \ar[d]^{\ttet} \\
H^1_{\syn}(R,\cS_n(1)) \ar[rr]^-{c^1} & & H^1_{\et}(K,\bZ/p^n(1)).
} \end{equation}
By Theorem \ref{prop1-2-1} and Theorem \ref{lem1-2-1}, 
the diagram \eqref{eq3-6} is written as
\begin{equation}\label{eq3-7}
\xymatrix{
K_2^M(\uR)/p^n \ar[rr] \ar[d]_{\tau^\syn} & &
 K_2^M(\uK)/p^n \ar[d]^{\tau^\et} \\
R^\times/p^n \ar[rr] & & K^\times/p^n.
} \end{equation}
where the horizontal arrows are the natural maps
and $\tau^\syn$ (resp. $\tau^\et$) is induced
from $\ttsyn$ (resp. $\ttet$).
Thus the commutativity of \eqref{eq3-6} is reduced to
the explicit calculations of $\tau^\syn$ and $\tau^\et$.
By Lemma \ref{lem3-3-0} we may replace $K_2^M(\uR)/p^n$ with $K_2^M(R[u]_{(p)})/p^n$ and $K_2^M(\uK)/p^n$ with $K_2^M(K(u))/p^n$. Then the explicit formula of $\tau^\et$ has been shown in \cite{A0} Theorem 4.4. Therefore it is enough to show that the formula of $\tau^\syn$ agrees with it. However this follows from \eqref{tau0}--\eqref{tau8}.
This finishes the proof of the commutativity of the diagram \eqref{eq3-7} and hence the diagrams \eqref{eq3-6}, \eqref{eq3-4} and \eqref{eq2-1}. This completes the proof of Theorem \ref{thm2-1}.

\newpage
\section{Main result on Elliptic surfaces over $\boldsymbol{p}$-adic fields}
\label{knownsect}
We mean by an {\it elliptic surface} over a commutative ring $A$ a projective
flat morphism $\pi:X\to C$ with a section $e:C\to X$
such that $X$ and $C$ are projective smooth over $\Spec(A)$ of relative dimension $2$ and $1$, respectively
and such that the general fiber of $\pi$ is an elliptic curve.

For a field $F$, we denote the absolute Galois group ${\mathrm{Gal}}(\ol{F}/F)$ by $G_F$.
We often write $H^*(F,-)$ for the continuous Galois cohomology $H^*_\cont(G_F,-)$.

\subsection{Split multiplicative fiber}\label{ellsurfsect}
Let \[ \tT_\Z:=\coprod\P^1_\Z \] be the disjoint union of
copies of $\P^1_\Z$ indexed by $\Z/m$.
Gluing the section $0 : \Spec(\bZ) \hra \P^1_\Z$ of the $i$-th $\P^1_\Z$ with
$\infty: \Spec(\bZ) \hra \P^1_\Z$ of the $(i+1)$-st $\P^1_\Z$, we obtain a connected  proper curve
$T_\Z$ over $\Spec\Z$ whose normalization is $\tT_\Z\to T_\Z$.
We call $T_\Z\otimes_\Z A$ the {\it N\'eron polygon} or the {\it N\'eron $m$-gon} over a ring $A$ (cf.\ \cite{DR} II.1.1).
%($m$ indicates the number of irreducible components of $T_\Z$)
Let $\pi:X\to C$ be an elliptic surface over a ring $A$.
For an $A$-valued point $P\in C(A)$,
we call the fiber $\pi^{-1}(P)$ {\it split multiplicative of type $I_m$} over $A$,
if there is a closed subscheme 
$D^\dag\subset \pi^{-1}(P)$ which is isomorphic to
the N\'eron $m$-gon $T_\Z\times_\Z A$.
We call $\pi^{-1}(P)$ {\it multiplicative}, if $\pi^{-1}(P)\otimes_A\ol{k}$ is split multiplicative over $\ol{k}$
for any geometric point $\Spec(\ol k)\to\Spec (A)$.
A singular fiber which is not multiplicative is called an {\it additive} fiber.
When $A=\ol{F}$ is an algebraically closed field,
all the singular fibers are classified by Kodaira and N\'eron (cf.\ \cite{Si} IV \S8).
In particular, we note
\[ \pi^{-1}(P)\mbox{ is multiplicative }\Longleftrightarrow
~H^1_\et(\pi^{-1}(P),\zp)\neq0~
(\Longrightarrow~H^1_\et(\pi^{-1}(P),\zp)=\zp). \]
\subsection{Formal Eisenstein seires}
\label{sect5-2}
Let $p\geq 5$ be a prime number.
Let $K$ be a finite {\it unramified} extension of $\qp$ and let $R$ be its integer ring.
Let $\pi_R:\XR\to \CR$ be an elliptic surface over $R$.
Let $\varSigma=\{P_1,\cdots,P_s\}$ be the set of all $R$-rational points
 for which the fiber $\DR_i:=\pi^{-1}(P_i)$ are split multiplicative fibers of type $I_{r_i}$.
Put $\DR:=\sum_{i=1}^s \DR_i$ and $\UR:=\XR-\DR$.
We assume  \[ p\not| \,6r_1\cdots r_s. \]
Let $t_i$ be a uniformizer of $\O_{\CR,P_i}$.
Let $\iota_i:\Spec (R\psfti) \to \CR-\{P_i\}$ be the punctured neighborhood of $P_i$, and let $X_i$ be the fiber: 
\[ \xymatrix{ X_i \ar[r] \ar[d] \ar@{}[rd]|\square & \UR \ar[d]_{\pi_R} \\
\Spec (R\psfti) \ar[r]^{\iota_i} & \CR-\{P_i\}. } \]
Then $X_i$ is isomorphic to a Tate elliptic curve over $R\psfti$.
More explicitly, let $q\in R\psfti$ be the unique power series such that $\ord_{t_i}(q)=r_i$ and 
\begin{equation}\label{jtate}
j(X_i)=\frac{\,1\,}{q}+744+196884q+\cdots.
\end{equation}
Then there is an isomorphism of $R\psfti$-schemes
\begin{equation}\label{tateisom0}
X_i \simeq E_q,
\end{equation}
which is unique up to the translation and the involution $u\mapsto u^{-1}$ (cf.\ \cite{DR} VII.2.6).
Put \[ a_i:=t_i^{r_i}\cdot j(X_i)\big\vert_{t_i=0}= \frac{\,t_i^{r_i}\,}{q} \Big|_{t_i=0}\in R^\times. \]
Put $K_i:=K(\hspace{-4pt}\sqrt[r_i]{a_i})$, and let $R_i$ be its integer ring.
Note that $K_i$ is also unramified over $\qp$ as $r_i$ is prime to $p$.
There is $q_i\in t_i\cdot R_i\psrti^\times$ such that $q_i^{r_i}=q$, and we have $R_i\psfti= R_i\psfqi$.
Let $\kappa_i$ be the composition of natural maps
\begin{multline}\label{drlocal1}
\kappa_i : \vg(\XR,\Omega^2_{\XR/R}(\log \DR))\lra \vg(X_i,\Omega^2_{X_i/R})\simeq \vg(E_q,\Omega^2_{E_q/R})\\
\lra R_i \psfqi \frac{dq_i}{q_i}\frac{dx}{2y+x}\simeq R_i\psfqi,
\end{multline}
where the isomorphism in the middle is induced by \eqref{tateisom0}.
We define a $\zp$-submodule $\Eis(\XR,\DR)_\zp$ of $\vg(\XR,\Omega^2_{\XR/R}(\log \DR))$ by
\[ \omega \in \Eis(\XR,\DR)_\zp
\Longleftrightarrow
\kappa_i(\omega) \,\hbox{{\it is a formal power series of Eisenstein type for each $i$} (\S\ref{sect2-2})}\]
and call it the space of {\it formal Eisenstein series}.
As an immediate consequence of
the main result on Tate curve (Theorem \ref{thm2-1}), the image of the syntomic cohomolgy
group is contained in the space of formal Eisenstein series:
\begin{equation}\label{key2eq}
H^2_\syn(\XR(\DR),\cS_\zp(2))
\lra \Eis(\XR,\DR)_\zp\subset \vg(\XR,\Omega^2_{\XR/R}(\log \DR)).
\end{equation}
We also introduce a $\zp$-submodule $\Eis^{(n)}(\XR,\DR)_\zp
\subset \vg(\XR,\Omega^2_{\XR/R}(\log \DR))$ such that
$\omega\in \Eis^{(n)}(\XR,\DR)_\zp$ if and only if
each $\kappa_i(\omega)$\,($1 \le i \le s$) satisfies (E1) (\S\ref{sect2-2}) and 
\par\medskip\noindent
\;\;(E2)${}^{(n)}$ \;\;
 $a^{(j)}_k\in k^2\zp$ for all $j$ and $1\leq k \leq n$.
\par\medskip\noindent
Obviously one has
\[ \Eis(\XR,\DR)_\zp=\bigcap_{n \ge 1} \ \Eis^{(n)}(\XR,\DR)_\zp. \]

\subsection{Main result on elliptic surfaces}\label{sect5-3}
The boundary maps on \'etale, syntomic and de Rham cohomology induce a commutative diagram
\begin{equation}\label{mc0}
\xymatrix{
H^2_\et(\UbK,\zp(2))^{G_K} \ar[rr]^-{\partial_\et} && \bigoplus_{i=1}^s \ \zp[D_i] \ar@{=}[d]\\
H^2_\syn(\XR(\DR),\cS_\zp(2)) \ar[rr]^-{\partial_\syn} \ar[u]^a 
\ar[d]_{\eqref{key2eq}} && \bigoplus_{i=1}^s \ \zp[\DR_i] \ar@{=}[d] \phantom{\Big|}\\
\Eis(\XR,\DR)_\zp \ar[rr]^-{\partial_\dR} && \bigoplus_{i=1}^s \ \zp[\DR_i]. }
\end{equation}
Here $\partial_\dR$ denotes the the Poincare residue map.
$\partial_\et$ denotes the composite map
\[ H^2_\et(\UbK,\zp(2))\lra H^3_{\DbK}(\XbK,\zp(2)) \simeq \bigoplus_{i=1}^s \ \zp[\DR_i] \]
where the first arrow is the boundary localization sequence of the etale cohomology and
the last isomorphism is obtained from the Poincare-Lefschetz duality. 
See Proposition \ref{prop1-1-1} for $\partial_\syn$.
Note that $\partial_\et$ is injective modulo torsion.

\stepcounter{thm}
\begin{lem}\label{lemkey2}
We have
\stepcounter{equation}
\begin{equation}\label{key3eq}
\partial_\et\big( H^2_\et(\UbK,\zp(2))^{G_K}\big)\subset 
\partial_\dR\big(\Eis(\XR,\DR)_\zp\big).
\end{equation}
\end{lem}
\begin{pf}
It is enough to show that the map $a$ in \eqref{mc0} is surjective.
Put $\SR:=\CR-\varSigma$.
Let $\UF$ and $\SF$ denote the fibers of $\UR \to \Spec(R)$ and $\SR \to \Spec(R)$ over the closed point $s$ of $\Spec(R)$.
By a result of Tsuji \cite{Ts2} Theorem 5.1, there is an exact sequence
\begin{equation}\label{lemkey2eq1}
H^2_\syn(\XR(\DR),\cS_\zp(2))
\to H^2_\et(\UK,\zp(2))
\os{\delta_U}{\to} \plim{n}\O(\UF)^\times/p^n\simeq \plim{n}\O(\SF)^\times/p^n,
\end{equation}
where the last isomorphism follows from the fact that $\UF\to \SF$ is proper with geometrically connected fibers.
On the other hand, we assert that there is an exact sequence
\begin{equation}\label{lemkey2eq2}
0\lra H^2_\et(\SK,\zp(2))\os{\pi_K^*}{\lra}
H^2_\et(\UK,\zp(2))\lra
H^2_\et(\UbK,\zp(2))^{G_K}\lra0,
\end{equation}
where we put $\SK:=\SR \otimes_R K$. Indeed, the Hochschild-Serre spectral sequences
\begin{align}
\label{HS1}
E_{2,U}^{p,q} &=H^p(K,H^q_\et(\UbK,\zp(2)))\Longrightarrow
E_U^{p+q}=H^{p+q}(\UK,\zp(2)) \\
\label{HS2}
E_{2,S}^{p,q} &=H^p(K,H^q_\et(\SbK,\zp(2)))\Longrightarrow
E_S^{p+q}=H^{p+q}(\SK,\zp(2))
\end{align}
and the isomorphisms $H^q_\et(\SbK,\zp(2)) \simeq H^q_\et(\UbK,\zp(2))$ for $q\le 1$
 (cf.\ Lemma \ref{aux1} below) yield a short exact sequence
\[ 0 \lra H^2_\et(\SK,\zp(2))\os{\pi_K^*}{\lra} H^2_\et(\UK,\zp(2)) \lra E_{\infty,U}^{0,2}\lra0. \]
It remains to show $E_{2,U}^{0,2}=E_{\infty,U}^{0,2}$.
Since $\cd(K)=2$, it is enough to show that the edge homomorphism $E_{2,U}^{0,2} \to E_{2,U}^{2,1}$ is zero.
There is a commutative diagram with exact rows
\[ \xymatrix{
E_{2,U}^{0,2} \ar[r]^\alpha & E_{2,U}^{2,1} \ar[r]^\beta & E_{U}^3\\
% @AAA@A{\simeq}AA@AA{\bigcup}A\\
E_{2,S}^{0,2} \ar[r]^\gamma \ar[u] & E_{2,S}^{2,1} \ar[r]^\delta \ar[u]^{\wr\hspace{-1pt}} & E_{S}^3. \ar@{^{(}->}[u]
}\]
Therefore $\alpha=0$ $\Longleftrightarrow$ $\beta$ is injective $\Longleftrightarrow$ $\delta$ is injective
$\Longleftrightarrow$ $\gamma=0$. The last condition follows from $E_{2,S}^{0,2}=0$, which completes the proof of
 the exact sequence \eqref{lemkey2eq2}.

Finally we show that $a$ is surjective.
By \eqref{lemkey2eq1} and \eqref{lemkey2eq2},
it is enough to show
\[ \Image(H^2_\et(\UK,\Z/p^n(2))
\os{\delta_U}{\to} \O(\UF)^\times/p^n)=
\Image(H^2_\et(\SK,\Z/p^n(2))
\os{\delta_S}{\to} \O(\SF)^\times/p^n). \]
The inclusion `$\supset$' follows from the fact $\pi^*_s \circ \delta_U=\delta_S\circ \pi^*_K$, and
`$\subset$' follows from the fact $e^*_s\circ \delta_U=\delta_S \circ e^*_K$, where $e_R:\CR\to \XR$ denotes the section,
 and $\pi_F$ and $e_F$\,($F=K,s$) denote $\pi_R \otimes_R F$ and $e_R \otimes_R F$, respectively.
\end{pf}

Since $\partial_\et$ is injective modulo torsion, the right hand side of \eqref{key3eq}
gives an upper bound of the rank of $H^2_\et(\UbK,\zp(2))^{G_K}$.
However it is in general impossible to compute it
even though one can compute $\Eis^{(n)}(\XR,\DR)_\zp$ for a particular $n$.
The main result allows us to obtain an ``computable" upper bound under some mild conditions.
\addtocounter{thm}{5}
\begin{thm}\label{key3}
Let $\pi_R : \XR\to \CR$ and $\DR$ be as before.
Let \[ f_p:\bigoplus_{i=1}^s \ \zp [\DR_i]\lra \bigoplus_{i=1}^s \ \F_p [\DR_i] \] be 
the natural map taking the residue class modulo $p$, and let $\ol{\partial_\dR}$ be the composit map
\[ \Eis(\XR,\DR)_\zp\os{\partial_\dR}{\lra} 
\bigoplus_{i=1}^s \ \zp [\DR_i]\os{f_p}{\lra} 
\bigoplus_{i=1}^s \ \F_p [\DR_i]. \]
Assume $p\hspace{-3pt}\not|\,6r_1\cdots r_s$ and further the following conditions{\rm:}
\begin{description}
\item[(A)]
$H^3_\et(\XbK,\zp)$ is torsion-free.
\item[(B)]
$H^2_\et(\XbK,\Z/p(2))^{G_K}=0$.
\end{description}
Then we have
\stepcounter{equation}
\begin{equation}\label{conj1seq3}
\dim_\qp \ H^2_\et(\UbK,\qp(2))^{G_K} \leq \dim_{\F_p} \ \ol{\partial_\dR} \big( \Eis(\XR,\DR)_\zp\big).
%\leq \dim_{\F_p} \ \ol{\partial_\dR} \big( \Eis^{(n)}(\XR,\DR)_\zp\big)
\end{equation}
%for any $n\geq1$(the last term is a computable upper bound).
More precisely, for $v\in \bigoplus_i \ \zp[\DR_i]$
\begin{equation}\label{conj1seq4}
f_p(v)\not\in  \ol{\partial_\dR} \big( \Eis(\XR,\DR)_\zp\big)
\Longrightarrow
v\not\in \partial_\et\big( H^2_\et(\UbK,\qp(2))^{G_K}\big).
\end{equation}
\end{thm}
\begin{pf}
Applying $f_p$ on \eqref{key3eq}, we obtain
\[ \ol{\partial_\et} \big(H^2_\et(\UbK,\zp(2))^{G_K}\big) \subset \ol{\partial_\dR} \big( \Eis(\XR,\DR)_\zp\big). \]
Therefore all we have to do is to show that
the quotient group
\begin{equation}\label{tor-free00}
\zp^{\op s}\big/\partial_\et\big( H^2_\et(\UbK,\zp(2))^{G_K}\big)
\end{equation}
is torsion-free (note $ H^2_\et(\UbK,\qp(2))^{G_K}
\cong \partial_\et\big( H^2_\et(\UbK,\qp(2))^{G_K}\big)$).

Put $V_\zp:= H^2_\et(\XbK,\zp(2))/\langle D\rangle\ot\zp(1)$
where $\langle D\rangle$ denotes the subgroup of
the cycle classes of irreducible components of $\DbK$.
There is an exact sequence
\begin{equation}\label{exact0}
0\lra V_\zp
\lra H^2_\et(\UbK,\zp(2))\os{\partial_\et}{\lra}
\zp^{\op s}\lra 0
\end{equation}
under the condition {\bf (A)} by Lemma \ref{aux5} below.
Moreover $H^2_\et(\XbK,\zp(2))$ is torsion-free by Lemma \ref{aux2} below and the natural surjective map
\[ H^2_\et(\XbK,\zp(2)) \lra V_\zp \]
has a $G_K$-equivariant splitting (cf.\ Lemma \ref{aux3} below).
We thus obtain an exact sequence
\begin{equation}
0\lra H^2_\et(\UbK,\zp(2))^{G_K}\os{\partial_\et}{\lra}
\zp^{\op s} \lra H^1(K,H^2_\et(\XbK,\zp(2))).
\end{equation}
The condition {\bf (B)} together with an exact sequence
\[ 0\lra H^2_\et(\XbK,\zp(2))\os{p}{\lra}
H^2_\et(\XbK,\zp(2))\lra
H^2_\et(\XbK,\Z/p(2))\lra 0 \]
implies that $H^1(K,H^2_\et(\XbK,\zp(2)))$
is torsion-free. Hence so is the quotient \eqref{tor-free00}.
\end{pf}

Finally we give the following useful criterion for the conditions
{\bf (A)} and {\bf (B)}.
\addtocounter{thm}{5}
\begin{prop}\label{key4}
Consider the following conditions:
\begin{description}
\item[(A')]
the elliptic surface $\pi: \XbK \to \CbK$
has at least one additive {\rm(}singular{\rm)} fiber,
\item[(B')]
Let $Y=\XR_s$ be the special fiber of $\XR$ over the closed point $s \in \Spec(R)$, and let $C:\vg(Y,\Omega^2_Y)\to\vg(Y,\Omega^2_Y)$ be the Cartier operator.
Then $C-\id$ is injective.
\end{description}
Then 
{\bf (A')} $\Longrightarrow$ {\bf (A)}, and {\bf (B')} $\Longrightarrow$ {\bf (B)}.
\end{prop}
\begin{pf}
See Lem.\ref{aux4} below for the first assertion.
The second assertion follows from the torsion $p$-adic Hodge theory by Fontaine and Messing.
One has the isomorphism
$$
H^2_\et(\XbK,\Z/p(2))^{G_K}\cong\ker\big(
\vg(Y,\Omega^2_{Y/\F})\os{\varphi_2-1}{\lra}
H^2_\dR(Y/\F)
\big)
$$
where $1$ denotes the natural inclusion (cf.\ \cite{BM} Theorem 3.2.3).
Let $v:H^2_\dR(Y/\F)\to H^0(\Omega^2_Y/d\Omega^1_Y)$ be the natural map
and $C:H^0(\Omega^2_Y/d\Omega^1_Y)\os{\sim}{\to} H^0(\Omega^2_Y)$ be the Cartier operator.
Then $C\circ v \circ \varphi_2=\id$ and
$C\circ v \circ 1=C$. Hence the injectivity of $C-\id$ implies that of $\varphi_2-1$. 
\end{pf}
\subsection{Auxiliary results on Betti cohomology groups}
In the previous section, we used a number of results on \'etale cohomology groups of elliptic surfaces over $\ol K$,
which we prove here. We take the base-change by a fixed embedding $\ol K \hra \C$
and then replace the \'etale cohomology with the Betti cohomology $H^*=H^*_B$. For simplicity, we denote $X \otimes_K \C$, $S \otimes_K \C$ and $U \otimes_K \C$ by $X$, $S$ and $U$, respectively.
\begin{lem}\label{aux1}
Suppose that $\pi:X\to C$ is not iso-trivial.
Then $H^q(S,\Z) \simeq H^q(U,\Z)$ for $q \le 1$.
\end{lem}
\begin{pf}
The case $q=0$ is clear.
As for the case $q=1$, the Leray spectral sequence
\stepcounter{equation}
\begin{equation}\label{leray1}
E_2^{pq}=H^p(S,R^q\pi_*\Z_U)\Longrightarrow
H^{p+q}(U,\Z)
\end{equation}
yields an exact sequence
\[ 0\lra H^q(S,\Z)\lra H^q(U,\Z)\lra \vg(S,R^1\pi_*\Z). \]
Since $\vg(S,R^1\pi_*\Z)\hra H^1(X_t,\Z)^{\pi_1(S)}=0$,
the assertion follows.
\end{pf}
\stepcounter{thm}
\begin{lem}\label{aux2}
Assume that $H^3(X,\Z)$ is $p$-torsion-free. Then $H^2(X,\Z)$ is $p$-torsion-free as well.
\end{lem}
\begin{pf}
By the universal coefficient theorem
$H^1(X,M)\simeq\Hom(H_1(X,\Z),M)$ for an abelian group $M$.
Therefore if $H_1(X,\Z)=H^3(X,\Z)$ is $p$-torsion-free, then
the map $H^1(X,\Z)\to H^1(X,\Z/p)$ is surjective.
Equivalently, the multiplication by $p$ is injective on $H^2(X,\Z)$.
\end{pf}
\begin{lem}\label{aux3}
Let $\langle D\rangle\subset H^2(X,\Z)$ be the subgroup 
generated
by the cycle classes of irreducible components of $D=\sum_{i=1}^s D_i$.
Assume $p\geq5$.
Then the inclusion map
\[ \langle D\rangle\ot\zp\lra H^2(X,\zp) \]
has a natural splitting.
\end{lem}
\begin{pf}
Let us first assume that $X$ is minimal (hence $D_i$ is a N\'eron polygon).
Let $D_i^{(k)}$\,($k \ge 1$) be the irreducible components of $D_i$, and let \[ \bigoplus_{i,k} \ \Z \, [D_i^{(k)}] \lra H^2(X,\Z) \] be the map sending $[D_i^{(k)}]$ to its cycle class. Let \[ H^2(X,\Z) \lra \bigoplus_{i,k} \ H^2(D_i^{(k)},\Z) \simeq \bigoplus_{i,k} \ \Z[D_i^{(k)}] \;\; \hbox{ and } \;\; H^2(X,\Z) \os{e^*}\to H^2(C,\Z) \simeq \Z[C] \] be the natural pull-back maps. Then the composite map
\addtocounter{equation}{2}
\begin{equation}\label{intersectionmap}
\bigoplus_{i,k} \ \Z[D_i^{(k)}] \lra H^2(X,\Z) \lra\bigoplus_{i,k} \ \Z[D_i^{(k)}] \op \Z[C]
\end{equation}
is described as follows
\[ %\begin{equation}\label{splitting}
[D_i^{(k)}] \longmapsto \sum_{j,l} \ (D_i^{(k)},D_{j}^{(l)}) [D_{j}^{(l)}] +(D_i^{(k)},e(C))[C]
\] %\end{equation}
where $(-,-)$ denotes the intersection pairing. Therefore by an elementary computation on the intersection matrix,
we see that the composite map \eqref{intersectionmap} induces an injective map
\[
 \langle D\rangle\simeq
\bigoplus_{i,k\geq1} \ \Z[D_i^{(k)}]\bigg/\langle [D_i]-[D_j];i<j\rangle \;\;
\hra \;\; \bigoplus_{i\geq1,k\geq2} \ \Z[D_i^{(k)}] \op \Z[C], \]
whose matrix has determinant prime to 6.
Next we consider a general $X$.
Take the minimal model $\mu:X\to X_0$.
Then $H^2(X,\Z)$ is a direct sum of $H^2(X_0,\Z)$ and
the cycle classes of exceptional divisors.
Since $\langle \mu(D)\rangle\ot\zp$ is a direct summand of $H^2(X_0,\zp)$,
$\langle D\rangle\ot\zp$ is a direct summand of $H^2(X,\zp)$. 
This completes the proof.
\end{pf}
\stepcounter{thm}
\begin{lem}\label{aux4}
Assume that there is at least one additive singular fiber $E$.
Then $H_1(X,\Z)\simeq H_1(C,\Z)$. In particular
$H_1(X,\Z)=H^3(X,\Z)$ is torsion-free.
\end{lem}
\begin{pf}
It is enough to show that $e_*:H_1(C,\Z)\to H_1(X,\Z)$ is surjective.
Let $S^0\subset S$ be a Zariski open (non-empty) subset such that $U^0:=\pi^{-1}(S^0) \to S^0$ is a smooth fibration.
By the fibration exact sequence $1\to \pi_1(X_t)\to\pi_1(U^0)\to\pi_1(S^0)\to1$, which is split
by $e$, one has a surjection $\pi_1(X_t) \times \pi_1(S_0)^{\mathrm{ab}}
\to\pi_1(U_0)^{\mathrm{ab}}$.
Since $\pi_1(U_0)^{\mathrm{ab}}\to \pi_1(X)^{\mathrm{ab}}
=H_1(X,\Z)$ is surjective, it is enough to show that the map
$\pi_1(X_t)=H_1(X_t,\Z)\to H_1(X,\Z)$ is zero,
which follows from the fact that it factors through $H_1(E,\Z)=0$.
\end{pf}
\begin{lem}\label{aux5}
Let $\partial:H^2(U,\Z)\to H^3_D(X,\Z)=\Z^{\op s}$ be the boundary
map on the Betti cohomology arising from the localization exact sequence
\[ \cdots\lra H^2(U,\Z)\os{\partial}{\lra} \Z^{\op s}
\lra H^3(X,\Z)\lra H^3(U,\Z)\lra\cdots. \]
Then $\partial \ot\Q$ is surjective. In particular, if $H^3(X,\Z)$ is $p$-torsion-free,
then $\partial\ot\zp$ is surjective.
\end{lem}
\begin{pf}
We show that $H^3(X,\Q)\lra H^3(U,\Q)$ is injective.
The case $s=0$ (i.e., $X=U$) is obvious.
Assume $s>0$. One has 
$H^3(U,\Q)=H^1(S,R^2\pi_*\Q_U)=H^1(S,\Q)$ from the Leray spectral sequence \eqref{leray1}.
Since
\[ H^2(C,R^1\pi_*\Q_X) \lisom H^2_c(S,R^1\pi_*\Q_U) =H^0(S,R^1\pi_*\Q_U)^\lor=0, \]
one also has $H^3(X,\Q)=H^1(C,\Q)$ from the Leray spectral sequence for $\pi:X\to C$.
Thus the map $H^3(X,\Q) \to H^3(U,\Q)$ is identified with the natural restriction map
$H^1(C,\Q) \to H^1(S,\Q)$, which is clearly injective.
\end{pf}

%%%%%%%%%%%%%%%%%%%%%%%%%%%%%%%%%%%%%%%%%%%%%%
\newpage
\section{Application to Beilinson's Tate conjecture for $K_2$}\label{application-sect}
By the main result on elliptic surfaces (Theorem \ref{key3}) we obtain a computable upper bound of the Galois fixed part of \'etale cohomology group. We can now obtain a number of non-trivial examples of Beilinson's Tate conjecture for $K_2$ and also non-zero elements of the Selmer group of Bloch-Kato.

\subsection{Beilinson's Tate conjecture for $K_2$}\label{beilinson-sect}
\subsubsection{Example 1}\label{exmpsect1-1}
Let $k\geq1$ be an integer. Let $\zeta_k$ be a primitive $k$-th root of unity, and put $F=\Q(\zeta_k)$. Let $\O_F$ be the ring of integers in $F$. We consider an elliptic surface $\pi:\cX\to \P^1_\O$ over $\O:=\O_F[1/6k]$ whose generic fiber over $F(t)$ is given by
\begin{equation}\label{expm1}
3Y^2+X^3+(3X+4t^k)^2=0
\end{equation}
(see the beginning of \S \ref{knownsect} for the notion of elliptic surface over a ring). We assume that $\ol{X}:=\cX\otimes_\O\ol{F}$ is relatively minimal over $\P^1_{\ol{F}}$. The surface $\ol{X}$ is a rational surface if $1\leq k \leq 3$,
a K3 surface if $4\leq k \leq 6$ and the Kodaira dimension $\kappa(\ol{X})=1$ if $k\geq 7$.
Letting $t=0\in \P^1_\O(\O)$ and $t=\zeta_k\in \P^1_\O(\O)$
be the $\O$-rational points, we put $\cD_0:=\pi^{-1}(0)$
and $\cD_i:=\pi^{-1}(\zeta_k^i)$ ($1\leq i\leq k$)
which are relative normal crossing divisors over $\O$.
Then $\cD_i$ ($1\leq i\leq k$) is split multiplicative of type $I_1$ over $\O$, and $\cD_0$ is multiplicative of type $I_{3k}$ over $\O$. If $\sqrt{-3}\in F$ then $\cD_0$ is split. The singular fibers of $\ol{X}$ are $D_i:=\cD_i\otimes_\O \ol{F}$ and possibly $D_\infty:=\pi^{-1}(\infty)$. The type of the fiber $D_\infty$ is as follows:
\begin{center}
\begin{tabular}{c|ccc}
$k$ mod 3&1&2&3\\
\hline
$D_\infty$&IV*&IV&(smooth)
\end{tabular}
\end{center}
Put $\cU:=\cX-(\cD_1+\cdots+\cD_k)$, $U_F:=\cU\otimes_\O F$ and $\ol{U}:=U_F\otimes_F
\ol{F}$. We then discuss the Beilinson-Tate conjecture for $K_2(U_F)$, namely the
surjectivity of the higher Chern class map  
$$
K_2(U_F\otimes _FL)\ot\Q_p\lra H^2_\et(\ol{U},\Q_p(2))^{G_L}\qquad
(G_L:=\mathrm{Gal}(\ol{L}/L))
$$
for a finite extension $L$ of $F$.
\stepcounter{thm}
\begin{claim}\label{bc-claim}
Let $L\supset K\supset F$ and $[L:K]<\infty$. Then we have
\[ H^2_\et(\ol{U},\qp(2))^{G_K}=H^2_\et(\ol{U},\qp(2))^{G_L}.  \]
\end{claim}
\begin{pf}
The exact sequence \eqref{exact0} yields a commutative diagram ($V:=H^2(\ol X,\qp(2))$)
\[ \xymatrix{
 0 \ar[r] & H^2_\et(\ol{U},\qp(2))^{G_L} \ar[r]^-{\partial_\et} & \qp^k\ar[r] 
& H^1_\cont(G_L,V)\\
 0 \ar[r] & H^2_\et(\ol{U},\qp(2))^{G_K}\ar[r]^-{\partial_\et}\ar[u] & \qp^k\ar[r] 
 \ar@{=}[u]
& H^1_\cont(G_K,V)\ar[u]_{\mathrm{res}_{L/K}}.
} \]
The assertion follows from the fact that the right vertical arrow is injective.
\end{pf}

\medskip

Fix a prime $q$ of $F$ above $p$ which is prime to $6k$.
Let us denote by $D_q\subset G_F$ 
the decomposition group of $q$.
Let $K=F_q$ be the completion and $R\subset K$ the ring of integers.
Then $K$ is a finite unramified extension of $\Q_p$ and $D_q\cong G_K:=\mathrm{Gal}(\ol{K}/K)$.
We put $X_R:=\cX\otimes_\O R$, $U_R:=\cU\otimes_\O R$ and $D_R:=\sum_{i=1}^k\cD_i\otimes_\O R$.

Let us see the conditions {\bf (A)'} and {\bf (B)'}.
As we have seen above, if $(k,3)=1$ then $\ol{X}$ has the additive fiber $D_\infty$,
namely {\bf (A)'} is satisfied.
Next we see the condition {\bf (B)'}.
Let $\F=R/pR$ be the residue field.
Let $f(X)=-1/3(X^3+(3X+4t^k)^2)$ and 
$\sum_k a_kt^k$ the coefficient of $X^{p-1}$ in 
$f(X)^{(p-1)/2}$. Write $k=3\ell+a$ with $1\leq a\leq 3$.
$\vg(Y,\Omega^2_{Y/\F})$ has a basis
$\{t^idtdX/Y~;~0\leq i\leq \ell-1\}$.
One has
\begin{align*}
t^idt\frac{dX}{Y}&= t^if(X)^{(p-1)/2}dt\frac{dX}{Y^p}\\
&\equiv \sum_{k=1}^\ell a_{kp-i-1}t^{kp-1}dt\frac{X^{p-1}dX}{Y^p}
\quad(\mbox{ in }\vg(\Omega^2_{Y/\bF}/d\Omega^1_{Y/\bF}))\\
&=\sum_{k=1}^\ell a_{kp-i-1}C^{-1}(t^{k-1}dt\frac{dX}{Y})
\end{align*}
where $a_i:=0$ if $i<0$.
Thus the Cartier operator $C$ is described by a matrix
\[ A= \begin{pmatrix}
a_{p-1}&\cdots&a_{p-\ell}\\
\vdots&&\vdots\\
a_{\ell p-1}&\cdots&a_{\ell p-\ell}
\end{pmatrix}. \]
Hence {\bf (B)'} is satisfied if and only if there is no nontrivial solution of
\[
A\begin{pmatrix}
\alpha_1^p\\
\vdots\\
\alpha_\ell^p
\end{pmatrix}
=
\begin{pmatrix}
\alpha_1\\
\vdots\\
\alpha_\ell
\end{pmatrix}
\quad (\alpha_1,\cdots,\alpha_\ell\in \F).
\]

\medskip

We can now apply Theorem \ref{key3}.
For example, let $k=5$, $p=11$, $K=\Q_{11}$ and $R=\Z_{11}$.
Then the conditions {\bf (A)'} and {\bf (B)'} are satisfied.
A straightforward calculation (with the aid of computer) shows that
$\ol{\partial_\dR}\big(\Eis^{(99)}(X_R,D_R)_{\Z_{11}}\big)$ is generated by
\stepcounter{equation}
\begin{equation}\label{expm1-1}
- 4[D_1] - 4[D_4] + [D_5],\quad
 3[D_1] + 3[D_2] + [D_4],\quad
 [D_1] - 3[D_2] + [D_3]
\end{equation}
in $\op_{i=1}^5\F_{11}[D_i]$. We thus have
\begin{equation}\label{11-3}
\dim_{\Q_{11}} \ H^2_\et(\ol{U},\Q_{11}(2))^{G_K}
\leq 3.
\end{equation}
This remains true if we replace $K$ with any finite extension $L$ by Claim \ref{bc-claim}.
The same result can be checked for several $(k,p)$'s.

\medskip

One has
\begin{equation}\label{expm1-3}
\partial_\et \big(H^2_\et(U_{\ol{F}},\zp(2))^{G_F}\big)\subset \bigcap_{q|p} \ \partial_\et\big(
H^2_\et(U_{\ol{F}},\zp(2))^{D_q}\big)
\end{equation}
where $q$ runs over all primes of $F$ above $p$.
If each $\partial_\et
\big(H^2_\et(U_{\ol{F}},\zp(2))^{D_q}\big)$ is a direct summand
of $\zp^{\op k}$ then
the right hand side of \eqref{expm1-3}
is a direct summand as well.
Then one has
\[ \rank_\zp \, \bigcap_{q|p} \ \partial_\et \big(H^2_\et(\ol{U},\zp(2))^{D_q}\big)
\leq\dim_{\bF_p} \ \bigcap_{q|p} \ \ol{\partial_\et} \big(H^2_\et(\ol{U},\zp(2))^{D_q}\big) \]
and it is bounded by using the formal Eisenstein series.
For example, in case $k=5$, $p=11$, one has
\[ \bigcap_{q|p} \ \ol{\partial_\et} \big(H^2_\et(\ol{U},\zp(2))^{D_q}\big) \subset\bF_p([D_1]+[D_2]+\cdots+[D_5]), \]
hence
\begin{equation}\label{leq1f}
\rank_\zp \, H^2_\et(\ol{U},\zp(2))^{G_F}\leq
\rank_\zp \, \bigcap_{q|p} \ \partial_\et \big(H^2_\et(\ol{U},\zp(2))^{D_q}\big)\leq 1.
\end{equation}
The same thing is true if we replace $F$ with any finite extension $L$ by Claim \ref{bc-claim},
and same thing can be checked for several $(k,p)$'s.

We claim that the equalities hold in \eqref{leq1f}.
To do this, it is enough to show that 
$\partial K_2(U_F)\ot\Q\ne0$.
If $k=1$ then $X=X_1$ is the elliptic modular surface
for $\vg_1(3)$. Therefore it follows
from Beilinson's theorem that the image of
the boundary map on $K_2(X_1-D_1)$ is 1-dimensional.
Since there is the finite surjective map $X\to X_1$ induced from
$t\mapsto t^k$, one has $\dim_\Q \ \partial K_2(U_F)\ot\Q\geq1$.
(Explicitly, the symbol
\begin{equation}
\left\{
\frac{\sqrt{-3}(Y-X-4)-4(t^k-1)}{\sqrt{-3}(Y+X+4)-4(t^k-1)},
\frac{8(-\sqrt{-3}Y+3X+4t^k)(t^k-1)}{(X+4)^3}
\right\}
\end{equation}
in $\vg(\ol{U},\cK_2)\ot\Q$
has non-zero boundary.)

Finally let $F'\supset \Q(\sqrt{-3},\zeta_k)$ and
put $U'_{F'}=U_{F'}-D_0$.
Then one can also show
\[%\begin{multline}
\dim_{\Q} \, \partial K_2(U'_{F'})\ot\Q=
\dim_\qp \, H^2_\et(\ol{U'},\qp(2))^{G_{F'}}
= \dim_\qp \, \bigcap_{q|p} \ H^2_\et(\ol{U'},\qp(2))^{D_q}= 2
\]%\end{multline}
where $q$ runs over all primes of $F'$ above $11$.
(The symbol \eqref{symbol1} and
\begin{equation}
\left\{
\frac{\sqrt{-3}Y-3X-4t^k}{-8t^k},
\frac{\sqrt{-3}Y+3X+4t^k}{8t^k}
\right\}\in\vg(\ol{U'},\cK_2)\ot\Q
\end{equation}
span the 2-dimensional boundary image.)
This remains true if we replace $F'$
with any its finite extension or $U'_{F'}$ with
$U'_{F'}-$(any other fibers).
The same results can be checked for several $(k,p)$'s.

Summarizing above we obtain
\addtocounter{thm}{6}
\begin{thm}\label{expm1thm}
Let $F$ be a number field such that $F\supset\Q(\sqrt{-3},\zeta_k)$.
Let $\pi_F:X_F\to \P^1_F$ be the elliptic surface whose generic fiber
is defined
by \eqref{expm1}. Put $U_F:=X_F-(D_0+\cdots+D_k)-$(any other fibers).
Then $K_2(U_F)\ot\Q_p\to
H^2_\et(\ol{U},\Q_p(2))^{G_F}$ is surjective and
\stepcounter{equation}
\begin{equation}
\label{expm1-2}
\dim_\qp \, H^2_\et(\ol{U},\Q_p(2))^{G_F}=
\dim_\qp \, \bigcap_{q|p} \ H^2_\et(\ol{U},\qp(2))^{D_q}= 2
\end{equation}
if $(k,p)$ is one of the following cases:
\begin{center}
\begin{tabular}{c|c|c|c|c}
$k$&$5$&$7$&$11$&$13$\\
\hline
$p$&$7\leq p\leq97$&$\begin{array}{c} 5\leq p\leq109 \\ \hbox{and } \;p\neq7 \end{array}$
&$\begin{array}{c} 5\leq p\leq61 \\ \hbox{and } \; p\neq11 \end{array}$
&$\begin{array}{c} 5\leq p\leq 73 \\ \hbox{and } \; p\neq13\end{array}$
\end{tabular}
\end{center}
Explicitly $H^2_\et(\ol{U},\Q_p(2))^{G_F}$ is spanned by the image of the following symbols
\begin{equation}\label{symbol1}
\xi_1=\left\{
\frac{\sqrt{-3}(Y-X-4)-4(t^k-1)}{\sqrt{-3}(Y+X+4)-4(t^k-1)},
\frac{8(-\sqrt{-3}Y+3X+4t^k)(t^k-1)}{(X+4)^3}
\right\}
\end{equation}
\begin{equation}\label{symbol2}
\xi_2=\left\{
\frac{\sqrt{-3}Y-3X-4t^k}{-8t^k},
\frac{\sqrt{-3}Y+3X+4t^k}{8t^k}
\right\}
\end{equation}
{\rm(}Note $\partial(\xi_1)=[D_1]+[D_2]+\cdots+[D_k]$ and
$\partial(\xi_2)=k[D_0])$.{\rm)}
\end{thm}
\subsubsection{Example 2}\label{exmpsect1-2}
Let $k\geq1$ be an integer and $p$ a prime number such that $(p,6k)=1$.
Let $F$ be a number field such that $F\supset\Q(\sqrt{-1},\zeta_k)$.
Let $\pi_F:X_F\to \P^1_F$ be the elliptic surface whose generic fiber
is defined by
\begin{equation}\label{expm2}
3Y^2=2X^3-3X^2+t^k.
\end{equation}
The surface $X_{\ol{F}}$ is a rational surface if $1\leq k \leq 6$,
a K3 surface if $7\leq k \leq 12$ and $\kappa(X_{\ol{F}})=1$ if $k\geq 13$.
We put $D_0:=\pi^{-1}(0)$,
$D_i:=\pi^{-1}(\zeta_k^i)$ ($1\leq i\leq k$) and
$U_F:=X_F-(D_0+D_1+\cdots+D_k)$.
Then $D_i$ are split multiplicative fibers.
The other singular fiber of $X_{\ol{F}}$ is
$\pi^{-1}(\infty)$ whose type is as follows:
\begin{center}
\begin{tabular}{c|cccccc}
$k$ mod 6&1&2&3&4&5&6\\
\hline
$\pi^{-1}(\infty)$&II*&IV*&${\mathrm I}_0^*$&IV&II&(smooth)
\end{tabular}
\end{center}
In the same way as in \S \ref{exmpsect1-1} one can show the following theorem.
\addtocounter{thm}{4}
\begin{thm}\label{expm2thm}
Let $F$ be a number field such that $F\supset\Q(\sqrt{-1},\zeta_k)$.
Let $\pi_F:X_F\to \P^1_F$ and $D_i$ be as above.
Put $U_F:=X_F-(D_0+\cdots+D_k)-($any other fibers$)$.
Then $K_2(U_F)\ot\Q_p\to
H^2_\et(U_{\ol{F}},\Q_p(2))^{G_F}
$ is surjective and
\stepcounter{equation}
\begin{equation}
\label{expm2-22}
\dim_\qp \ H^2_\et(U_{\ol{F}},\Q_p(2))^{G_F}=\dim_\qp \ \bigcap_{q|p} \ H^2_\et(U_{\ol{F}},\Q_p(2))^{D_q}= 2
\end{equation}
if $(k,p)$ is one of the following cases:
\begin{center}
\begin{tabular}{c|c|c|c|c|c}
$k$&$7$&$8$&$9$&$11$&$13$\\
\hline
$p$&$\begin{array}{cc} 5\leq p\leq109 \\ \hbox{and }\; p\neq7 \end{array}$ &$5\leq p\leq 43$ &$5\leq p\leq 47$ &$\begin{array}{cc} 5\leq p\leq61 \\ \hbox{and }\; p\neq11 \end{array}$ &$\begin{array}{cc} 5\leq p\leq73 \\ \hbox{and }\; p\neq13  \end{array}$ 
\end{tabular}
\end{center}
Moreover $H^2_\et(\ol{U},\Q_p(2))^{G_F}$ is generated by the image of symbols
$$
\xi_1=\left\{
\frac{Y-\sqrt{-1}X}{Y+\sqrt{-1}X},-\frac{t^k}{2X^3}
\right\},\quad
\xi_2=\left\{
\frac{Y-(X-1)}{Y+(X-1)},-\frac{t^k-1}{2(X-1)^3}
\right\}.
$$
{\rm(}Note
$\partial(\xi_1)=k[D_0]$ and
$\partial(\xi_2)=[D_1]+[D_2]+\cdots+[D_k].${\rm)}
\end{thm}

\subsection{Bloch-Kato's Selmer group}\label{selmer-sect}
In the seminal paper \cite{BK2}, Bloch and Kato defined the $f$-part and the $g$-part in the continuous Galois cohomology
\[ H^1_f\subset H^1_g\subset H^1_\cont(G_k,V) \]
for a continuous $\Q_p$-representation $V$ of $G_k$, where $k$ is a $p$-adic local field or a number field. In particular, the $f$-part is often called the {\it Selmer group} of Bloch-Kato. Let $X$ be a projective smooth variety over $k$. Then the higher Chern character map on Quillen's $K$-theory induces a $p$-adic {\it regulator map}
\begin{equation}\label{reg-g}
\varrho_{1,2}:K_1(X)^{(2)}\ot\qp\lra H^1_g(G_k,H^2(\ol{X},\Q_p(2)))
\end{equation}
(\cite{sato-saito}, (3.3.1), Lemma 3.5.2). Let $\cX$ be a regular proper flat scheme over the ring of integers $\O_k$ which has a finite surjective map $h:\cX\otimes_{\O_k} k\to X$. We define the {\it integral part} of $K$-theory as
\[ K_1(X)^{(2)}_\Z:=\Image\big(K_1(\cX)^{(2)}\to K_1(\cX\otimes_{\O_k}k)^{(2)}\os{h_*}{\to} K_1(X)^{(2)}\big), \]
which is independent of the choice of $\cX$ (\cite{Scholl}). Then the image of the integral $K$-theory under $\rho_{1,2}$ is contained in the $f$-part (\cite{Ni}, \cite{Sa}), and we obtain the following integral regulator map:
\begin{equation}\label{reg-f}
\xymatrix{
\varrho_{1,2}:K_1(X)^{(2)}_\Z\ot\qp\ar[r]& H^1_f(G_k,H^{2}(\ol{X},\Q_p(2))).}
\end{equation}
When $k$ is a number field, Bloch and Kato conjecture that \eqref{reg-g} and \eqref{reg-f} are bijective, but it is a widely open problem.
\par\medskip
%We discuss the regulator map in case $(i,j)=(1,2)$.
To consider the surjectivity of these regulator maps, we decompose each of them into decomposable part and indecomposable part. The latter is much more difficult to compute in general. For a finite field extension $L/k$, let $\pi_L$ be the composite of the product and the norm map of $K$-groups:
\[  L^\times \otimes \Pic(X\otimes_k L) \lra K_1(X\otimes_k L)^{(2)} \lra K_1(X)^{(2)}. \]
We define $K_1^{\dec} (X)^{(2)}$, the {\it decomposable} part of $K_1(X)^{(2)}$, as the $\bQ$-subspace of $K_1(X)^{(2)}$ generated by the image of
\[\begin{CD} \displaystyle \bigoplus_{L/k} \ L^\times \otimes \Pic(X\otimes_k L) 
@>{\sum \, \pi_L}>> K_1 (X)^{(2)}, \end{CD}\]
where $L$ runs through all finite extensions of $k$. 
We put $K_1^\ind(X)^{(2)}:=K_1(X)^{(2)}/
K_1^{\dec}(X)^{(2)}$ and call it the {\it indecomposable $K_1$} of $X$.
We put
\[ K_1^{\dec}(X)^{(2)}_\Z:=K_1^{\dec}(X)^{(2)}\cap K_1(X)^{(2)}_\Z,
\quad K_1^\ind(X)^{(2)}_\Z:=\Image (K_1(X)^{(2)}_\Z \to K_1^\ind(X)^{(2)}). \]
One easily sees $K_1(X)^{(2)}_\Z\cong K_1^{\dec}(X)^{(2)}_\Z\op K_1^{\ind}(X)^{(2)}_\Z$.
Put \[ V^\ind:=H^2_\et(\ol{X},\qp(2))/\NS(\ol{X})\ot\qp(1). \]
Then the regulator map $\varrho=\varrho_{1,2}$ induces a commutative diagram 
$$
\xymatrix{
K^\ind_1(X)^{(2)}\ot\qp\ar[r]^-{\varrho}& H^1_g(G_k,V^\ind)\\
K^\ind_1(X)^{(2)}_\Z\ot\qp\ar[u]\ar[r]^-{\varrho}& H^1_f(G_k,V^\ind)
\ar[u].}
$$
This is the difficult part of the regulator map $\varrho$.
It is very difficult even to ask whether it vanishes or not in general.
We apply the results in the previous section to show the non-vanishing in case $X$
is an elliptic surface.
The strategy is as follows.

Let $\pi:X\to C$ be an elliptic surface over a number field or a $p$-adic local field
$k$ and $D_i=\pi^{-1}(P_i)$ split multiplicative fibers
over $P_i\in C(k)$. 
Let $\widetilde{D}_i\to D_i$ be the normalization.
Then the cokernel of $K_1(\widetilde{D}_i)^{(2)}\to K'_1(D_i)^{(1)}$ is a $\Q$-vector space of rank one and the generator is described in the following way.
Let $D_i^\dag\subset D_i$ be the N\'eron polygon. 
Let $D_i^\dag=\sum E_j$ be the irreducible decomposition and $Q_j$ the intersection points in $E_j$ and $E_{j+1}$.
Let $f_j$ be a rational function on $E_j$ such that ${\mathrm{div}}_{E_j}(f_j)=Q_j-Q_{j-1}$.
There is a localization exact sequence
\[ 0 \lra K'_1(D_i^\dag)^{(1)} \lra K_1(D_i^\dag- {\textstyle\coprod_j \, Q_j})^{(1)} \os{d}{\lra} \bigoplus_j \ K_0(Q_j) \]
and the element
\[ (f_j)_j\in \bigoplus_j \ K_1(E_j-\{Q_j,Q_{j+1}\})^{(1)} \cong K_1(D_i^\dag-{\textstyle\coprod_j \, Q_j})^{(1)} \]
lies in the kernel of $d$. Hence it defines an element 
$\xi'_i\in K'_1(D_i^\dag)^{(1)}$.
Let $\xi_i\in K'_1(D_i)^{(1)}$ be the image of $\xi'_i$
via the push-forward
$K'_1(D_i^\dag)^{(1)}\to K'_1(D_i)^{(1)}$.
Then $\xi_i$ gives a generator
of the cokernel of $K'_1(\widetilde{D}_i)^{(1)}
\to K'_1(D_i)^{(1)}$.

Put $D=\sum_{i=1}^sD_i$.
Let $K_1(X)_D$ be the subgroup of $K_1(X)^{(2)}$ generated by
the images of $K_1(\widetilde{D}_i)^{(1)}$.
Let $\langle D\rangle\subset\NS(X_{\ol{F}})$ be the subgroup
generated by irreducible components of $D$.
Put $V:=H^2_\et(\ol{X},\Q_p(2))/\langle D\rangle\ot\Q_p(1)$.
One has a commutative diagram
\begin{equation}\label{diagram2}
\begin{CD}
@.K_2(U)\ot\Q_p
@>{\partial}>>\bigoplus_{i=1}^s \ \Q_p\cdot\xi_i@>>>K_1(X)^{(2)}
/K_1(X)_{D}\ot\Q_p\\
@.@VVV@VV{\cong}V@VV{\varrho}V\\
0@>>>H^2_\et(\ol{U},\Q_p(2))^{G_k}@>{\partial_\et}>>
\bigoplus_{i=1}^s \ \Q_p[D_i]@>>>
H^1_g(G_k,V)
\end{CD}
\end{equation}
with exact rows.
We also write the image of $\xi_i$ in 
$K_1(X)^{(2)}/K_1(X)_D$ by $\xi_i$.
We have
\begin{equation}\label{diagram4}
\sum a_i\varrho (\xi_i)\neq0 \; \mbox{ in } \; H^1_g(G_k,V)
\Longleftrightarrow
\sum a_i[D_i] \not\in \partial_\et H^2_\et(\ol{U},\Q_p(2))^{G_k}.
\end{equation}
Let $Z\subset X$ be an irreducible curve and $\tilde{Z}\to Z$ the normalization.
Let \[ v_Z:H^1(G_k,H^2(\ol X,\qp(2)))\lra
H^1(G_k,H^2(\tZ,\qp(2)))
=\Q\ot\plim{n} \ k^\times/p^n \]
be the pull-back map.
Assume that
\begin{equation}\label{orthogonal}
v_Z(\sum a_i \varrho(\xi_i))=0 \quad
\text{ for each generator $Z$ of }\NS(\ol{X})_\Q.
\end{equation}
Since the intersection form on $\NS(\ol X)_\Q$ is non-degenerate, one has
 \[ \NS^{\bot} \op \NS(\ol X)_\qp \isom H^2_\et(\ol X,\qp(1))\]
 and $\NS^{\bot} \otimes \qp(1) \simeq V^\ind$, where $\NS^\bot$ denotes the orthogonal complement of $\NS(\ol X)_\qp$
in $H^2_\et(\ol X,\qp(1))$.
Therefore \eqref{orthogonal} implies that
$\varrho(\sum a_i \xi_i)$ is contained in
\[ \hbox{ the image of } \; H^1(\qp,\NS^\bot\otimes \qp(1)) \to H^1(\qp,V) \]
and hence
\begin{equation}\label{diagram40}
\sum a_i\varrho (\xi_i)\neq0 \; \mbox{ in } \; H^1_g(G_k,V^{\mathrm{ind}})
\end{equation}
under \eqref{diagram4}.
If $\sum a_i \xi_i$ is integral, then we obtain a non-zero element of the Selmer group $H^1_f(G_k,V^{\mathrm{ind}})$ of Bloch-Kato.

\subsubsection{Example 1}
Recall the elliptic surface \eqref{expm1} in \S\ref{exmpsect1-1}.
If $k$ is odd, 
the N\'eron-Severi group $\NS(\ol{X})\ot\Q$
is generated by irreducible components of singular fibers and the section
$e(\P^1)$ (\cite{stiller2}  Example 3).
Let $\xi_i\in K_1(X)^{(2)}$ be the element which is defined from
$$
f_i=\frac{Y-(X+4)}{Y+(X+4)}\big|_{D_i}\in K_1(D_i^{\mathrm{reg}})^{(1)}
$$
in the above way.
It is easy to check that it satisfies \eqref{orthogonal}.
One can also check that if $(k,6)=1$ then each $\xi_i$ is integral
(i.e. $\xi_i\in K_1(X)^{(2)}_\Z$).
By Theorem \ref{expm1thm}, we have
\addtocounter{thm}{6}
\begin{thm}\label{expm1thm1}
Let $\pi_F:X_F\to \P^1_F$ and $D_i$ be as in 
Theorem \ref{expm1thm}. Then the composition
\[ \left(\bigoplus_{i=1}^k \ \Q_p\xi_i\right)\bigg/ \Q_p(\xi_1+\cdots+\xi_k) \os{\varrho}{\lra} H^1_f(G_F,V^\ind) \lra \bigoplus_{q|p} \ H^1_f(G_{F_q},V^\ind) \]
is injective
if $(k,p)$ is one of the following cases:
\begin{center}
\begin{tabular}{c|c|c|c|c}
$k$&$5$&$7$&$11$&$13$\\
\hline
$p$&$7\leq p\leq97$&$\begin{array}{c}5\leq p\leq109 \\ \hbox{and }\; p\neq7 \end{array}$
&$\begin{array}{c}5\leq p\leq61\\ \hbox{and }\;p\neq11\end{array}$
&$\begin{array}{c}5\leq p\leq 73\\ \hbox{and }\;p\neq13\end{array}$
\end{tabular}
\end{center}
\end{thm}

\subsubsection{Example 2}
Recall the elliptic surface \eqref{expm2} in \S \ref{exmpsect1-2}.
If $(k,30)=1$ the N\'eron-Severi group $\NS(\ol{X})\ot\Q$
is generated by the irreducible components of singular fibers and the section
$e(\P^1)$ (\cite{stiller2} Example 4).
Similarly to before, Theorem \ref{expm2thm} implies
\begin{thm}\label{expm2thm2}
Let $\pi_F:X_F\to \P^1_F$ and $D_i$ be as in 
Theorem \ref{expm2thm}. 
Then the composition 
\[ \left(\bigoplus_{i=1}^k \ \Q_p\xi_i\right)\bigg/ \Q_p(\xi_1+\cdots+\xi_k)
\os{\varrho}{\lra}
H^1_f(G_F,V^\ind)
\lra
\bigoplus_{q|p} \ H^1_f(G_{F_q},V^\ind) \]
is injective
if $(k,p)$ is one of the following cases:
\begin{center}
\begin{tabular}{c|c|c|c}
$k$&$7$&$11$&$13$\\
\hline
$p$&$\begin{array}{c}5\leq p\leq109\\ \hbox{and }\;p\neq7\end{array}$
&$\begin{array}{c}5\leq p\leq61\\ \hbox{and }\;p\neq11\end{array}$
&$\begin{array}{c}5\leq p\leq73\\ \hbox{and }\;p\neq13\end{array}$
\end{tabular}
\end{center}
\end{thm}

%%%%%%%%%%%%%%%%%%%%%%%%%%%%%%%%%%%%%%%%%%%%%%
\newpage

\section{An elliptic $\boldsymbol{K}$3 surface over $\boldsymbol{\qp}$ with finitely many torsion zero-cycles}
\label{mordellsect}
\setcounter{equation}{0}
In this section, we discuss another application to $p$-adic regulator map on $K_1$
and finiteness of torsion zero-cycles.  
Hereafter all cohomology groups of schemes are taken over the \'etale topology.

\medskip

To state the main result of this section, we start from a numerical condition
on a prime number $p\geq5$.
Let $q$ be an indeterminate, and define formal power series $E_1,E_{3,a},E_{3,b}$ in $\bZ\psrqq$ as
\begin{align*} E_1&:=1+6\sum_{k=1}^\infty \ \left(\frac{q^{3k-2}}{1-q^{3k-2}} -\frac{q^{3k-1}}{1-q^{3k-1}} \right), \\
 E_{3,a}&:=1-9\sum_{k=1}^\infty \ \left(\frac{(3k-2)^2q^{3k-2}}{1-q^{3k-2}}-\frac{(3k-1)^2q^{3k-1}}{1-q^{3k-1}}\right), \\
 E_{3,b}&:=\sum_{k=1}^\infty \ \frac{k^2(q^k-q^{2k})}{1-q^{3k}},
\end{align*}
which are obtained from the $q$-expansion of Eisenstein series for $\Gamma_1(3)$ at the cusp $\tau=\infty$. Let $t\in \bZ[2^{-1}]\psfqq$ satisfy $t^4=E_{3,a}/(E_1)^3$:
\[ t=1-\frac{27}{4}q+\frac{1053}{32}q^2-\frac{23085}{2^7}q^3+\frac{2130003}{2^{11}}q^4-\frac{ 49565277}{2^{13}}q^5+\cdots. \]
Put
\begin{equation}\label{fig}
g(q)=-\frac{27t}{4}E_{3,b}, \quad
f_1(q)=-\frac{27t}{4(t-1)}E_{3,b}, \quad
f_2(q)=-\frac{27t}{4(t+1)}E_{3,b}
\end{equation}
and express
\[ f_1(q)=1+\sum_{i=1}^\infty\frac{a_iq^i}{1-q^i}, \quad
f_2(q)=\sum_{i=1}^\infty\frac{b_iq^i}{1-q^i}, \quad
g(q)=\sum_{i=1}^\infty\frac{c_iq^i}{1-q^i} \]
with $a_i,b_i,c_i\in\zp$. Then we say that a prime number $p\geq5$ satisfies {\bf C(p)} if the following two conditions are satisfied:
\begin{description}
\item[C(p)-1.]
Let $k_p$ be the coefficient of $x^{p-1}t^{p-1}$ in the polynomial $(-3(x^3+(3x+4t^4)^2))^{(p-1)/2}$ (cf. \S \ref{sect6-1} {\bf Fact 1}). Then $k_p\not\equiv 1
\text{ mod }p.$
\item[C(p)-2.]
There is no solution $n\in \bZ_p$ which satisfy all of the following congruent relations:
\begin{equation}\label{eq6-2-5-cp}
\begin{cases} a_{ip}-b_{ip}+n c_{ip}\equiv 0\mod p^2\bZ_p & (i=1,2,\dotsc,p-1)\\ %\mbox{ and }
a_{p^2}-b_{p^2}+n c_{p^2}\equiv 0\mod p^4\bZ_p. \end{cases} \end{equation}
\end{description}

The following is the main result of this section.
\setcounter{thm}{2}
\begin{thm}\label{finite1}
Let $\pi:X_\Q\to \P^1_\Q$ be the elliptic K{\rm 3} surface
over $\Q$ whose general fiber $\pi^{-1}(t)$ is an elliptic curve given by
\stepcounter{equation}
\begin{equation}\notag %\label{expm1}
3Y^2+X^3+(3X+4t^4)^2=0.
\end{equation}
Let $p\ge 5$ be a prime number which satisfies {\bf C(p)}, and put $\XQp:=X_\Q\otimes_\Q\qp$ and $\ol X:= X \otimes_{\qp} \ol{\Q}_p$.
Then the p-adic regulator
\begin{equation}\label{finite1rho}
\varrho:K_1(\XQp)^{(2)}\ot\qp \lra H^1_g(\qp,H^2(\ol X,\qp(2)))
\end{equation}
is surjective.
\end{thm}
The condition {\bf C(p)} is in fact computable, and the authors checked that it holds for $p=7,11,19,23,31$ with the aid of computer.
They expect it is true for any $p$ with $p\equiv 3$ mod 4, but do not have any idea to prove it.
When $p\equiv 1$ mod 4, they do not know any example of $p$ which satisfies {\bf C(p)}, nor whether the conclusion of Theorem {\rm \ref{finite1}} (or Corollary {\rm\ref{finite2}} below) holds.

Theorem \ref{finite1} immediately implies the following finiteness result by \cite{sato-saito} Theorem 3.1.1.
\addtocounter{thm}{1}
\begin{cor}\label{finite2}
The $p$-primary torsion subgroup of $\CH_0(\XQp)$ is finite, if $p$ satisfies {\bf C(p)}
{\rm(}e.g. $p=7,11,19,23,31${\rm)}.
\end{cor}
See \S\ref{sect6-3} below for the finiteness of the full torsion part $\CH_0(\XQp)_\tors$. The proof of Theorem \ref{finite1} is divided into two parts, i.e., the surjectivity onto $H^1_g/H^1_f$ and $H^1_f$ (see \S \ref{exmpsect0}). The former one will be reduced to results of Flach and Mildenhall (see \S\ref{exmpsect1} below).
To show the latter one, we shall construct a new indecomposable element in $K_1(\XZp)^{(2)}$ and show the non-vanishing in $H^1_f$ under the condition {\bf C(p)}, where Theorem \ref{key3} will play an essential role (see Proposition \ref{claimstep2} below).
\begin{rem}
The elliptic surface defined by 
$3Y^2+X^3+(3X+4t)^2=0$ {\rm(}$t\in \P^1${\rm)}
is the universal family of elliptic curves over the modular curve
$X_1(3)$. The elliptic K{\rm 3} in Theorem {\rm \ref{finite1}} is a finite covering of this surface, but
it is no longer modular.
\end{rem}

\subsection{Preliminary facts}\label{sect6-1}
Before proving the theorem, we write down some facts on the $X_\Q$ and $X_{\ol \bQ}:=X_\Q \otimes_{\bQ} \ol {\bQ}$
which we shall need later (proof is left to the reader).

\begin{description}
\item[Fact 1.]
$X_\Q$ has good reduction at $p\geq5$.
For $p \ge 5$, the reduction $Y_p$ of $X_\Q$ at $p$ is ordinary if and only if $p\equiv1$ mod 4,
and super-singular if and only if $p\equiv 3$ mod 4.
\end{description}
Note that the Cartier operator $C:H^0(\Omega^2_{Y_p})\to H^0(\Omega^2_{Y_p})$ is given by
multiplication by $k_p$ in {\bf C(p)-1}. $Y_p$ is ordinary if and only if
$k_p\equiv 0$ mod $p$ and otherwise it is super-singular.

\begin{description}
\item[Fact 2.]
The functional $j$-invariant is $27(9-8t^4)^3/((1-t^4)t^{12})$.
There are 5 multiplicative fibers over $t=0,~\pm1,~\pm\sqrt{-1}$,
and one additive fiber over $t=\infty$:
\end{description}
\begin{center}
\begin{tabular}{c|c|c|c}
$t$ &0&$\pm1,\pm\sqrt{-1}$&$\infty$\\
\hline
$\pi^{-1}(t)$&$I_{12}$&$I_1$&$IV^*$
\end{tabular}
\end{center}

\begin{description}
\item[Fact 3.]
The N\'eron-Severi group $\NS(X_{\ol{\Q}})$ has rank 20.
Explicitely it is generated by
the following irreducible curves:
\begin{enumerate}
\item[{\rm (i)}]
irreducible components of $\pi^{-1}(0)$ ($12$-components),
\item[{\rm (ii)}]
irreducible components of $\pi^{-1}(\infty)$ ($7$-components),
\item[{\rm (iii)}]
the section at infinity $e(\P^1)=E$,
\item[{\rm (iv)}]
the section $C$ defined by $x=-4t^4/3$ and $y=8t^6/9$.
\end{enumerate}
Note that $\pi^{-1}(0)\equiv\pi^{-1}(\infty)$
is the only relation in $\NS(X_{\ol{\Q}})$ among the above. 
\end{description}
\begin{description}
\item[Fact 4.]
Put $\ol {Y_p} :=Y_p \otimes_{\F_p} {\ol{\F}_p}$.
The N\'eron-Severi group $\NS(\ol {Y_p})$ has rank $20$ if
$p\equiv1$ mod $4$, and rank $22$ if $p\equiv3$ mod $4$.
\end{description}
When $p\equiv1$ mod $4$, one has $\NS(X_{\ol \Q})\ot\Q
=\NS(\ol {Y_p})\ot\Q$. When $p\equiv3$ mod $4$, one has $\NS(X_{\ol \Q})\ot\Q
\ne\NS(\ol {Y_p})\ot\Q$ and
we will describe curves on $\ol {Y_p}$ which do not come
from $X_{\ol \Q}$ in \S \ref{exmpsect1}.

\subsection{Proof of Theorem \ref{finite1}}\label{exmpsect0}
Let $\XZp\to\Spec(\zp)$ be a projective smooth model of $\XQp=X_\Q \otimes_\Q\Q_p$, and let $Y \to\Spec(\F_p)$ be its special fiber. 
Since $K_1^\dec(\XQp)^{(2)}\ot\qp$ is onto $H^1(\qp,\NS(\ol X)\ot\qp(1))$ (cf.\ \cite{S} Lemma 3.6), we may replace the target with $H^1_g(\qp,V)$, where we put
\[ V^\ind:=H^2(\ol X,\qp(2))/\NS(\ol X)\ot\qp(1). \]
We shall prove it in the following steps.

\begin{description}
\item[Step 1.]
$K_1(X)^{(2)}\ot\qp\to H^1_g(\qp,V^\ind)/H^1_f(\qp,V^\ind)$ is surjective for any $p\geq5$.
\item[Step 2.]
$K_1(\cX)^{(2)}\ot\qp\to H^1_f(\qp,V^\ind)$ is surjective for $p\geq 5$ satisfying {\bf C(p)}.
\end{description}
\subsubsection{Proof of Step 1}\label{exmpsect1}
By the Tate conjecture for $Y$ (\cite{ArSw}) and the same arguments as in \cite{LS} Theorem 5.1,
 there is a canonical isomorphism
\[ H^1_g(\qp,V^\ind)/H^1_f(\qp,V^\ind)\simeq\NS(Y)/\NS(\XQp)\ot\qp \]
and the composition
\[  K_1(\XQp)^{(2)}\ot\qp \to H^1_g(\qp,V^\ind)/
H^1_f(\qp,V^\ind)\simeq\NS(Y)/\NS(\XQp)\ot\qp \]
is given by the boundary map $K_1(\XQp)^{(2)}\to
\NS(Y)\ot\Q$ arising from the localization exact sequence in $K$-theory.
\begin{prop}
Let $C_i$ $(i=1,2)$ be elliptic curves over $\Q$ defined by equations
\[ 3y^2+x^4+1=0\quad\hbox{and}\quad u^2+4v^4+3=0, \]
respectively. Then there is a dominant rational map
\[\xymatrix{ f:C_1\times C_2 \ar@{.>}[r] & X_0 }\]
of degree $8$ given by $(x,y)\times(u,v)\mapsto (X,Y,t)=((ux)^4,u^6x^4y,uvx)$.
\end{prop}
\begin{pf}
Straight-forward.
\end{pf}
Note that $C_1$ and $C_2$ are isomorphic to each other up to twist.
Let $S \to C_1\times C_2$ be a birational transformation
such that $f$ is extended to a morphism $\bar{f}:S\to X$.
Put \[ V^\ind(-):= H^2(-,\qp(2))/\NS(-)\ot\qp(1). \]
Then one has
\[ \xymatrix{ V^\ind(\ol X) \; 
\ar@{^{(}->}[r] & V^\ind(\ol S) \lisom V^\ind(\ol{C_1\times C_2}). } \]
We show that the first map is bijective.
Indeed, we have $\dim_\qp \, V^\ind(\ol X)=2$ by {\bf Fact 3}.
On the other hand, since $C_1$ and $C_2$ are isomorphic to a CM elliptic curve up to twist, 
$V^\ind(\ol{C_1\times C_2})$ is also 2-dimensional. 
Hence $V^\ind(\ol X)=V^\ind(\ol S)$.

In order to show {\bf Step 1} it is enough to
show that
\[ K_1(C_1\times C_2)^{(2)}\ot\qp \to H^1_g(\qp,
V^\ind(\ol{C_1\times C_2}))/ H^1_f(\qp,V^\ind(\ol {C_1\times C_2})) \]
is surjective. We may replace $\qp$ with arbitrary finite extension $K/\qp$ by a standard norm argument. Fix $K$ such that $C_{1,K} \simeq C_{2,K}(=:C)$ with smooth reduction $\cC_s$ and such that
\[ \End(C)=\End(\ol C) \quad \hbox{ and } \quad \End(\cC_s)=\End(\ol{\cC_s}). \]
We show that
\begin{align*}
K_1(C\times C)^{(2)}\ot\qp &\lra H^1_g(K,V^\ind(\ol{C\times C}))/ 
H^1_f(K,V^\ind(\ol{C\times C}))\\
&\simeq \End(\ol{\cC_s})\ot\qp/\End(\ol C)\ot\qp
\end{align*}
is surjective. If $\cC_s$ is ordinary (i.e., $p\equiv 1$ mod 4), the target is zero. If $\cC_s$ is super-singular (i.e., $p\equiv 3$ mod 4), it is generated by the Frobenius endomorphism and its composition with the CM endomorphism. The surjectivity then follows from \cite{flach} Proposition 2.1 or \cite{mildenhall} Theorem 5.8. This completes the proof of {\bf Step 1}.
\subsubsection{Proof of Step 2}\label{exmpsect2}
This is a crucial step. We first show that
$H^1_f(\qp,V^\ind)$ is $1$-dimensional over $\qp$.
\begin{lem}\label{lem6-2-2}
There is an isomorphism
\[ H^1_f(\qp,V^\ind)\simeq H^2(X,\O_X). \]
%\simeq H^3_\syn(\XZp,{\mathscr S}_\qp(2))/K_1^{\dec}(\XZp)^{(2)}\ot\qp where $K_1^\dec(\XZp)^{(2)}:=K_1^\dec(X)^{(2)} \cap K_1(\XZp)^{(2)}$.
\end{lem}
\begin{pf}
We first note the following isomorphisms:
\begin{align*}
 H^1_f(\qp,H^2(\ol X, \qp(2))) & \us{\;c^3}{\isom} H^3_\syn(\XZp,\cS_\qp(2)) \\ & \us{\;\delta^3}{\lisom} \Coker\big(1-f_2 : \vG(X,\Omega_{X/\qp}^2) \to H^2_{\dR}(X/\qp)\big) \\ &\;\; =: H^2_{\dR}(X/\qp)/(1-f_2),
\end{align*}
where $f_2$ denotes the map induced by $f_2$ in Definition \ref{def1-1-1}. The first isomorphism, induced by the map $c^3$ in Theorem \ref{thm1-1-1}, is due to Langer-Saito \cite{LS} Theorem 6.1, a special case of the $p$-adic point conjecture raised by Schneider \cite{Sch2} (cf.\ \cite{Ne} III Theorem (3.2)). The second isomorphism $\delta^3$ is the connecting map induced by the definition of $\cS_n(2)$ (see \cite{Sch2} p.\ 242 or \cite{Ne} III (3.1.1) for details). We denote the composite of these isomorphisms by $\alpha$.
\par
Let $L$ be a finite extension of $\qp$ for which $\NS(X_L)$ has rank $20$, i.e., $\NS(X_L)\cong \NS(\ol X)$ (cf.\ {\bf Fact 3}). We construct a commutative diagram as follows:
\addtocounter{equation}{2}
\begin{equation}\label{eq6-2-1}\xymatrix{\NS(X_L) \otimes L \ar[d]_{\tr_{L/\qp} \circ \varrho_{\dR,L}} \ar@{->>}[r]^-{\varpi_{L/\qp}}  & H^1_f(\qp,\NS(\ol X)\otimes \qp(1)) \ar[d] \\
H^2_{\dR}(X/\qp)/(1-f_2) \ar[r]_-\alpha \ar@<-1pt>@{}[r]^-\sim &  H^1_f(\qp,H^2(\ol X, \qp(2))), }\end{equation}
where the map $\varrho_{\dR,L}$ denotes the de Rham Chern class with values in $H^2_{\dR}(X_L/L)$ and $\tr_{L/\qp}$ denotes the trace map $H^2_{\dR}(X_L/L) \to H^2_{\dR}(X/\qp)$. The arrow $\varpi_{L/\qp}$ is defined as the composite (`$\exp$' denotes the exponential map defined in \cite{BK2} Definition 3.10)
\begin{align*} \NS(X_L) \otimes L & \us{\id \,\otimes\, \exp}\isom \NS(X_L) \otimes H^1_f(L, \qp(1)) \\ & \; \us{\cup}\isom H^1_f(L,\NS(\ol X)\otimes \qp(1)) \us{\Cor_{L/\qp}}\lra H^1_f(\qp,\NS(\ol X)\otimes \qp(1)), \end{align*}
which is surjective by a standard norm argument. In view of the right exactness of $H^1_f(\qp,-)$ (cf.\ \cite{Ne} III (1.7.3)), once we show the commutativity of the diagram \eqref{eq6-2-1}, the assertion of the lemma will follow from the fact that the Hodge Chern class map $\NS(X_L) \otimes L  \to H^1(X_L,\Omega^1_{X_L/L})$ is surjective. Therefore it remains to show the commutativity of \eqref{eq6-2-1}.
\par
Let $\fo$ be the integer ring of $L$, and put $\cX\hspace{-2pt}_{\fo}:=\cX \otimes_\zp \fo$. There is a commutative diagram of pairings
\[\xymatrix{\quad H^2_{\dR}(X_L/L)  \;\; \times \;\;\;\qquad L \;\qquad \ar@<18mm>[d]^{\delta^1}_{\wr\hspace{-1pt}} \ar[r]^-{\cup} & H^2_{\dR}(X_L/L) \ar[d]^{\delta^3} \\ H^2_{\syn}(\cX\hspace{-2pt}_\fo, \cS_\qp(1)) \times H^1_{\syn}(\fo, \cS_\qp(1)) \ar@<13mm>[u]^\iota \ar[r]^-{\cup} \ar@<-13mm>[d]_{c^2} \ar@<18mm>[d]^{c^1} &  H^3_\syn(\cX\hspace{-2pt}_\fo,\cS_\qp(2)) \ar[d]^{c^3} \\
H^2_\et(\ol X,\qp(1))^{G_L} \times \; H^1(L,\qp(1)) \;  \ar[r]^-{\cup} & H^1(L,H^2(\ol X, \qp(2))), }\]
where $\iota$ (resp.\ $\delta^1$) denotes the natural map (resp.\ connecting map) induced by the definition of $\cS_n(1)$ (cf.\ Definition \ref{def1-1-1}, \cite{Ne} III (3.1.1)). This diagram commutes by the definition of the product structure of syntomic cohomology and the compatibility of $c^i$'s with product structures (cf.\ Theorem \ref{thm1-1-1}). Moreover we have
\[ c^1 \circ \delta^1 = \exp \quad \hbox{ and } \quad c^3 \circ \delta^3=\alpha \]
by the definitions of these maps, and we have the following commutative diagram of the $1$st Chern class maps by \cite{Ts1} Proposition 3.2.4\,(3), Lemma 4.8.9:
\[\xymatrix{ & \Pic(\cX\hspace{-2pt}_\fo) \ar[ld]_{\varrho_{\dR}} \ar[d]_{\varrho_\syn} \ar[rd]^{\varrho_\et} \\
 H^2_{\dR}(X_L/L) & \ar[l]_-{\;\iota} H^2_{\syn}(\cX\hspace{-2pt}_\fo, \cS_\qp(1)) \ar[r]^-{c^2} & H^2_\et(\ol X,\qp(1))^{G_L}.}\]
 One can easily deduce the commutativity of \eqref{eq6-2-1} from these facts. %The details are straight-forward and left to the reader.
\end{pf}
\par
By this lemma, it is enough to show that there is an element $\xi\in K_1(\XZp)^{(2)}$ such that $\varrho(\xi)\neq 0$ in $H^1_f(\qp,V^\ind)$. Let $\DZp_1:=\pi^{-1}(1)$ and $\DZp_2:=\pi^{-1}(-1)$ be the multiplicative fiber over $\zp$ which are N\'eron $1$-gon ({\bf Fact 2}). Put $\DZp:=\DZp_1+\DZp_2$ and $\UZp:=\XZp-\DZp$. We consider rational functions
\[ f_1:=\frac{Y-(X+4)}{Y+(X+4)}\Big|_{\DZp_1}, \qquad f_2:=\frac{Y-(X+4)}{Y+(X+4)}\Big|_{\DZp_2} \]
on $\DZp_1$ and $\DZp_2$ respectively. They define elements $\xi'_i\in K'_1(\DZp_i)^{(1)}$ by Quillen's localization exact sequence
\[ 0\lra K'_1(\DZp_i)^{(1)}\lra K'_1(\DZp_i^{\mathrm{reg}})^{(1)} \os{j_*\mathrm{div}}{\lra}\Q. \]
Here $\DZp_i^{\mathrm{reg}}$ denotes the smooth locus of $\DZp_i$, which is isomorphic to $\G_{m,\zp}$, and $j : \tDZp_i\to \DZp_i$ denotes the normalization. We denote by $\xi_i\in K_1(\XZp)^{(2)}$ the image of $\xi'_i$ via the natural map $K'_1(\DZp_i)^{(1)}\to K_1(\XZp)^{(2)}$. The goal is to prove
\begin{equation}\label{eqstep2}
\varrho(\xi_1-\xi_2)\ne 0 \; \text{ in } \; H^1_f(\qp,V^\ind).
\end{equation}
\addtocounter{thm}{2}
\begin{prop}\label{claimstep2}
If $p$ satisfies {\bf C(p)}, then we have
$$[D_1]-[D_2]\not\in \partial_\et(H^2(\ol U,\qp(2))^{G_{\qp}})$$ where $\partial_\et$ is the boundary map in \eqref{mc0}.
Hence we have 
\stepcounter{equation}
\begin{equation}\label{eq2step2}
\varrho(\xi_1-\xi_2)\neq 0 \text{ in }H^1_f(\qp,V)
\end{equation}
by the commutative diagram \eqref{diagram2}
where $V:=H^2(\ol X,\qp(2))/\langle D_1\rangle\ot\qp(1)$.
\end{prop}
We note that one can further show $H^2(\ol U,\qp(2))^{G_{\qp}}=0$ though we do not need it.
\par\medskip
\begin{pf}
It follows from {\bf C(p)-1} and {\bf Fact 2} that the conditions {\bf (A)'} and {\bf (B)'} in Proposition \ref{key4} and hence {\bf (A)} and {\bf (B)} in Theorem \ref{key3} are satisfied. We apply \eqref{conj1seq4}. It is enough to show 
\[ [D_1]-[D_2]\not\in \ol{\partial_\dR}\big(\Eis^{(p^2)}(\XR,\DR)_\zp\big)\] under {\bf C(p)}. The space $\vg(\XZp,\Omega^2_{\XZp/\zp}(\log \DZp))$ is a free $\zp$-module of rank 3 generated by
\[ \omega_0:=dt \frac{dX}{Y},\quad \omega_1:=\frac{dt}{t-1}\frac{dX}{Y},\quad \omega_2:=\frac{dt}{t+1}\frac{dX}{Y}, \]
which satisfies
\[ \partial_\dR(\omega_0)=0,\quad \partial_\dR(\omega_1)=\DZp_1,\quad \partial_\dR(\omega_2)=\DZp_2. \]
Let $\XZp^*$ be the `tubular neighborhood' of $\DZp_1$:
\[ \xymatrix{\XZp^* \ar[r] \ar[d] \ar@{}[rd]|{\square} & \XZp-\DZp_1 \ar[d] \\ \Spec(\bZ_p\psftt) \ar[r] & \bP^1_{\bZ_p}-\{1\}. } \]
We fix an isomorphism between $\XZp^*$ and the Tate elliptic curve
\[ E_q:y^2+xy=x^3+a_4(q)x+a_6(q) \;\; \hbox{ over \; $\zp\psfqq$} \]
given by
\[ x+\frac{1}{12}=-\frac{\,E_1^2 \cdot (X+3)\,}{12},\qquad y+\frac{\,1\,}{2}x=\frac{\,E_1^3\cdot Y\,}{24}, \qquad t^4=\frac{E_{3,a}}{(E_1)^3}=\frac{E_{3,a}}{E_{3,a}+27E_{3,b}}, \]
\[ \left(\quad\Longrightarrow \qquad j(q)=\frac{27(9-8t^4)^3}{(1-t^4)t^{12}}, \qquad j(q^3)=\frac{27(1+8t^4)^3}{(1-t^4)^3t^4} \quad\right)\]
where $E_1$, $E_{3,a}$ and $E_{3,b}$ are as in {\bf C(p)}.
% Note that they come from the $q$-expansions of the Eisenstein series for $\Gamma_1(3)$ at the cusp $\tau=\infty$.
Let $\iota$ be the composite map
\[ \iota:\vg(\XZp, \Omega^2_{\XZp}(\log \DZp))\to \vg(\XZp^*,\Omega^2_{\XZp^*})\simeq \vg(E_q,\Omega^2_{E_q})\to\zp\psfqq \frac{du}{u}\frac{dq}{q} \]
where $u$ denotes the Tate parameter of $E_q$. We note
\[ \iota\left(\frac{dt}{t}\frac{dX}{Y} \right)=-\frac{27}{4}\,E_{3,b} \frac{du}{u}\frac{dq}{q}, \qquad \iota\left(\frac{d(t^4-1)}{t^4-1}\frac{dX}{Y}\right) =E_{3,a}\frac{du}{u}\frac{dq}{q} \]
and
\begin{align*} & \phantom{\Longleftrightarrow \quad} \iota(\omega_1)=f_1(q)\frac{du}{u}\frac{dq}{q},\quad
\iota(\omega_2)=f_2(q)\frac{du}{u}\frac{dq}{q},\quad
\iota(\omega_0)=g(q)\frac{du}{u}\frac{dq}{q}, \end{align*}
where $f_i(q)$ and $g(q)$ are as in \eqref{fig}. Assume 
$[D_1]-[D_2]\in \ol{\partial_\dR}\big(\Eis^{(p^2)}(\XZp,\DZp)_\zp\big)$,
namely there is an $n\in \Z_p$ such that $\omega_1-\omega_2+n\omega_0\in \Eis^{(p^2)}(\XZp,\DZp)_\zp$. This means that
\begin{equation}\label{cr1}
 a_i- b_i+n c_i\equiv 0 \mod i^2\zp
\end{equation}
for all $i\leq p^2$ with $p|i$, which contradicts to {\bf C(p)-2}. 
\end{pf}

It remains to check \eqref{orthogonal} for $\xi_1-\xi_2$.
If $Z$ is a curve in (i) or (ii) (see Fact 3 in \S\ref{sect6-1}), one clearly has $v_Z(\varrho(\xi_i))=1$.
If $Z=E$ is the section of infinity,
then one has
\[ \frac{Y-(X+4)}{Y+(X+4)}\Big|_{D_i\cap E}=1 \]
and hence $v_Z(\varrho(\xi_i))=1$.
If $Z=C$ is the section in (iv), then one has
\[ \frac{Y-(X+4)}{Y+(X+4)}\Big|_{D_i\cap C}=
\frac{8/9-(-4/3+4)}{8/9+(-4/3+4)}=-\frac{1}{2}. \]
Hence $v_Z(\varrho(\xi_1))=v_Z(\varrho(\xi_2))=-1/2$ and $v_Z(\varrho(\xi_1-\xi_2))=1$.
This completes the proof of \eqref{eqstep2} and hence {\bf Step 2}.

\subsection{Finiteness of torsion in $\CH_0(\XQp)$}\label{sect6-3}
We end this paper by showing that the torsion part of $\CH_0(\XQp)$ is finite.
Since we have proved the finiteness of the $p$-primary torsion part,
it remains to show that the $\ell$-primary torsion part is finite for any $\ell\neq p$
 and zero for almost all $\ell\neq p$.
In view of the isomorphism \eqref{eq6-3-0} below, Bloch's exact sequence (cf.\ \cite{CTR1} (2.1))
\[ 0 \lra H^1_\zar(X,\cK_2) \ot \ql/\zl \lra N^1\hspace{-1.5pt}H^3(\XQp,\ql/\zl(2)) \lra \CH_0(\XQp)\{ \ell \} \lra 0 \]
%($\cK_2$ denotes the Zariski sheafification of $K_2$, and $N^\bullet$ denotes the coniveau filtration)
and the isomorphism $H^1_\zar(X,\cK_2) \otimes \bQ \simeq K_1(\XQp)^{(2)}$ (cf.\ \cite{So}), 
it is enough to show that the regulator map % the $\ell$-adic regulator map
\begin{equation}\label{eq6-3-1}
 K_1(\XQp)^{(2)}\ot\ql \lra H^1(\qp,H^2(\ol X,\ql(2)))
\end{equation}
is surjective for any $\ell\neq p$ and that $H^3(\XQp,\qzl(2))$ is divisible for almost all $\ell\neq p$.

Assume $\ell \ne p$ in what follows. Put $\ol Y:=Y \otimes_{\F_p} {\ol {\F_p}}$. There are isomorphisms
\[ H^1(\qp,H^2(\ol X,\ql(2)))\simeq H^0(\F_p,H^2(\ol Y,\ql(1)))\simeq \NS(Y)_{\ql} \]
by the Tate conjecture for $Y$ (\cite{ArSw}) and a similar argument as for \eqref{eq6-3-2} below.
The map \eqref{eq6-3-1} is identified with the boundary map
\begin{equation}\label{eq6-3-3}
 K_1(\XQp)^{(2)}\ot\ql \lra \NS(Y)_{\ql}.
\end{equation}
This map is surjective by {\bf Step 1}, which shows that \eqref{eq6-3-1} is surjective.

Next we show that $H^3(\XQp,\qzl(2))$ is divisible for almost all $\ell$.
Since $\ol X$ is a $K$3 surface, we have $H^1(\ol X,\zl)=H^3(\ol X,\zl)=0$ and $H^2(\ol X,\zl)$ is torsion-free for any $\ell$.
Hence we have \[ H^2(\ol X,\qzl)=H^2(\ol X,\zl)\ot\qzl \quad \hbox{ and } \quad
 H^1(\ol X,\qzl)=H^3(\ol X,\qzl)=0, \] and moreover
\begin{equation}\label{eq6-3-0}
 H^3(\XQp,\qzl(2))=H^1(\qp, H^2(\ol X,\qzl(2)))
\end{equation}
by a Hochschild-Serre spectral sequence. We compute the right hand side as follows. There is a short exact sequence
\begin{align*}
0\to H^1(\F_p,H^2(\ol X,\qzl(2)))
 &\to H^1(\qp,H^2(\ol X,\qzl(2))) \\
 &\to H^2(\ol Y,\qzl(1))^{G_{\F_p}}\to 0.
\end{align*}
By the divisibility of $H^2(\ol X,\qzl)$ and a standard argument on weights (cf.\ \cite{CTSS} \S2, \cite{De}), the first term is zero. Thus we have
\begin{equation}\label{eq6-3-2}
H^1(\qp,H^2(\ol X,\qzl(2))) \simeq H^2(\ol Y,\qzl(1))^{G_{\F_p}},
\end{equation}
and we are reduced to showing that the right hand side is divisible for almost all $\ell$.
Put \[ N_\ell:=H^2(\ol Y, \zl(1))^{G_{\F_p}}, \] and note the following fact due to Deligne \cite{De}:
\begin{enumerate}
\item[($*$)] {\it the characteristic polynomial of the geometric Frobenius $\varphi$
acting on $H^2(\ol Y, \ql(1))$ is independent of $\ell (\ne p)$}.
\end{enumerate}
Since $H^2(\ol Y,\zl(1))$ is torsion-free,
 it is easy to see that $H^2(\ol Y, \zl(1))/N_\ell$ is torsion-free as well, for any $\ell$.
% by the exact sequence
% \[\xymatrix{ 0 \ar[r] & N_\ell \ar[r] & H^2(\ol Y, \zl(1)) \ar[r]^-{1-\varphi} & H^2(\ol Y, \zl(1)). }\]
Hence there is a short exact sequence
\[ 0\lra N_\ell \ot\qzl\lra 
H^2(\ol Y, \qzl(1))  \lra 
(H^2(\ol Y, \zl(1))/N_\ell)\ot\qzl\lra 0 \]
for any $\ell$. By ($*$), we have
\[ ((H^2(\ol Y, \zl(1))/N_\ell)\ot\qzl)^{G_{\F_p}} =0 \]
for almost all $\ell$. For such $\ell$, we have
\[ H^2(\ol Y, \qzl(1))^{G_{\F_p}}=N_\ell\ot\qzl, \]
which is divisible. This completes the proof of the finiteness of $\CH_0(\XQp)_\tors$.
\addtocounter{thm}{4}
\begin{rem}
Here is a more systematic {\rm(}but essentially the same{\rm)} proof of the finiteness result in this subsection.
By the surjectivity of \eqref{eq6-3-3} and a result of Spiess {\rm\cite{Sp}} Proposition {\rm 4.3}
 {\rm(}see also {\rm\cite{SS2}} for a generalization{\rm)}, we have
\[ \CH_0(\XQp)_\tors \simeq \CH_0(Y)_\tors \quad \hbox{up to the  $p$-primary torsion part.} \]
The right hand side is finite by Colliot-Th\'el\`ene--Sansuc--Soul\'e {\rm\cite{CTSS}}.
\end{rem}

\noindent
Department of Mathematics, Hokkaido University,
Sapporo 060-0810,
JAPAN

\smallskip

\noindent
{\it E-mail} : \textbf{asakura@math.sci.hokudai.ac.jp}

\bigskip

\noindent
Graduate school of Mathematics, Nagoya University,
Nagoya 464-8602,
JAPAN

\smallskip

\noindent
{\it Email} : \textbf{kanetomo@math.nagoya-u.ac.jp}
\end{document}